\documentclass[a4paper]{article}
\usepackage{amsmath,amssymb,amsfonts,bm,natbib}
\usepackage{array,graphicx}
\pdfoutput=1

\newlength{\Figsize}
\newcommand{\ci}{\mathrm{i}}
\newcommand{\ce}{\mathrm{e}}
\newcommand{\sign}{\operatorname{sign}}

\begin{document}

\title{Hopf bifurcation with zero frequency and imperfect $SO(2)$ symmetry}
\author{F. Marques$^*$
 \and A. Meseguer\thanks{Departament de F{\'\i}sica Aplicada, Universitat
  Polit\`ecnica de Catalunya, 08034 Barcelona, Spain}
 \and Juan M. Lopez\thanks{School of Mathematical and Statistical
   Sciences, Arizona State University, Tempe AZ 85287, USA}
 \and J. R. Pacheco$^{\dagger,}$\thanks{Environmental Fluid Dynamics
   Laboratories, Department of Civil Engineering and Geological
   Sciences, University of Notre Dame, Notre Dame, Indiana  46556, USA}
 \and Jose M. Lopez$^*$
}

\maketitle

\begin{abstract} 

Rotating waves are periodic solutions in $SO(2)$ equivariant dynamical
systems. Their precession frequency changes with parameters and it may
change sign, passing through zero. When this happens, the dynamical
system is very sensitive to imperfections that break the $SO(2)$
symmetry and the waves may become trapped by the imperfections,
resulting in steady solutions that exist in a finite region in
parameter space. This is the so-called pinning phenomenon. In this
study, we analyze the breaking of the $SO(2)$ symmetry in a dynamical
system close to a Hopf bifurcation whose frequency changes sign along
a curve in parameter space. The problem is very complex, as it
involves the complete unfolding of high codimension. A detailed
analysis of different types of imperfections indicates that a pinning
region surrounded by infinite-period bifurcation curves appears in all
cases. Complex bifurcational processes, strongly dependent on the
specifics of the symmetry breaking, appear very close to the
intersection of the Hopf bifurcation and the pinning region. Scaling
laws of the pinning region width, and partial breaking of $SO(2)$ to
$Z_m$, are also considered. Previous and new experimental and
numerical studies of pinned rotating waves are reviewed in light of
the new theoretical results.

\end{abstract}

\section{Introduction}

Dynamical systems theory plays an important role in many areas of
mathematics and physics because it provides the building blocks that
allow us to understand the changes many physical systems experience in
their dynamics when parameters are varied. These building blocks are
the generic bifurcations (saddle-node, Hopf, etc.) that any arbitrary
physical system experiences under parameter variation, regardless of
the physical mechanisms underlying the dynamics. When one single
parameter of the system under consideration is varied, codimension-one
bifurcations are expected. If the system depends on more parameters,
higher codimension bifurcations appear and they act as organizing
centers of the dynamics.

The presence of symmetries changes the nature and type of bifurcations
that a dynamical system may undergo. Symmetries play an important role
in many idealized situations, where simplifying assumptions and the
consideration of simple geometries result in dynamical systems
equivariant under a certain symmetry group. Bifurcations with symmetry
have been widely studied \citep{GoSc85,GSS88,ChIo94,GoSt02,ChLa00,CrKn91}. 
However, in any real system, the symmetries are only approximately
fulfilled, and the breaking of the symmetries, due to the presence of
noise, imperfections and/or other phenomena, is always present.  There
are numerous studies of how imperfect symmetries lead to dynamics that
are unexpected in the symmetric problem, e.g.\ 
\citep{Kee87,CaHo92,KHD95,HiKn96,DHK97,LaWu00}. However, a complete
theory is currently unavailable.

One observed consequence of imperfections in systems that support
propagating waves is that the waves may become trapped by the
imperfections \citep[e.g., see][]{Kee87,WeBu03,ThKn06_prl,ThKn06_njp}.
In these various examples, the propagation direction is typically
biased. However, a more recent problem has considered a case where a
rotating wave whose sense of precession changes sign is pinned by
symmetry-breaking imperfections \citep{AHHP08}. We are unaware of any
systematic analysis of the associated normal form dynamics for such a
problem and this motivates the present study.

When a system is invariant to rotations about an axis (invariance
under the $SO(2)$ symmetry group), $SO(2)$-symmetry-breaking Hopf
bifurcations result in rotating waves, consisting of a pattern that
rotates about the symmetry axis at a given precession frequency
without changing shape.  This frequency is parameter dependent, and in
many problems, when parameters are varied, the precession frequency
changes sign along a curve in parameter space. What has been observed
in different systems is that in the presence of imperfections, the
curve of zero frequency becomes a band of finite width in parameter
space. Within this band, the rotating wave becomes a steady
solution. This is the so-called pinning phenomenon. It can be
understood as the attachment of the rotating pattern to some
stationary imperfection of the system, so that the pattern becomes
steady, as long as its frequency is small enough so that the
imperfection is able to stop the rotation. This pinning phenomenon
bears some resemblance to the frequency locking phenomena, although in
the frequency locking case we are dealing with a system with two
non-zero frequencies and their ratio becomes constant in a region of
parameter space (a resonance horn), whereas here we are dealing with a
single frequency crossing zero.

In the present paper, we analyze the breaking of $SO(2)$ symmetry in a
dynamical system close to a Hopf bifurcation whose frequency changes
sign along a curve in parameter space. The analysis shows that
breaking $SO(2)$ symmetry is much more complex than expected,
resulting in a bifurcation of high codimension (about nine). Although
it is not possible to analyze in detail such a complex and
high-codimension bifurcation, we present here the analysis of five
different ways to break $SO(2)$ symmetry. This is done by introducing
into the normal form all the possible terms, up to and including
second order, that break the symmetry, and analyzing each of these
five terms separately. Three of these particular cases have already
been analyzed in completely different contexts unrelated to the
pinning phenomenon \citep{Gam85,Wag01, BDV08, SaWa10}. In the present
study, we extract the common features that are associated with the
pinning. In all cases, we find that a band of pinning solutions
appears around the zero frequency curve from the symmetric case, and
that the band is delimited by curves of infinite-period
bifurcations. The details of what happens when the infinite-period
bifurcation curves approach the Hopf bifurcation curve are different
in the five cases, and involve complicated dynamics with several
codimension-two bifurcations occurring in a small region of parameter
space as well as several global bifurcations.

Interest in the present analysis is two-fold. First of all, although
the details of the bifurcational process close to the zero-frequency
Hopf point are very complicated and differ from case to case, for all
cases we observe the appearance of a pinning band delimited by
infinite-period bifurcations of homoclinic type that, away from the
small region of complicated dynamics, are SNIC bifurcations
\citep[saddle-node on an invariant circle bifurcation,
  e.g.\ see][]{Str94}. Secondly, some of the scenarios analyzed are
important \emph{per se}, as they correspond to the generic analysis of
a partial breaking of the $SO(2)$ symmetry, so that after the
introduction of perturbations, the system still retains a discrete
symmetry (the $Z_2$ case is analyzed in detail).

The paper is organized as follows. In section \S\ref{Sec_perfect_Hopf}
the properties of a Hopf bifurcation with $SO(2)$ symmetry with the
precession frequency crossing through zero are summarized, and the
general unfolding of the $SO(2)$ symmetry breaking process is
discussed. The next sections explore the
particulars of breaking the symmetry at order zero
(\S\ref{Sec_epsilon}), one (\S\ref{Sec_bar_z}) and two
(\S\ref{Sec_quadratic}). Sections
\S\ref{Sec_bar_z} and \S\ref{Sec_Z3} are
particularly interesting because they consider the symmetry-breaking
processes $SO(2)\to Z_2$ and $SO(2)\to Z_3$ which are readily realized
experimentally. Section \S\ref{Sec_general} extracts the general
features of the pinning problem from the analysis of the specific
cases carried out in the earlier sections. Section
\S\ref{Sec_experiments} presents comparisons with experiments and
numerical computations in two real problems in fluid dynamics,
illustrating the application of the general theory developed in the
present study. Finally, in \S\ref{conclusions}, conclusions and
perspectives are presented.

\section{Hopf bifurcation with $SO(2)$ symmetry and zero frequency}
\label{Sec_perfect_Hopf}

The normal form for a Hopf bifurcation is
\begin{equation}\label{HopfNF}
 \dot z=z(\mu+i\omega-c|z|^2),
\end{equation}
where $z$ is the complex amplitude of the bifurcating periodic
solution, $\mu$ is the bifurcation parameter, and $\omega$ and $c$ are
functions of $\mu$ and generically at the bifurcation point
($\mu=0$) both are different from zero. It is the non-zero character
of $\omega$ that allows one to eliminate the quadratic terms in $z$ in
the normal form. This is because the normal form $\dot z=P(z,\bar z)$
satisfies \citep[e.g., see][]{HaIo11}
\begin{equation}\label{NFevalcond}
 P(\ce^{-\ci\omega t}z,\ce^{\ci\omega t}\bar z)=\ce^{-\ci\omega t}P(z,\bar z),
\end{equation}
where $P$ is a low order polynomial that captures the dynamics in a
neighborhood of the bifurcation point.  If $\omega=0$, this equation
becomes an identity and $P$ cannot be simplified.  The case $\omega=0$
is a complicated bifurcation and it depends on the details of the
double-zero eigenvalue of the linear part $L$ of $P$; as $z=x+\ci y$
is complex, the matrix of $L$ using the real coordinates $(x,y)$ is a
real $2\times2$ matrix. If $L$ is not completely degenerate, that is
\begin{equation}
 L=\begin{pmatrix} 0 & 1 \\ 0 & 0 \end{pmatrix},
\end{equation}
then we have the well-studied Takens--Bogdanov bifurcation, whereas
the completely degenerate case,
\begin{equation}\label{CompletDegen}
 L=\begin{pmatrix} 0 & 0 \\ 0 & 0 \end{pmatrix},
\end{equation}
is a high-codimension bifurcation that has not been completely
analyzed.

\begin{figure}
  \begin{center}
    \includegraphics[width=0.8\linewidth]{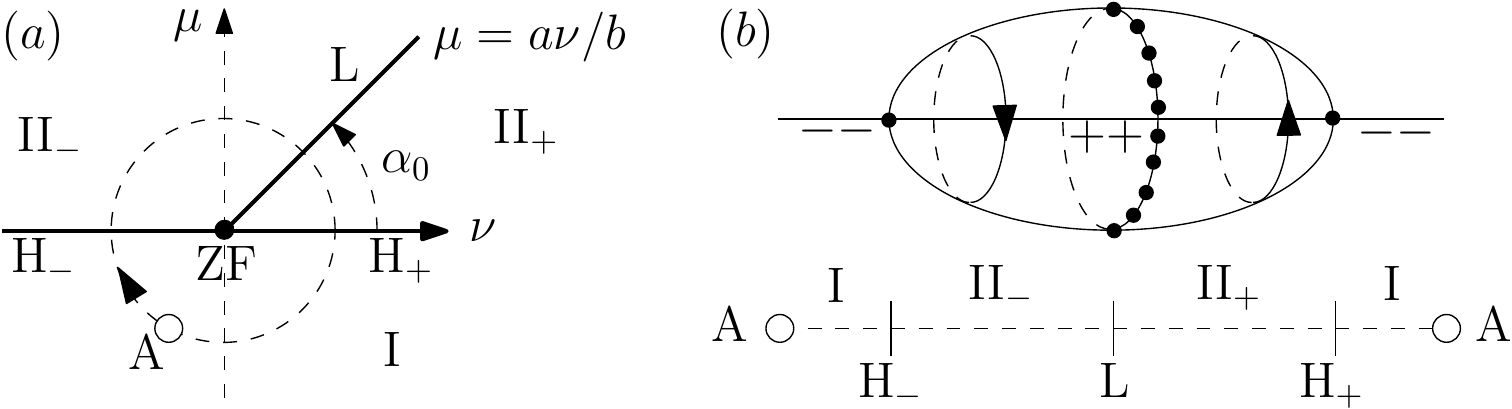}
  \end{center}
  \caption{Hopf bifurcation with $SO(2)$ symmetry and zero frequency;
    part $(a)$ shows the bifurcation diagram, where the thick lines are
    bifurcation curves, and part $(b)$ shows the bifurcations along the
    path A shown in $(a)$. The fixed point curve is labeled with the
    signs of its eigenvalues. In regions II$_-$ and II$_+$ the limit
    cycles, born at the Hopf bifurcations H$_-$ and H$_+$, rotate in
    opposite senses. L is the line where the limit cycle becomes an
    invariant curve of fixed points.}
  \label{zeroHopf}
\end{figure}

If the system has $SO(2)$ symmetry, it must also satisfy
\begin{equation}\label{NFsymSO2}
 P(\ce^{\ci m\theta}z,\ce^{-\ci m\theta}\bar z)=\ce^{\ci m\theta}P(z,\bar z),
\end{equation}
where $Z_m$ is the discrete symmetry retained by the bifurcated
solution. When the group $Z_m$ is generated by rotations of angle
$2\pi/m$ about an axis of $m$-fold symmetry, as is usually the case
with $SO(2)$, then the group is also called $C_m$. Equations
\eqref{NFevalcond} and \eqref{NFsymSO2} are completely equivalent and
have the same implications for the normal form structure. Advancing in
time is the same as rotating the solution by a certain angle ($\omega
t=m\theta$); the bifurcated solution is a rotating wave. Therefore, if
$\omega$ becomes zero by varying a second parameter, we still have the
same normal form \eqref{HopfNF}, due to \eqref{NFsymSO2}, with
$\omega$ replaced by a small parameter $\nu$:
\begin{equation}\label{zeroHopfNF}
 \dot z=z(\mu+i\nu-c|z|^2).
\end{equation}
The Hopf bifurcation with $SO(2)$ symmetry and zero
frequency is, in this sense, trivial. Introducing the modulus and phase
of the complex amplitude $z=r\ce^{\ci\phi}$, the normal form becomes
\begin{equation}\label{zeroHopfNF2}
\begin{aligned}
 & \dot r=r(\mu-ar^2),\\
 & \dot\phi=\nu-br^2,
\end{aligned}
\end{equation}
where $c=a+\ci b$, and let us assume for the moment that $a$ and $b$
are positive. The bifurcation frequency in \eqref{zeroHopfNF2} is now
the small parameter $\nu$. The bifurcated solution $RW_m$ exists only
for $\mu>0$, and has amplitude $r=\sqrt{\mu/a}$ and frequency
$\omega=\nu-b\mu/a$. The limit cycle $RW_m$ becomes an invariant set
of steady solutions along the straight line $\mu=a\nu/b$ (labeled L in
figure~\ref{zeroHopf}) where the frequency of $RW_m$ goes to zero; the
angle between L and the Hopf bifurcation curve (the horizontal axis
$\mu=0$) is $\alpha_0$. The bifurcation diagram and a schematic of the
bifurcations along a one-dimensional path is also shown in
figure~\ref{zeroHopf}. The bifurcation point $\mu=\nu=0$, labeled ZF
(zero-frequency Hopf point) in figure~\ref{zeroHopf}$(a)$, is a
codimension-two bifurcation. It coincides with the generic Hopf
bifurcation, except that it includes a line L along which the
bifurcated solution has zero frequency.

Assuming $c\ne0$, we can simplify \eqref{zeroHopfNF2} by scaling $z$ so
that $|c|=1$; we will write
\begin{equation}\label{tiltangle}
  c=a+\ci b=\ci\ce^{-\ci\alpha_0}=\sin\alpha_0+\ci\cos\alpha_0,\quad
  b+\ci a=\ce^{\ci\alpha_0},
\end{equation}
which helps simplify subsequent expressions. The case $a$ and $b$ both
positive, which we will analyze in detail in the following sections,
corresponds to one of the fluid dynamics problems that motivated the
present analysis \citep[see][and \S\ref{Sec_experiments}\ref{TCpinning}]
{AHHP08,PLM11}. For other signs of $a$ and $b$, analogous conclusions
can be drawn. It is of particular interest to consider the subcritical
case $a<0$ as it corresponds to the other fluid dynamics problem
analyzed here \citep[see][and \S\ref{Sec_experiments}\ref{conrot_pinning}]
{MMBL07,LoMa09}. By reversing time and changing the sign of $\mu$ and
$\nu$, we obtain exactly the same normal form \eqref{zeroHopfNF2} but
with the opposite sign of $a$ and $b$. By changing the sign of $\phi$
and $\nu$, we obtain \eqref{zeroHopfNF2} with the opposite sign of
$b$. Therefore, all possible cases corresponding to different signs of
$a$ and $b$ can be reduced to the case where $a$ and $b$ are both
positive.

\subsection{Unfolding the Hopf bifurcation with zero frequency}

If the $SO(2)$ symmetry in the normal form \eqref{zeroHopfNF} is
completely broken, and no symmetry remains, then the restrictions
imposed on the normal form by \eqref{NFsymSO2} disappear completely
and all the terms in $z$ and $\bar z$ missing from \eqref{zeroHopfNF}
will reappear multiplied by small parameters. This means that the
normal form will be
\begin{equation}\label{NFcompleteunfold}
 \dot z=z(\mu+i\nu-c|z|^2)+\epsilon_1+\epsilon_2\bar z+\epsilon_3\bar{z}^2+
      \epsilon_4z\bar z+\epsilon_5 z^2,
\end{equation}
where additional cubic terms have been neglected because we assume
$c\ne0$ and that $cz|z|^2$ will be dominant. As the $\epsilon_i$ are
complex, we have a problem with 12 parameters. Additional
simplifications can be made in order to obtain the so-called
hypernormal form; this method is extensively used by \citet{Kuz04},
for example. Unfortunately, many of the simplifications rely on having
some low-order term in the normal form being non-zero with a
coefficient of order one. For example, if $\omega\ne0$, it is possible
to make $c$ real by using a time re-parametrization. In our problem,
all terms up to and including second order are zero or have a small
coefficient, and so only a few simplifications are possible. These
simplifications are an infinitesimal translation of $z$ (two
parameters), and an arbitrary shift in the phase of $z$ (one
parameter). Using these transformations the twelve parameters can be
reduced to nine. In particular, one of either $\epsilon_4$ or
$\epsilon_5$ can be taken as zero and the other can be
made real. By rescaling $z$, we can make $c$ of modulus one, as in
\eqref{tiltangle}. A complete analysis of a normal form depending on
nine parameters, i.e.\ a bifurcation of codimension of about nine, is
completely beyond the scope of the present paper. In the literature,
only codimension-one bifurcations have been completely analyzed. Most
of the codimension-two bifurcations for ODE and maps have also been
analyzed, except for a few bifurcations for maps that remain
outstanding \citep{Kuz04}. A few codimension-three and very few
codimension-four bifurcations have also been analyzed
\citep{CLW94,drs97}, but to our knowledge, there is no systematic
analysis of bifurcations of codimension greater than two.

In the following sections, we consider the five cases, $\epsilon_1$ to
$\epsilon_5$, separately. A combination of analytical and numerical
tools allows for a detailed analysis of these bifurcations. We extract
the common features of the different cases when
$\epsilon_i\ll\sqrt{\mu^2+\nu^2}$, which captures the relevant
behavior associated with weakly breaking $SO(2)$ symmetry. In
particular, the $\epsilon_2$ case exhibits very interesting and rich
dynamics that may be present in some practical cases when the $SO(2)$
symmetry group is not completely broken and a $Z_2$ symmetry group,
generated by the half-turn $\theta \to \theta + \pi$, remains.

Some general comments can be made here about these five cases, which
are of the form
\begin{equation}\label{NFsinglepert}
 \dot z=z(\mu+i\nu-c|z|^2)+\epsilon z^q\bar{z}^{p-q},
\end{equation}
for integers $0\le q\le p\le 2$, excluding the case $p=q=1$ which is
$SO(2)$ equivariant and so $\epsilon$ can be absorbed into $\mu$ and
$\nu$. By changing the origin of the phase of $z$, we can modify the
phase of $\epsilon$ so that it becomes real and positive. Then, by
re-scaling $z$, time $t$, and the parameters $\mu$ and $\nu$ as
\begin{equation}\label{epsilon_scaling}
 (z,t,\mu,\nu)\to(\epsilon^{\delta}z,\epsilon^{-2\delta}t,
  \epsilon^{2\delta}\mu,\epsilon^{2\delta}\nu),\quad\delta=\frac{1}{3-p},
\end{equation}
we obtain \eqref{NFsinglepert} with $\epsilon=1$, effectively leading
to codimension-two bifurcations in each of the five cases. We expect
complex behavior for $\mu^2+\nu^2\lesssim\epsilon^2$, when the three
parameters are of comparable size, while the effects of small
imperfections breaking $SO(2)$ will correspond to
$\mu^2+\nu^2\gg\epsilon^2$. From now on $\epsilon=1$ will be assumed,
and we can restore the explicit $\epsilon$-dependence by reversing the
transformation \eqref{epsilon_scaling}. Three of the five normal forms
\eqref{epsilon_scaling} have been analyzed in the literature
(discussed below), focusing on the regions where $\mu$, $\nu$ and
$\epsilon$ are of comparable size; here we will also consider what
happens for $\mu^2+\nu^2\gg\epsilon^2$ which is particularly important
for the pinning phenomenon.

The normal forms corresponding to the $\epsilon_1$, $\epsilon_2$ and
$\epsilon_3$ cases have already been analyzed in contexts completely
different to the $SO(2)$ symmetry-breaking context considered
here. The context in which these problems were studied stems from
low-order resonances in perturbed Hopf problems. \citet{Gam85} studied
time-periodic forcing near a Hopf bifurcation point, analyzing the
problem using the Poincar\'e stroboscopic map. The normal forms
corresponding to the 1:1, 1:2 and 1:3 strong resonances coincide with
the normal forms we present below for cases with only the
$\epsilon_1$, $\epsilon_2$ and $\epsilon_3$ terms retained in
\eqref{NFcompleteunfold}, respectively.  Later, motivated by a problem
of a nonlinear oscillator with damping and quasi-periodic driving, a
series of papers extended the strong resonances results of
\citet{Gam85} by studying the semi-global bifurcations for
periodically and quasi-periodically perturbed driven damped
oscillators near a Hopf bifurcation \citep[see][and references
  therein]{Wag01, BDV08, SaWa10}.  The other two cases we consider,
with only the $\epsilon_4$ or the $\epsilon_5$ terms retained in
\eqref{NFcompleteunfold}, do not appear to have been studied
previously. They fall outside of the context in which the other three
were studied because they do not correspond to any canonical resonance
problem. We should point out that within the resonance context, the
three cases studied would not make sense to consider in combination
(they correspond to completely distinct frequency ratios and so would
not generically occur in a single problem). In contrast, within the
context motivating our study, all five cases correspond to different
ways in which the $SO(2)$ symmetry of a system may be broken, and in a
physical realization, all five could co-exist.  In the following
sections, we present a detailed analysis of all five cases.

\section{Symmetry breaking of $SO(2)$ with an $\epsilon$ term}
\label{Sec_epsilon}

The normal form to be analyzed is \eqref{NFsinglepert} with $p=q=0$ and
$\epsilon=1$:
\begin{equation}\label{complexNFepsilon}
 \dot z=z(\mu+i\nu-c|z|^2)+1.
\end{equation}
This case has been analyzed in \citet{Gam85,Wag01,BDV08,SaWa10}.

\begin{figure}
  \begin{center}
    \begin{tabular}{c@{\hspace{1cm}}c}
      $(a)$ & $(b)$ \\
      \raisebox{8mm}{\includegraphics[width=0.3\linewidth]{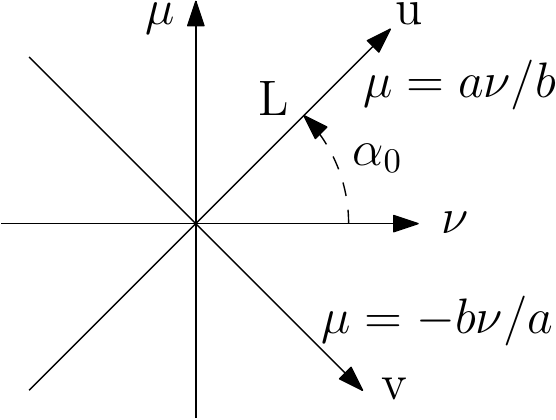}} &
      \includegraphics[width=0.36\linewidth]{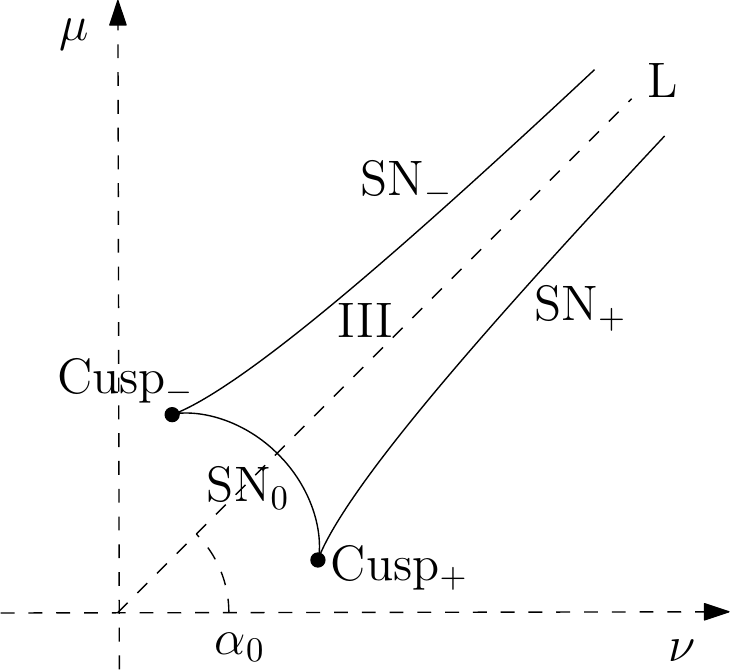}
    \end{tabular}
  \end{center}
  \caption{$(a)$ Coordinates $(u,v)$ in parameter space adapted to the
    line L, which coincides with zero frequency curve in the unperturbed
    $SO(2)$ symmetric case. $(b)$ Steady bifurcations of the fixed
    points corresponding to the normal form \eqref{complexNFepsilon}.
    SN$_\pm$ and SN$_0$ are saddle-node bifurcation curves and Cusp$_\pm$
    are cusp bifurcation points. In region III there exist three fixed
    points, and only one in the rest of parameter space.}
  \label{uv_and_epsilon_SN}
\end{figure}

It is convenient to introduce coordinates $(u,v)$ in parameter space,
rotated an angle $\alpha_0$ with respect to $(\mu,\nu)$, so that the
line L becomes the new coordinate axis $v=0$. Many features, that
are symmetric with respect to L will then simplify; e.g.\ the
distance to the bifurcation point along L is precisely $u$.
\begin{equation}\label{uv_munu}
  \begin{pmatrix} u \\ v \end{pmatrix}=\begin{pmatrix} a & b \\
  -b & a \end{pmatrix}\begin{pmatrix} \mu \\ \nu \end{pmatrix}=
  \begin{pmatrix} a\mu+b\nu \\ a\nu-b\mu \end{pmatrix},\quad
  \begin{pmatrix} \mu \\ \nu \end{pmatrix}=\begin{pmatrix} a & -b \\
  b & a \end{pmatrix}\begin{pmatrix} u \\ v \end{pmatrix},
\end{equation}
where $a=\sin\alpha_0$, $b=\cos\alpha_0$. Figure~\ref{uv_and_epsilon_SN}$(a)$
shows the relationship between the two coordinate systems.

\subsection{Fixed points and their local bifurcations}

The normal form \eqref{complexNFepsilon}, in terms of the modulus and
phase $z=r\ce^{\ci\phi}$, is
\begin{equation}\label{realNFepsilon}
\begin{aligned}
 & \dot r=r(\mu-ar^2)+\cos\phi,\\
 & \dot\phi=\nu-br^2-\frac{1}{r}\sin\phi.
\end{aligned}
\end{equation}
The fixed points are given by a cubic equation in $r^2$, and so we do
not have convenient closed forms for the corresponding roots (the
Tartaglia explicit solution is extremely involved). However, it is
easy to obtain the locus where two of the fixed points coalesce. The
parameter space is divided into two regions, region III with three
fixed points, and the rest of parameter space with one fixed point, as
seen in figure~\ref{uv_and_epsilon_SN}$(b)$. The curve separating both
regions is a saddle-node curve given by (see \ref{Appendix_epsilon})
\begin{equation}\label{cusp_curve}
  (u,v)=\frac{(3+3s^2,2\sqrt{3}s)}{(2+6s^2)^{2/3}},
  \quad s\in(-\infty,+\infty),
\end{equation}
in $(u,v)$ coordinates \eqref{uv_munu}. The saddle-noddle curve is
divided into three different arcs SN$_+$, SN$_-$ and SN$_0$ by two
codimension-two cusp bifurcation points, Cusp$_\pm$. SN$_-$
corresponds to values $s\in(-\infty,1)$, SN$_0$ to $s\in(-1,+1)$ and
SN$_+$ to $s\in(1,+\infty)$. The cusp points Cusp$_\pm$ have values
$s=\pm1$. The curves SN$_+$ and SN$_-$ are asymptotic to the line L,
and region II is the pinning region in this case. The two fixed points
that merge on the saddle-node curve have phase space coordinates
$z_2=r_2\ce^{\ci\phi_2}$, and the third fixed point is
$z_0=r_0\ce^{\ci\phi_0}$, where
\begin{equation}
  r_0=\Big(\frac{4}{1+3s^2}\Big)^{1/3},\quad
  r_2=\Big(\frac{1+3s^2}{4}\Big)^{1/6},
\end{equation}
and the phases are obtained from $\sin\phi=r(\nu-br^2)$,
$\cos\phi=r(ar^2-\mu)$.

\begin{figure}
  \begin{center}
    \begin{tabular}{m{0.32\linewidth}@{\hspace{0.15\linewidth}}
        m{0.38\linewidth}c}
      \centering $(a)$ & \centering $(b)$ & \\
      \includegraphics[width=\linewidth]{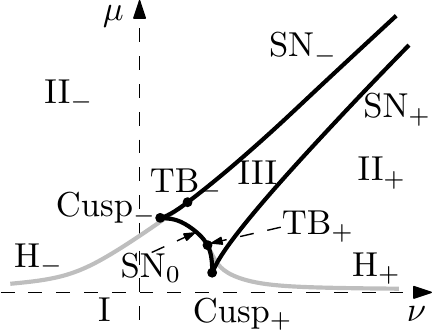} &
      \includegraphics[width=\linewidth]{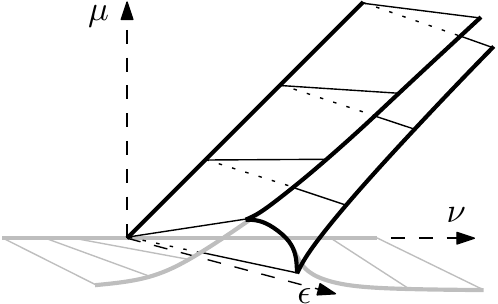} & \\
      \multicolumn{2}{c}{$(c)$} & \\
      \multicolumn{2}{c}{\includegraphics[width=0.44\linewidth]
        {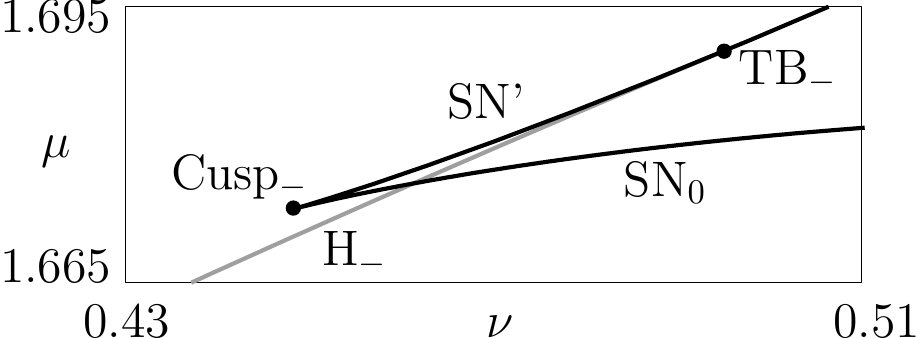}} & \\
    \end{tabular}\vspace*{-6pt}
  \end{center}
  \caption{$(a)$ Bifurcations of the fixed points corresponding to the
    normal form \eqref{complexNFepsilon}. SN$_\pm$ and SN$_0$ are
    saddle-node bifurcation curves, H$_\pm$ are Hopf bifurcation
    curves, Cusp$_\pm$ are cusp bifurcation points and TB$_\pm$ are
    Takens--Bogdanov bifurcation points. A perspective view of the
    corresponding codimension-three bifurcation in terms of
    $(\mu,\nu,\epsilon)$ is shown in $(b)$. $(c)$ is a zoom of $(a)$,
    showing that Cusp$_-$ and TB$_-$ are different. The parameter
    range in $(a)$ is $(\nu,\mu)\in[-3,6]\times[-1,6]$, for
    $\alpha=45^\text{o}$.}
  \label{epsilon_SN_Hopf}
\end{figure}

The Hopf bifurcations of the fixed points can be obtained by imposing
the conditions $T=0$ and $D>0$, where $T$ and $D$ are the trace and
determinant of the Jacobian of \eqref{complexNFepsilon}. These
conditions result in two curves of Hopf bifurcations (see
\ref{Appendix_epsilon} for details):
\begin{equation}
\begin{gathered}
  (\mu,\nu)=a^{1/3}(1-s^2)^{1/3}\bigg(2\,,
  \frac{b}{a}+\frac{s}{\sqrt{1-s^2}}\bigg), \\
  {\rm H}_-:\ s\in\Big(-1,-\sqrt{(1-b)/2}\Big),\qquad
  {\rm H}_+:\ s\in\Big(\sqrt{(1+b)/2},+1\Big).
\end{gathered}
\end{equation}
For $s\to\pm1$ both curves are asymptotic to the $\mu=0$ axis
($\nu\to\pm\infty$), the Hopf curve for $\epsilon=0$; the stable limit
cycles born at these curves are termed $C_-$ and $C_+$
respectively. The limit cycles (rotating waves) in II$_\pm$ rotate in
opposite directions, and III is the pinning region where the rotation
stops and we have a stable fixed point. Solutions with $\omega=0$,
that existed only along a single line in the absence of imperfections,
now exist in a region of finite width. Figure~\ref{epsilon_SN_Hopf}$(b)$ 
shows what happens when the $\epsilon$ dependence is restored; what we
have is that figure~\ref{epsilon_SN_Hopf}$(a)$ just scales with
$\epsilon$ as indicated in \eqref{epsilon_scaling}, and the pinning
region collapses onto the line L of the perfect case with $SO(2)$ symmetry.

The other ends of the H$_\pm$ curves are on the saddle-node curves of
fixed points previously obtained, and at these points $T=D=0$, so they
are Takens--Bogdanov points TB$_\pm$, as shown in
figure~\ref{epsilon_SN_Hopf}. The TB$_-$ and Cusp$_-$ codimension-two
bifurcation points are very close, as shown in the zoomed-in
figure~\ref{epsilon_SN_Hopf}$(b)$. In fact, depending on the angle
$\alpha_0$, the Hopf curve H$_-$ is tangent to, and ends at, either
SN$_-$ or SN$_0$. For $\alpha_0>60^\text{o}$, H$_-$ ends at SN$_0$,
and for $\alpha_0=60^\text{o}$ Cusp$_-$ and TB$_-$ coincide, and H$_-$
ends at the cusp point, a very degenerate case.

\begin{figure}
  \begin{center}
    \includegraphics[width=0.65\linewidth]{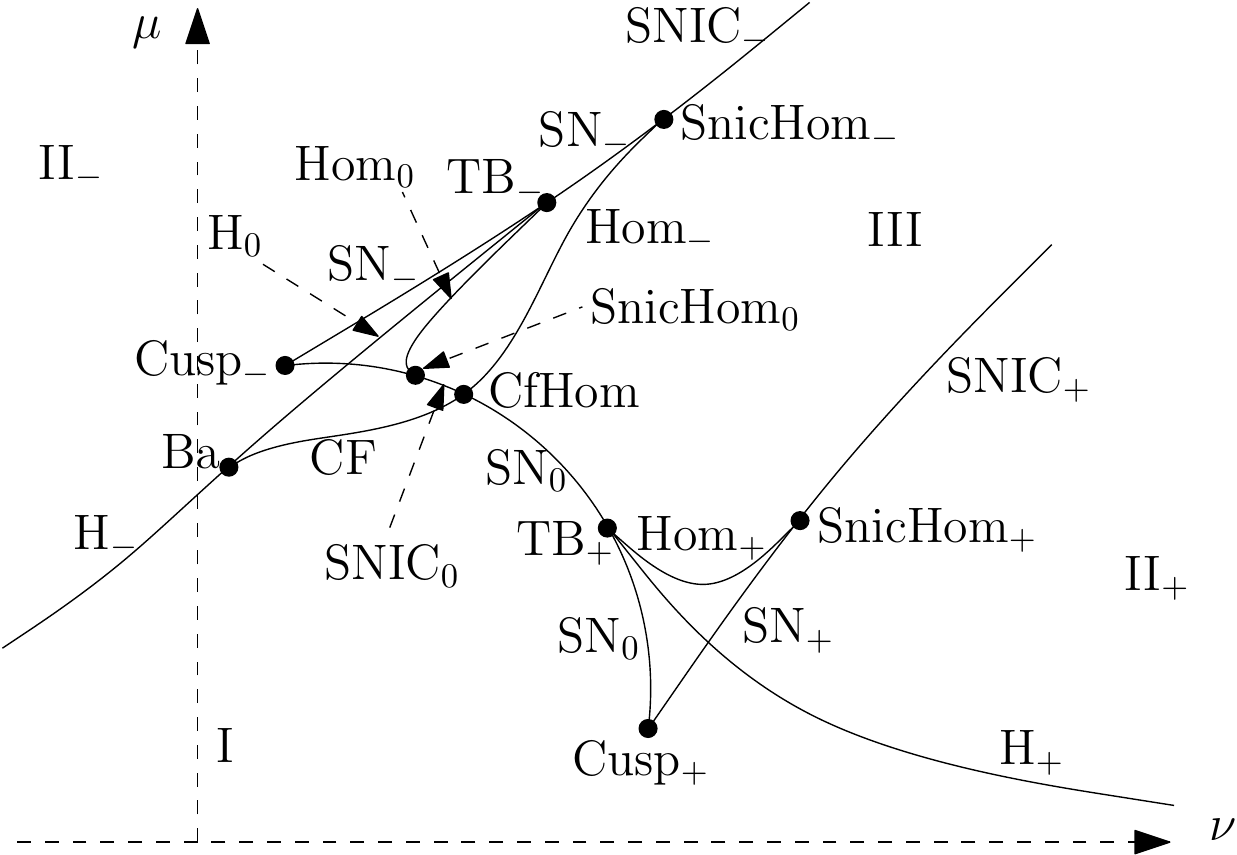}
  \end{center}
  \caption{Schematic of the bifurcations of the normal form
    \eqref{complexNFepsilon}. There are seven curves of global
    bifurcations, Hom$_\pm$, Hom$_0$ (homoclinic collisions of a limit
    cycle with a saddle), CF (a cyclic-fold), SNIC$_\pm$, and SNIC$_0$,
    and nine codimension-two points (black circles). The regions around
    the codimension-two points have been enhanced for clarity.}
  \label{epsilon_globals}
\end{figure}

\begin{figure}
  \begin{center}
    \begin{tabular}{m{15pt}m{0.6\linewidth}}
      $(a)$ & \includegraphics[width=\linewidth]{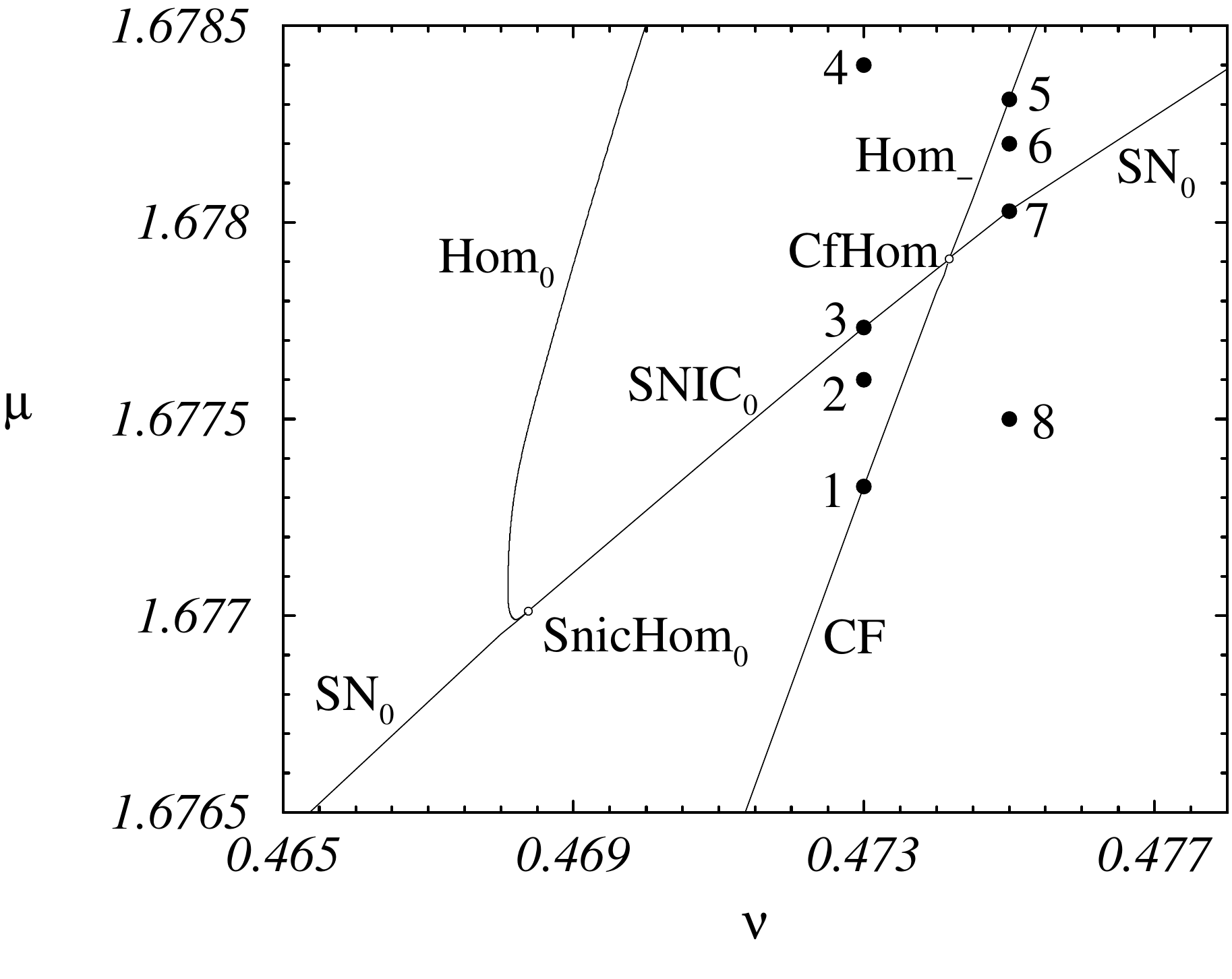}
    \end{tabular}\qquad\qquad

    \begin{tabular}{ccc}
      & \normalsize $(b)$: CfHom & \\[5pt]
      \includegraphics[width=0.3\linewidth]{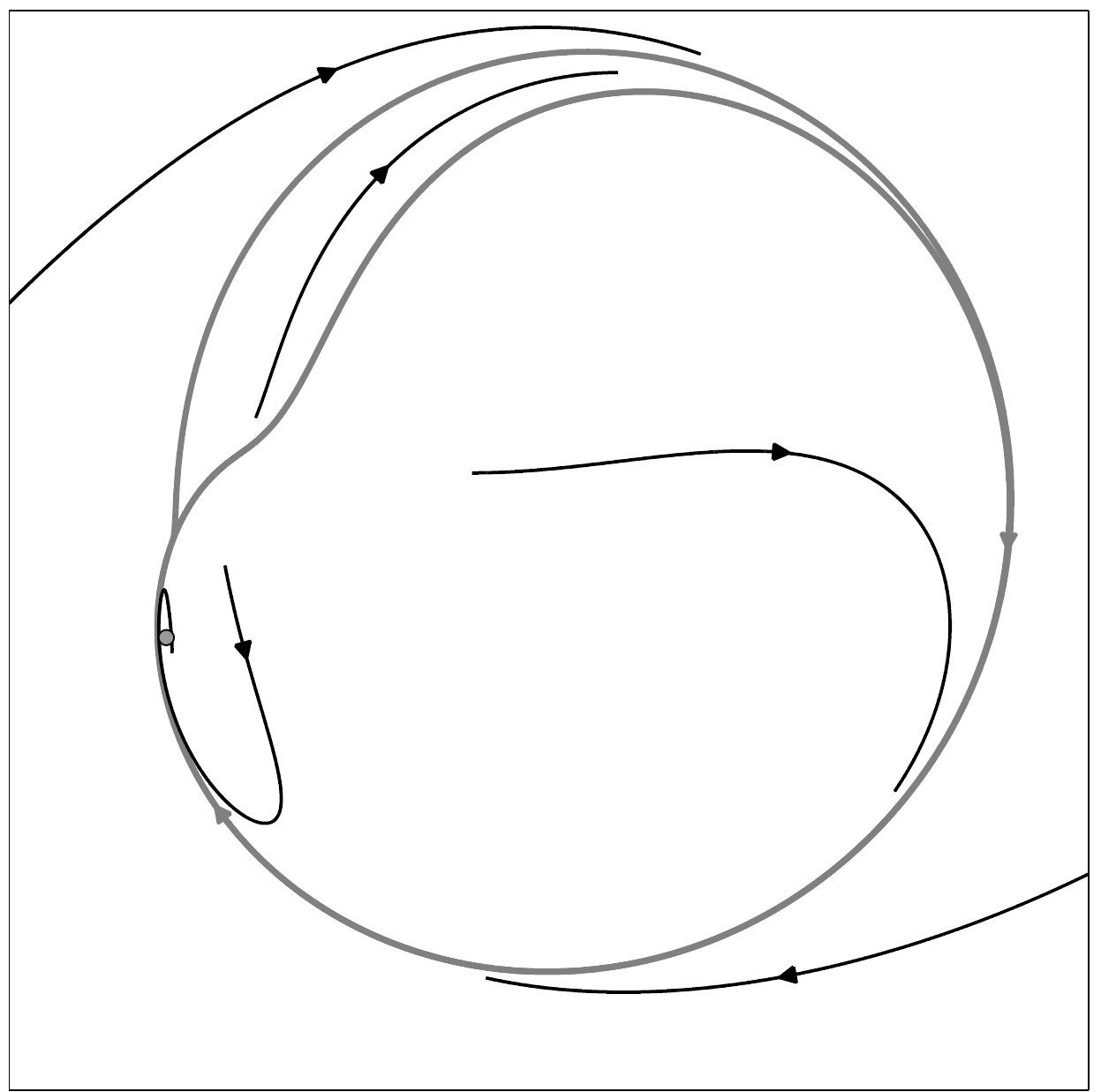} &
      \includegraphics[width=0.3\linewidth]{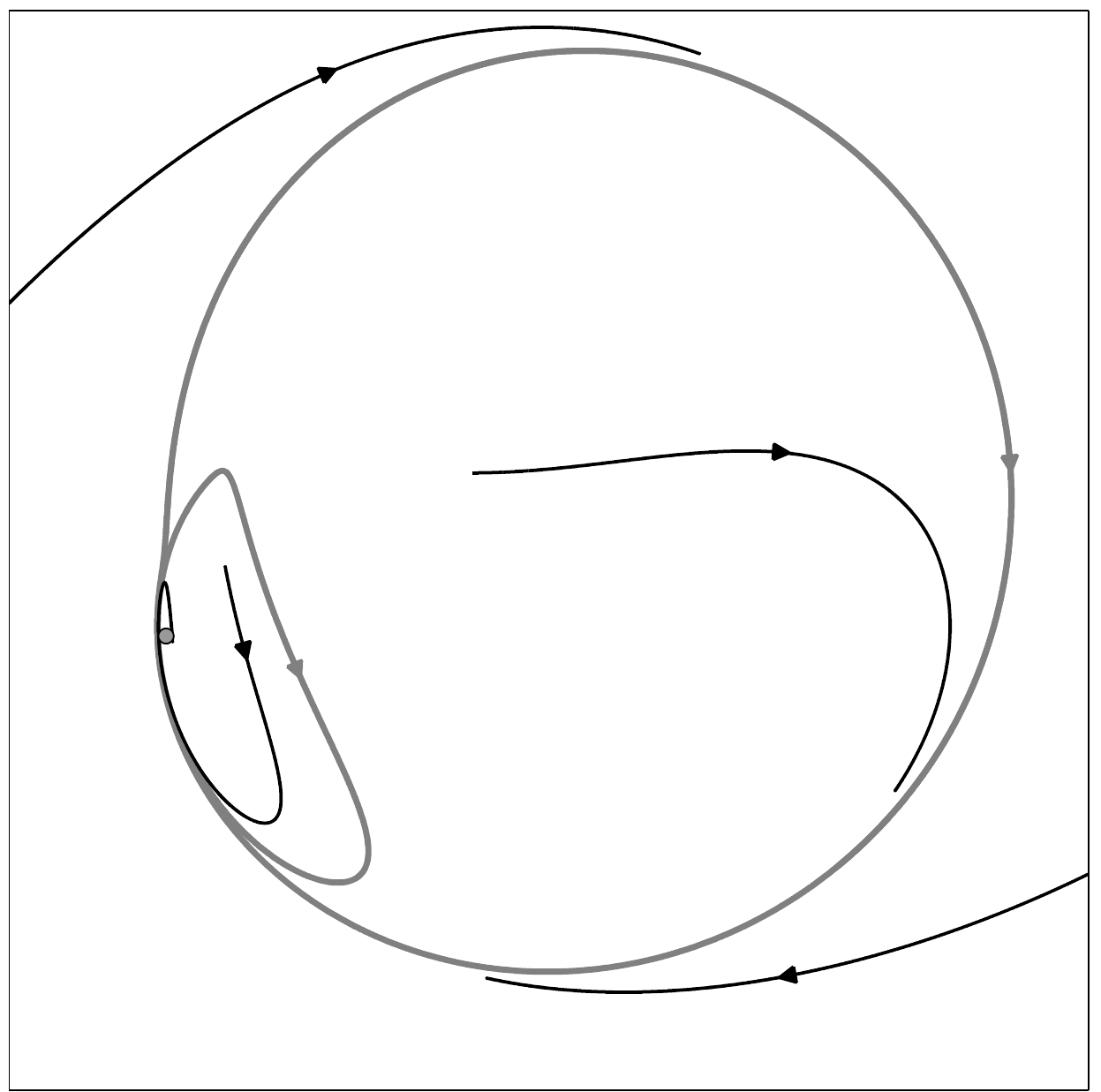} &
      \includegraphics[width=0.3\linewidth]{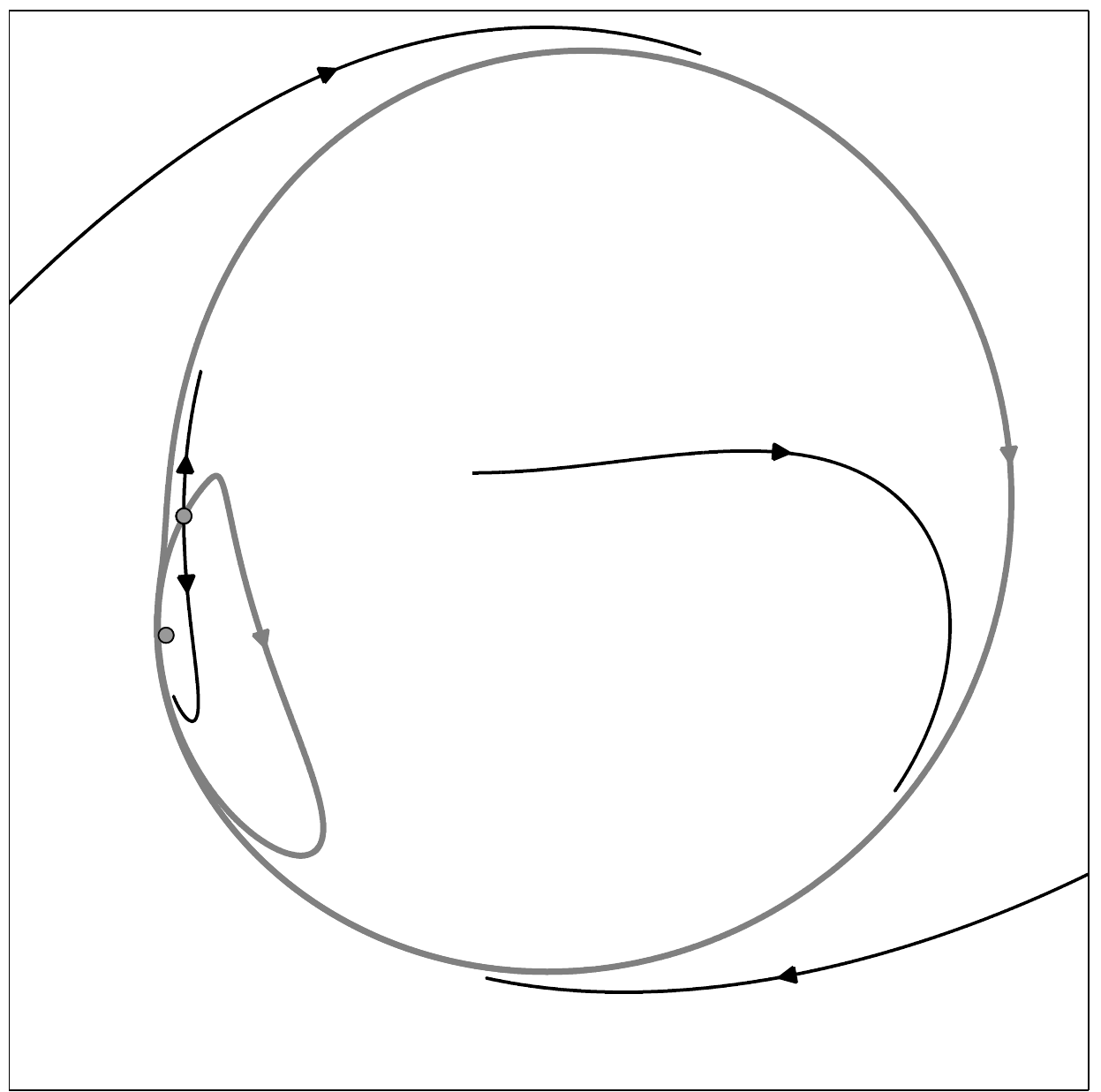} \\
      1: $\mu=1.67732854$ & 2: $\mu=1.6776$ & 3: $\mu=1.6777322$
    \end{tabular}\medskip
    
    \begin{tabular}{ccccc}
      \includegraphics[width=0.12\linewidth]{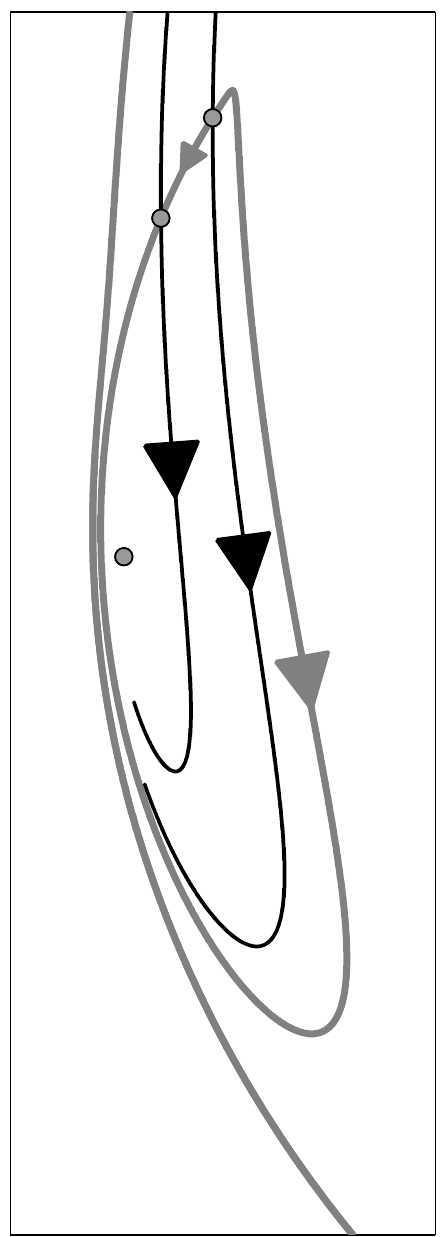} &
      \includegraphics[width=0.12\linewidth]{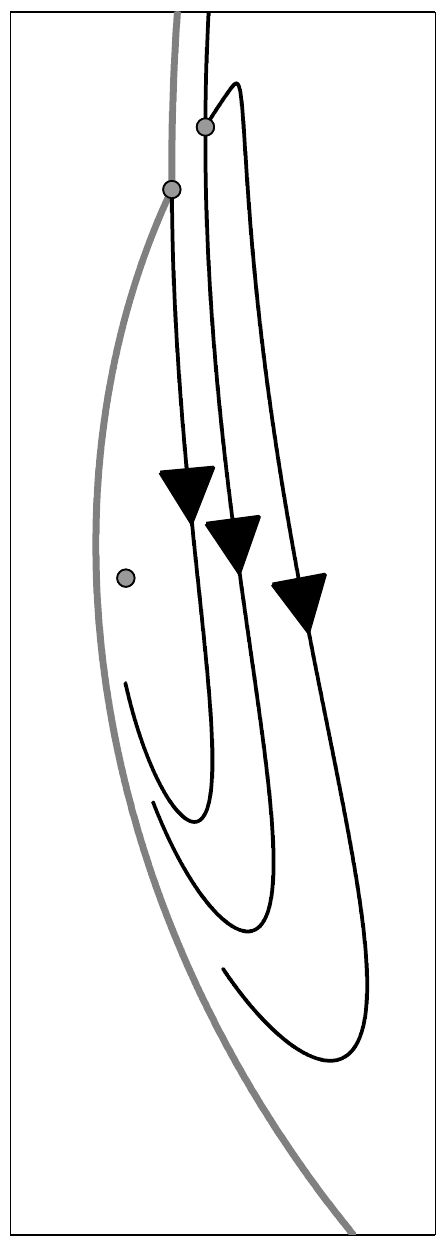} &
      \includegraphics[width=0.12\linewidth]{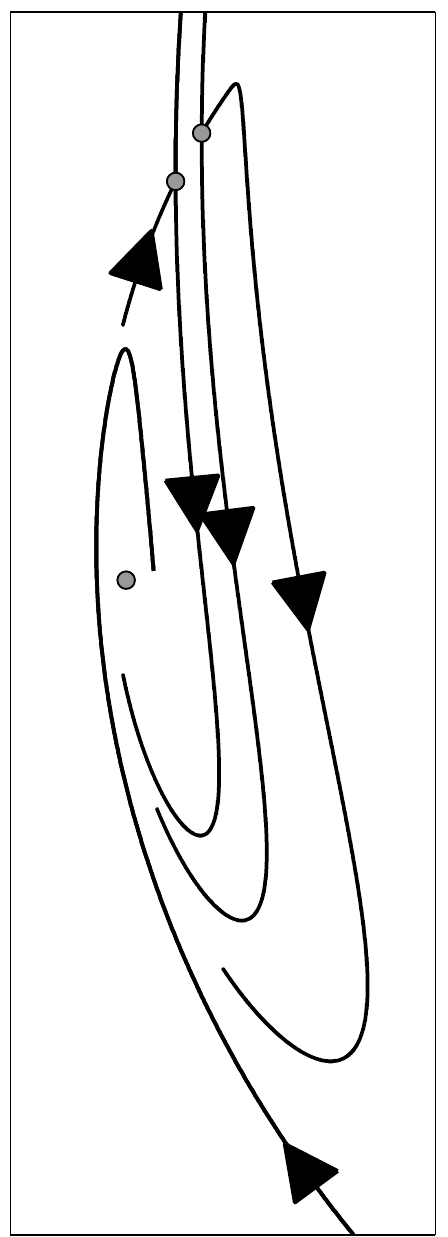} &
      \includegraphics[width=0.12\linewidth]{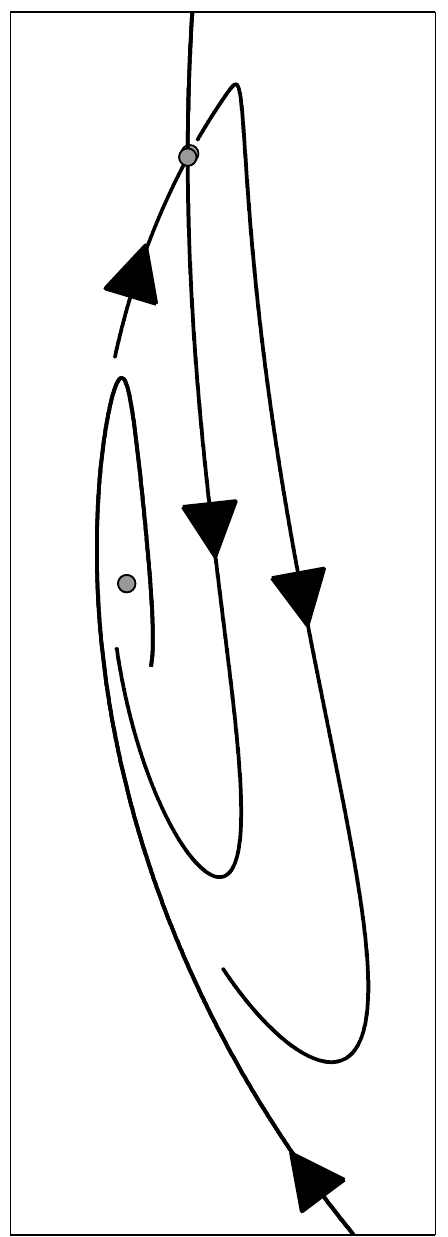} &
      \includegraphics[width=0.12\linewidth]{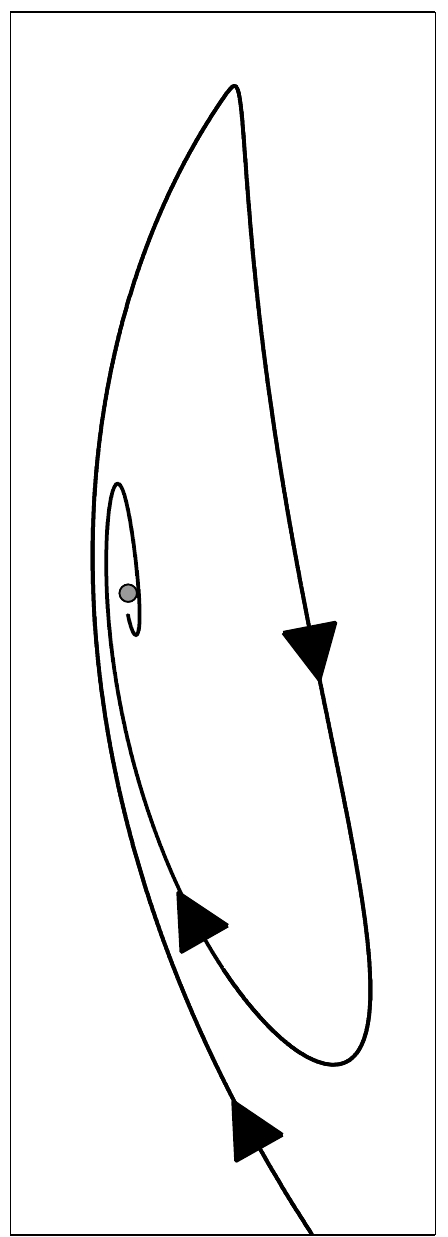} \\
      4: $\mu=1.6784$ & 5: $\mu=1.678312$ & 6: $\mu=1.6782$ & 
      7: $\mu=1.678028$ & 8: $\mu=1.6775$
    \end{tabular}
  \end{center}
  \caption{$(a)$ Zoom in parameter space around the codimension-two
    global bifurcation points SnicHom and CfHom. $(b)$ Phase portraits
    around CfHom; plots 1 to 4 at $\nu=0.473$, plots 5 to 8 at
    $\nu=0.475$, for $\mu$ as specified.}
  \label{CfHom}
\end{figure}

From the Takens--Bogdanov points, dynamical systems theory says that
two curves of homoclinic bifurcations emerge, resulting in global
bifurcations around these points. Moreover, the stable limit cycles in
regions II$_\pm$ do not exist in region III, so they must disappear in
additional bifurcations. These additional global bifurcations have
been explored numerically and using dynamical systems theory; these
are summarized in figure~\ref{epsilon_globals}. There are nine
codimension-two points organizing the dynamics of the normal form
\eqref{complexNFepsilon}. Apart from the cusp and Takens--Bogdanov
points already found, Cusp$_\pm$ and TB$_\pm$, there are five new
points: Ba, CfHom and three different Snic-Homoclinic bifurcations,
SnicHom$_\pm$ and SnicHom$_0$. On H$_-$, before crossing the SN$_0$
curve, the Hopf bifurcation becomes subcritical at the Bautin point
Ba, and from this point a curve of cyclic-folds CF appears. This curve
is the limit of the subcritical region where two periodic solutions,
$C_-$ and $C_0$, exist and they merge on CF. $C_0$ is the unstable
limit cycle born in the branch of H$_-$ between Cusp$_-$ and TB$_-$,
from now on termed the Hopf curve H$_0$. Inside the pinning region
these two periodic solutions disappear when they collide with a saddle
fixed point along the curves Hom$_0$ ($C_0$ collision) and Hom$_-$
($C_-$ collision), in the neighborhood of the Takens--Bogdanov point
TB$_-$, where the homoclinic curve Hom$_0$ is born. Away from the
TB$_-$, for increasing values of $\mu$ and $\nu$, the curve Hom$_-$
becomes tangent to and collides with the SN$_-$ curve, in a
SnicHom$_-$ bifurcation. Very close to SnicHom$_-$, SN$_-$ closely
followed by Hom$_-$ become indistinguishable from the SNIC$_-$
bifurcation, in the same sense as discussed in section
\S\ref{glob_bif_subsect}.

The two curves Hom$_-$ and Hom$_0$ on approaching SN$_0$, result in a
couple of codimension-two bifurcation points, SnicHom$_0$ and
CfHom. The arc of the curve SN$_0$ between the two new points
SnicHom$_0$ and CfHom, is a curve of saddle-node bifurcations taking
place on the limit cycle $C_0$, resulting in the SNIC$_0$ bifurcation
curve, as shown in figure~\ref{epsilon_globals}, and in more detail in
the numerically computed inset figure~\ref{CfHom}$(a)$. SnicHom$_0$ is
exactly the same bifurcation as SnicHom$_\pm$.

The cyclic fold bifurcation curve CF intersects the SNIC$_0$
bifurcation curve past the SnicHom$_0$ point, i.e.\ when the SN$_0$
curve is a line of SNIC bifurcations, at the point CfHom. On the SNIC
curve, one of the limit cycles born at CF undergoes a SNIC
bifurcation. At the point CfHom, the SNIC bifurcation happens
precisely when both limit cycles are born at CF: it is a saddle-node
bifurcation of fixed points taking place on a saddle-node bifurcation
of limit circles. After the CfHom point, the SNIC curve becomes an
ordinary saddle-node bifurcation curve, and there is an additional
homoclinic bifurcation curve emerging from this point CfHom, Hom$_-$
in figure~\ref{epsilon_globals}. The two limit cycles born at CF exist
only on one side of the CF line, so when following a closed path
around the CfHom point they must disappear. One of them undergoes a
SNIC bifurcation on the SNIC curve and disappears. The other limit
cycle collides with the saddle point born at the saddle-node curve
SN$_0$ (or SNIC$_0$) and disappears at the homoclinic collision
Hom$_-$. Numerically computed phase portraits illustrating these
processes around the CfHom point are shown in figure~\ref{CfHom}.

The other Takens--Bogdanov point TB$_+$ does not present any
additional complications. The homoclinic curve Hom$_+$ emerging from
it approaches and intersect the SN$_+$ curve in a SnicHom$_+$
codimension two point, as shown in figure~\ref{epsilon_globals}. For
large values of $\mu^2+\nu^2$ the stable limit cycles in regions
II$_\pm$ disappear at SNIC$_\pm$ (saddle-node on an invariant circle)
bifurcation curves. On these curves, a saddle-node bifurcation of
fixed points takes place on top of the limit cycle, and the cycle
disappears in an infinite-period bifurcation. What remains, and is
observable, is the stable fixed point born at the saddle-node.

The width of the pinning $w(d,\epsilon)$ region away from the origin
(large $s$) is easy to compute from \eqref{cusp_curve}:
\begin{equation}
  d=u\sim (3s^2/4)^{1/3},\quad w=2v\sim (16/\sqrt{3}s)^{1/3}
  \Rightarrow w=2/\sqrt{d}.
\end{equation}
Restoring the $\epsilon$ dependence, we obtain
$w(d,\epsilon)=2\epsilon/\sqrt{d}$. The pinning region becomes
narrower away from the bifurcation point, and the width is
proportional to $\epsilon$, the size of the imperfection.

\section{Symmetry breaking of $SO(2)$ to $Z_2$, with an $\epsilon\bar
z$ term}\label{Sec_bar_z}

The $\epsilon_2\bar z$ term in \eqref{NFcompleteunfold} corresponds to
breaking $SO(2)$ symmetry in a way that leaves a system with $Z_2$
symmetry, corresponding to invariance under a half turn. The normal
form to be analyzed is \eqref{NFsinglepert} with $p=1$, $q=0$ and
$\epsilon=1$:
\begin{equation}\label{complexNF}
 \dot z=z(\mu+i\nu-c|z|^2)+\bar z.
\end{equation}
The new normal form \eqref{complexNF} is still invariant to $z\to -z$,
or equivalently, the half-turn $\phi\to\phi+\pi$. This is all that
remains of the $SO(2)$ symmetry group, which is reduced to $Z_2$,
generated by the half-turn. In fact, the $Z_2$ symmetry implies that
$P(z,\bar z)$ in $\dot z=P(z,\bar z)$ must be odd: $P(-z,-\bar
z)=-P(z,\bar z)$, which is \eqref{NFsymSO2} for $\theta=\pi$, the half
turn. Therefore, \eqref{complexNF} is the unfolding corresponding to
the symmetry breaking of $SO(2)$ to $Z_2$. This case has also been
analyzed in \citet{Gam85,Wag01,BDV08,SaWa10}.

Writing the normal form \eqref{complexNF} in terms of the modulus and
phase of $z=r\ce^{\ci\phi}$ gives
\begin{equation}\label{realNF}
  \begin{aligned}
    & \dot r=r(\mu-ar^2)+r\cos 2\phi,\\
    & \dot\phi=\nu-br^2-\sin 2\phi.
  \end{aligned}
\end{equation}

\subsection{Fixed points and their bifurcations}

\begin{figure}
  \begin{center}
    \includegraphics[width=0.98\linewidth]{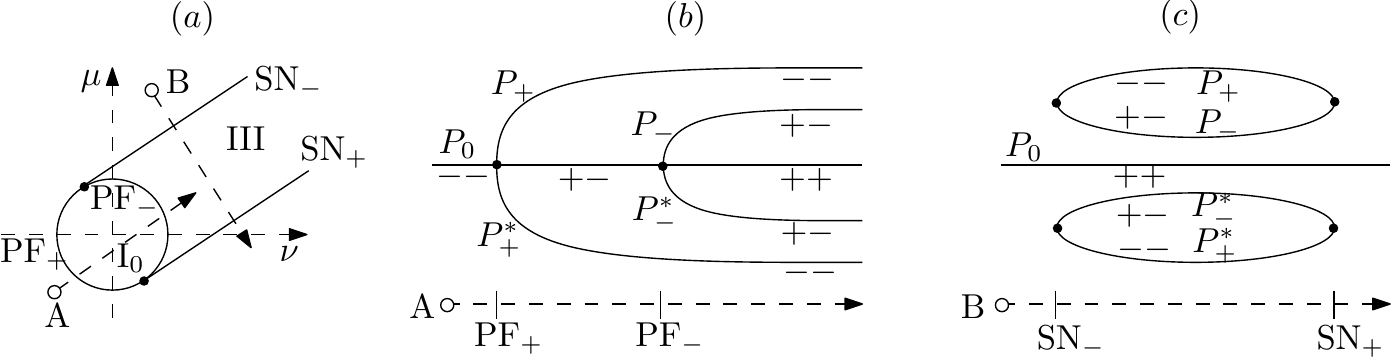}\vspace*{-8pt}
  \end{center}
  \caption{Steady bifurcations of fixed points corresponding to the
    normal form \eqref{realNF}: $(a)$ regions in
    parameter space delimited by the fixed points and their steady
    bifurcations; $(b)$ and $(c)$ bifurcations along the paths A and B
    shown in $(a)$, respectively.}
  \label{impHopf_steady}
\end{figure}

The normal form \eqref{complexNF}, or \eqref{realNF}, admits up to
five fixed points. One is the origin $r=0$, the trivial solution
$P_0$. The other fixed points come in two pairs of $Z_2$-symmetric
points: one is the pair $P_+=r_+\ce^{\ci\phi_+}$ and
$P^*_+=-r_+\ce^{\ci\phi_+}$, and the other pair is
$P_-=r_-\ce^{\ci\phi_-}$ and $P^*_-=-r_-\ce^{\ci\phi_-}$.
Coefficients $r_\pm$ and $\phi_\pm$ are given by
\begin{gather}
 r^2_\pm=a\mu+b\nu\pm\Delta,\quad \phi_\pm=(\alpha_0\pm\alpha_1)/2,\\
 \Delta^2=1-(a\nu-b\mu)^2,\quad \ce^{\ci\alpha_1}=a\nu-b\mu-\ci\Delta.
\end{gather}
The details of the computations are given in \ref{Appendix_bar_z}.
There are three different regions in the $(\mu,\nu)$-parameter plane:
region III, where there exist five fixed points, $P_0$, $P_\pm$ and
$P^*_\pm$; region I$_0$ where three fixed poins exist, $P_0$, $P_+$
and $P^*_+$; and the rest of the parameter space where only $P_0$
exists. These three regions are separated by four curves along which
steady bifurcations between the different fixed points take place, as
shown in figure~\ref{impHopf_steady}. Along the semicircle
\begin{equation}
  \text{PF}_+:~ \mu^2+\nu^2=1 \text{ and } a\mu+b\nu<0,
\end{equation}
the two symmetrically-related solutions $P_+$ and $P^*_+$ are born in
a pitchfork bifurcation of the trivial branch $P_0$. Along the
semicircle
\begin{equation}
  \text{PF}_-:~ \mu^2+\nu^2=1 \text{ and } a\mu+b\nu>0,
\end{equation}
the two symmetrically-related solutions $P_-$ and $P^*_-$ are born
in a pitchfork bifurcation of the trivial branch $P_0$. Along
the two half-lines
\begin{equation}
  \text{SN$_+$}:~ \mu=(a\nu-1)/b,\quad \text{SN$_-$}:~ \mu=(a\nu+1)/b,\quad
  \text{both with } a\mu+b\nu>0,
\end{equation}
a saddle-node bifurcation takes place. It is a double saddle-node, due
to the $Z_2$ symmetry; we have one saddle-node involving $P_+$ and
$P_-$, and the $Z_2$-symmetric saddle-node between $P^*_+$ and
$P^*_-$. Bifurcation diagrams along the paths A and B in
figure~\ref{impHopf_steady}$(a)$ are shown in parts $(b)$ and $(c)$ of
the same figure.

We can compare with the original problem with $SO(2)$ symmetry,
corresponding to $\epsilon=0$. In order to do that, the $\epsilon$
dependence will be restored in this paragraph.  The single line L
($\mu=\nu\tan\alpha_0$) where $\omega=0$ and nontrivial fixed points
exist in the perfect problem, becomes a region of width $2\epsilon$ in
the imperfect problem, where up to four fixed points exist, in
addition to the base state $P_0$; they are the remnants of the circle
of fixed points in the original problem.  Solutions with $\omega=0$,
that existed only along a single line in the absence of imperfections,
now exist in a region bounded by the semicircle $F_+$ and the
half-lines $a\nu-b\mu=\pm\epsilon$; this region will be termed the
pinning region. It bears some relationship with the frequency-locking
regions appearing in Neimark-Sacker bifurcations, in the sense that
here we also have frequency locking, but with $\omega=0$.  The width
of the pinning region is proportional to $\epsilon$, a measure of the
breaking of $SO(2)$ symmetry due to imperfections.

In the absence of imperfections ($\epsilon=0$) the $P_0$ branch looses
stability to a Hopf bifurcation along the curve $\mu=0$. Let us
analyze the stability of $P_0$ in the imperfect problem. Using
Cartesian coordinates $z=x+\ci y$ in \eqref{complexNF} we obtain
\begin{equation}\label{cartesianNF}
  \begin{pmatrix} \dot x \\ \dot y \end{pmatrix}=
  \begin{pmatrix} \mu+1 & -\nu \\ \nu & \mu-1 \end{pmatrix}
  \begin{pmatrix} x \\ y \end{pmatrix}-
  (x^2+y^2)\begin{pmatrix} ax-by \\ bx+ay \end{pmatrix}.
\end{equation}
The eigenvalues of $P_0$ are the eigenvalues of the linear part of
\eqref{cartesianNF}, $\lambda_\pm=\mu\pm \sqrt{1-\nu^2}$. There is a
Hopf bifurcation ($\Im\lambda_\pm\ne0$) when $\mu=0$ and $|\nu|>1$,
i.e.\ on the line $\mu=0$ outside region II; the Hopf frequency is
$\omega=\sign(\nu)\sqrt{\nu^2-1}$. The sign of $\omega$ is the same as
the sign of $\nu$, from the $\dot\phi$ equation in
\eqref{realNF}. Therefore, we have a Hopf bifurcation with positive
frequency along H$_+$ ($\mu=0$ and $\nu>1$) and a Hopf bifurcation
with negative frequency along H$_-$ ($\mu=0$ and $\nu<1$). The
bifurcated periodic solutions are stable limit cycles $C_+$ and $C_-$,
respectively.

\begin{figure}
  \begin{center}
    \includegraphics[width=.6\linewidth]{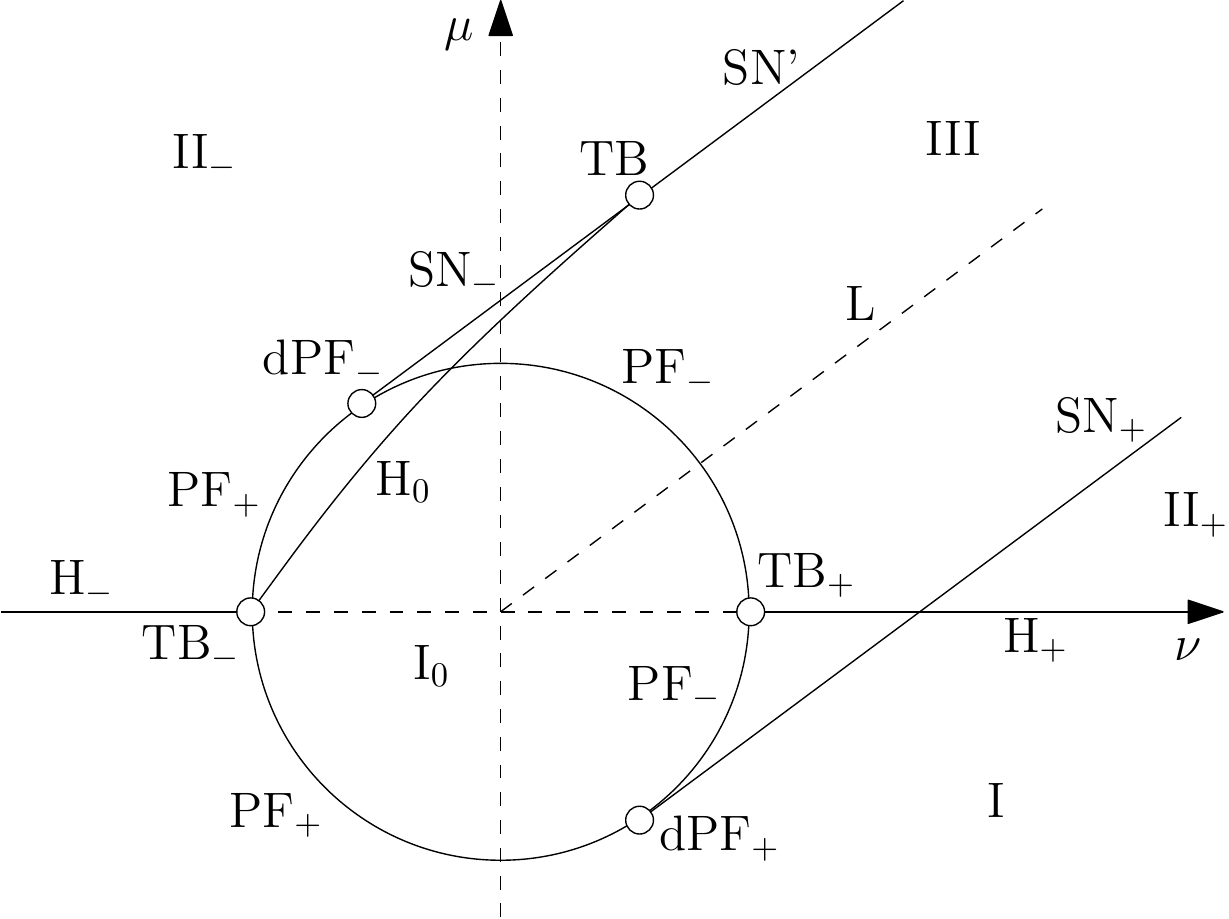}\vspace*{-12pt}
  \end{center}
  \caption{Local bifurcations of fixed points in the symmetry breaking
    of $SO(2)$ to $Z_2$ case. Codimension-one bifurcation curves: Hopf
    curves H$_\pm$ and H$_0$, pitchforks PF$_\pm$, and saddle-nodes SN$_\pm$.
    Codimension-two bifurcation points: degenerate pitchforks
    dPF$_\pm$, Takens--Bogdanov TB and TB$_\pm$. L is the
    zero-frequency curve in the $SO(2)$ symmetric case.}
  \label{bif_fixed_points}
\end{figure}

The Hopf bifurcations of the $P_\pm$ and $P^*_\pm$ points can be
studied analogously. The eigenvalues of the Jacobian of the right-hand
side of \eqref{cartesianNF} at a fixed point characterize the
different bifurcations that the fixed point can undergo. Let $T$ and
$D$ be the trace and determinant of $J$. The eigenvalues are given by
\begin{equation}
  \lambda^2-T\lambda+D=0\quad\Rightarrow\quad
  \lambda=\frac{1}{2}(T\pm\sqrt{Q}),\quad Q=T^2-4D.
\end{equation}
A Hopf bifurcation takes place for $T=0$ and $Q<0$. The equation $T=0$
at the four points $P_\pm$ results in the ellipse
$\mu^2-4ab\mu\nu+4a^2\nu^2=4a^2$ (see \ref{Appendix_bar_z} for
details). This ellipse is tangent to SN$_-$ at the point
$(\mu,\nu)=(2b,(b^2-a^2)/a)$.  The condition $Q<0$ is only satisfied
by $P_+$ and $P^*_+$ on the elliptic arc H$_0$ from $(\mu,\nu)=(0,-1)$
to $(2b,(b^2-a^2)/a)$:
\begin{equation}
  \mu=2ab\nu+2a\sqrt{1-a^2\nu^2},\quad \nu\in[-1,(b^2-a^2)/a].
\end{equation}
The elliptic arc H$_0$ is shown in figure~\ref{bif_fixed_points}.
Along this arc a pair of unstable limit cycles $C_0$ and $C^*_0$ are
born around the fixed points $P_+$ and $P^*_+$, respectively.

\subsection{Codimension-two points}

The local codimension-one bifurcations of the fixed points are now
completely characterized. There are two curves of saddle-node
bifurcations, two curves of pitchfork bifurcations and three Hopf
bifurcation curves. These curves meet at five codimension-two
points. The analysis of the eigenvalues at these points, and of the
symmetry of the bifurcating points ($P_0$, $P_+$ and $P^*_+$),
characterizes these points as two degenerate pitchforks dPF$_\pm$, two
Takens--Bogdanov bifurcations with $Z_2$ symmetry TB$_\pm$,
and a double Takens--Bogdanov bifurcation TB, as shown in
figure~\ref{bif_fixed_points}.

\begin{figure}
  \begin{center}
    \includegraphics[width=0.82\linewidth]{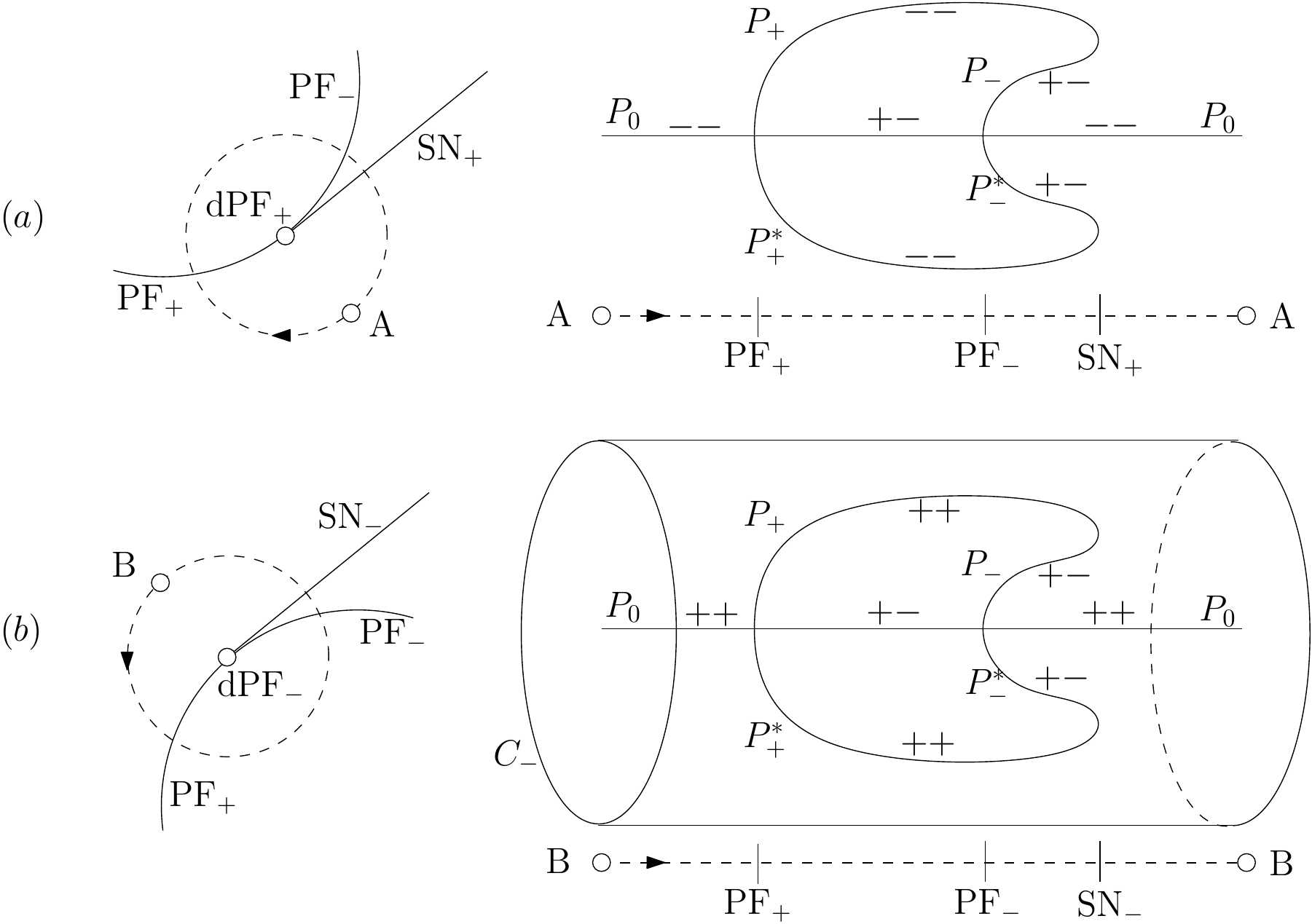}\vspace*{-6pt}
  \end{center}
  \caption{Schematics of the degenerate pitchfork bifurcations $(a)$
    dPF$_+$ and $(b)$ dPF$_-$. On the left, bifurcation curves emanating
    from dPF$_\pm$ in parameter space are shown, along with a closed
    one-dimensional path (dashed); shown on the right are schematics of
    the bifurcations along the closed path, starting and ending at A
    (dPF$_+$) and B (dPF$_-$).  The fixed point curves are labeled with
    the signs of their eigenvalues. $C_-$ is the periodic solution born
    at the curve H$_-$.}
  \label{barz_DPp}
\end{figure}

The degenerate pitchforks dPF$_\pm$ correspond to the
transition between supercritical and subcritical pitchfork
bifurcations. At these points, a saddle-node curve is born, and only
fixed points are involved in the neighboring dynamics. The only
difference between dPF$_+$ and dPF$_-$ is the stability of the base
state $P_0$; it is stable outside the circle $\mu^2+\nu^2=1$ at dPF$_+$
and unstable at dPF$_-$. Schematics of the bifurcations along a
one-dimensional path in parameter space around the dPF$_+$ and dPF$_-$
points are illustrated in figure~\ref{barz_DPp}. The main difference,
apart from the different stability properties of $P_0$, $P_+$ and
$P^*_+$, is the existence of the limit cycle $C_-$ surrounding the
three fixed points in case $(b)$, dPF$_-$.

\begin{figure}\setlength{\Figsize}{0.45\linewidth}
  \begin{center}
    \includegraphics[width=0.82\linewidth]{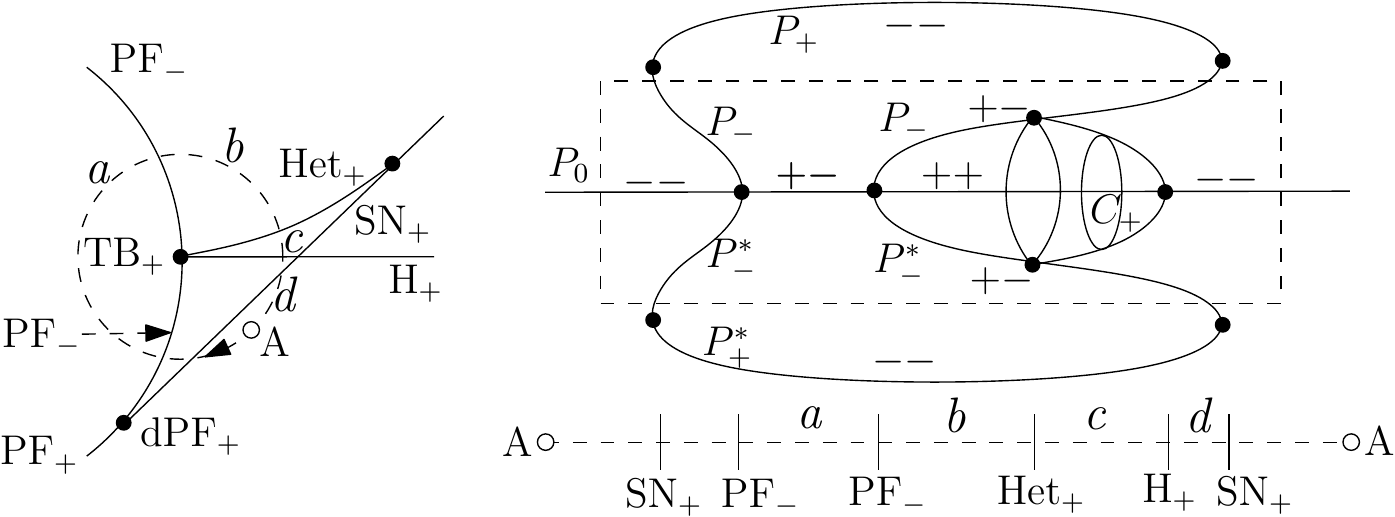} \\[20pt]
    \begin{tabular}{@{}cc@{}}
      $a$: $\mu=0.05,~\nu=0.9$ & $b$: $\mu=0.1,~\nu=1.1$ \\
      \includegraphics[width=\Figsize]{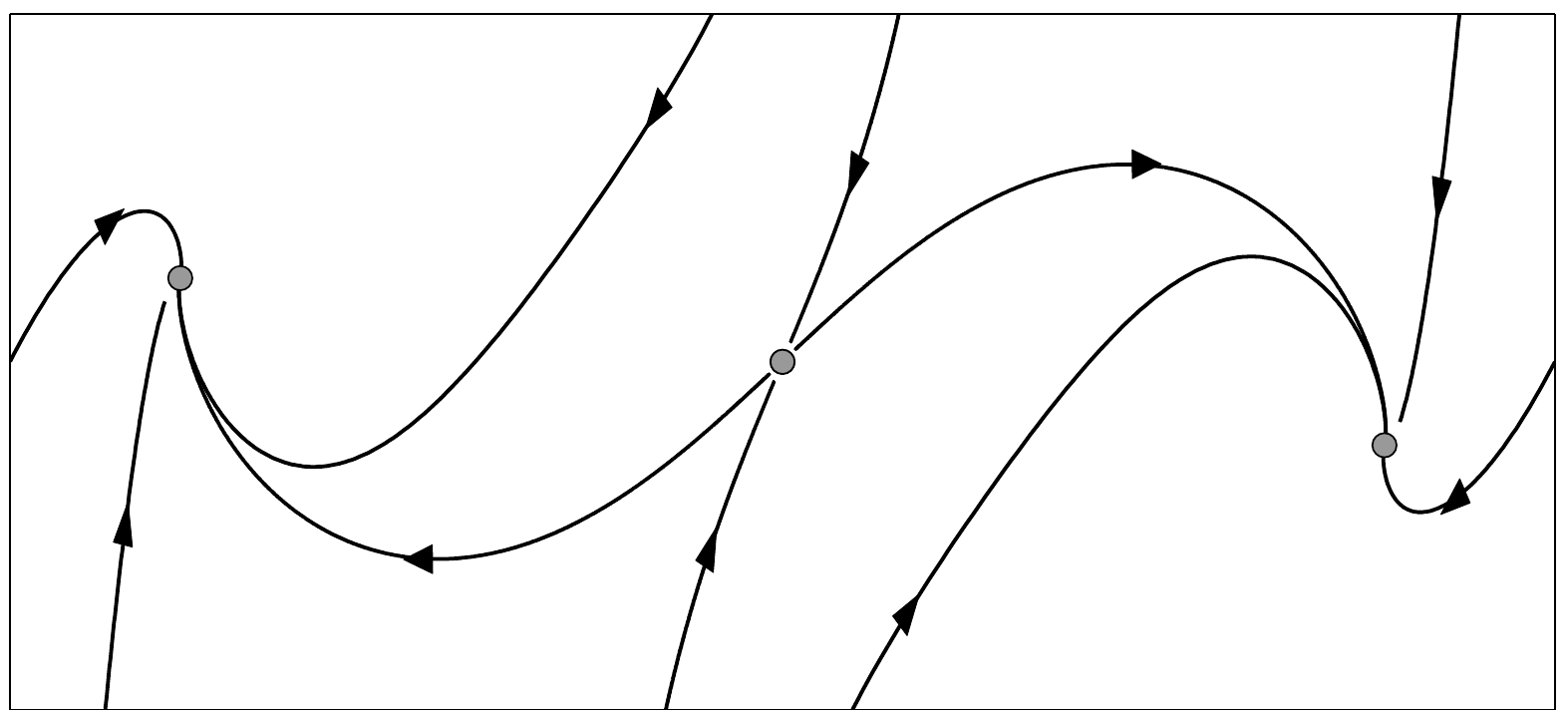} &
      \includegraphics[width=\Figsize]{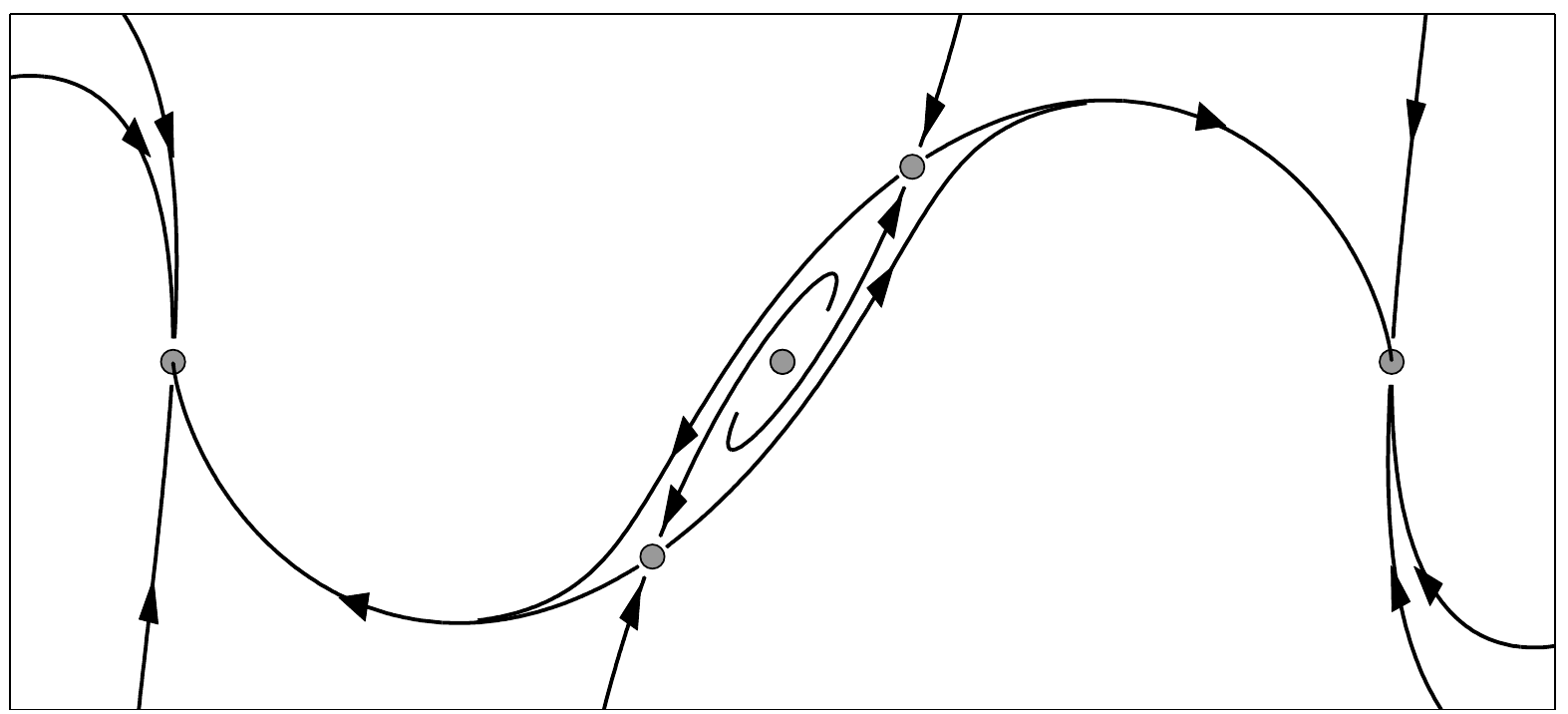} \\
      $c$: $\mu=0.025,~\nu=1.2$ & $d$:  $\mu=-0.05,~\nu=1.1$ \\
      \includegraphics[width=\Figsize]{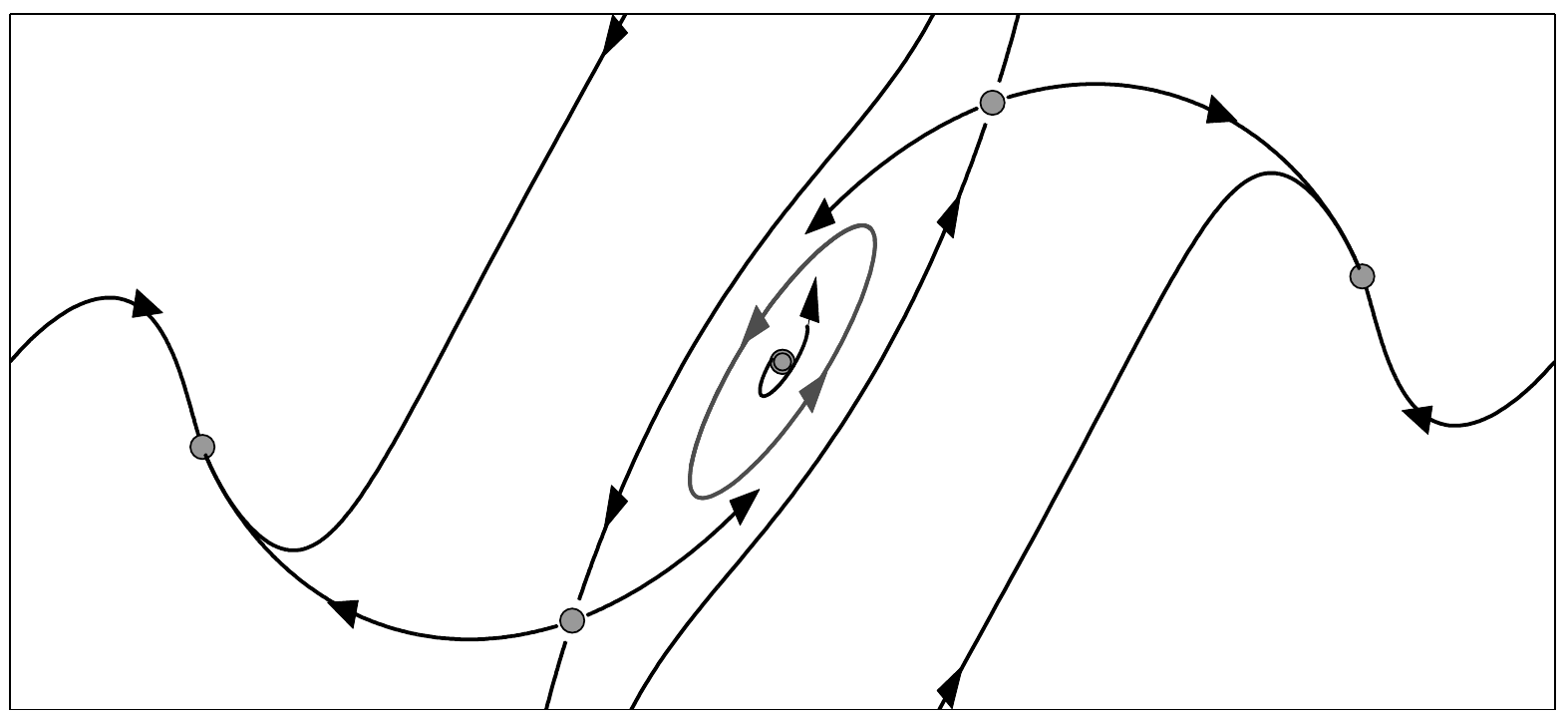} &
      \includegraphics[width=\Figsize]{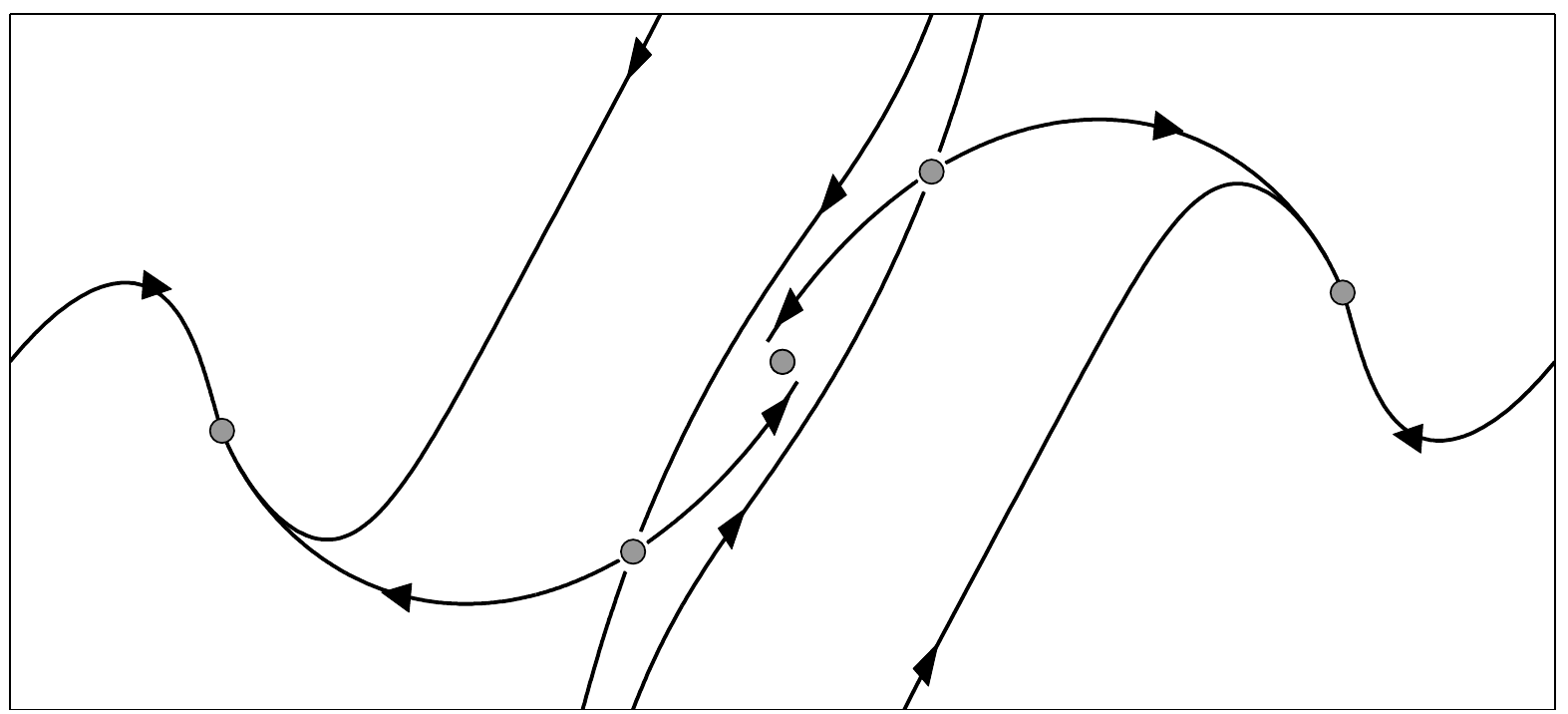} \\
    \end{tabular}\vspace*{-6pt}
  \end{center}
  \caption{Takens--Bogdanov bifurcation with $Z_2$ symmetry TB$_+$. The
    top left shows bifurcation curves emanating from TB$_+$ in parameter
    space, with a closed one-dimensional path. The top right shows a
    schematic of the bifurcations along the closed path, starting and
    ending at A.  The fixed point curves are labeled with the signs of
    their eigenvalues. $C_+$ is the periodic solution born at the curve
    H$_+$. The region inside the dashed rectangle on the right contains
    the states locally connected with the bifurcation TB$_+$. The bottom
    panels show four numerically computed phase portraits, at points
    labeled $a$, $b$, $c$ and $d$, for the specified parameter values.}
  \label{TB_sym_p}
\end{figure}

The Takens--Bogdanov bifurcation with $Z_2$ symmetry has two different
scenarios \citep{CLW94}, and they differ in whether one or two Hopf
curves emerge from the bifurcation point.  In our problem, bifurcation
point TB$_+$ has a single Hopf curve, H$_+$, while the TB$_-$ point
has two Hopf curves, H$_-$ and H$_0$, emerging from the bifurcation
point. The scenario TB$_+$ is depicted in figure~\ref{TB_sym_p},
showing the bifurcation diagram as well as the bifurcations along a
closed one-dimensional path around the codimension-two point. We have
also included the $P_+$ and $P^*_+$ solutions that merge with the
$P_-$ and $P^*_-$ fixed points along the saddle-node bifurcation curve
SN$_+$, although they are not locally connected to the codimension-two
point, in order to show all the fixed points in the phase space of
\eqref{complexNF}. A curve of global bifurcations, a heteroclinic
cycle Het$_+$ connecting $P_-$ and $P^*_-$, is born at TB$_+$. The
heteroclinic cycle is formed when the limit cycle $C_+$ simultaneously
collides with the saddles $P_-$ and $P^*_-$.

\begin{figure}\setlength{\Figsize}{0.3\linewidth}
  \begin{center}
    \includegraphics[width=0.90\linewidth]{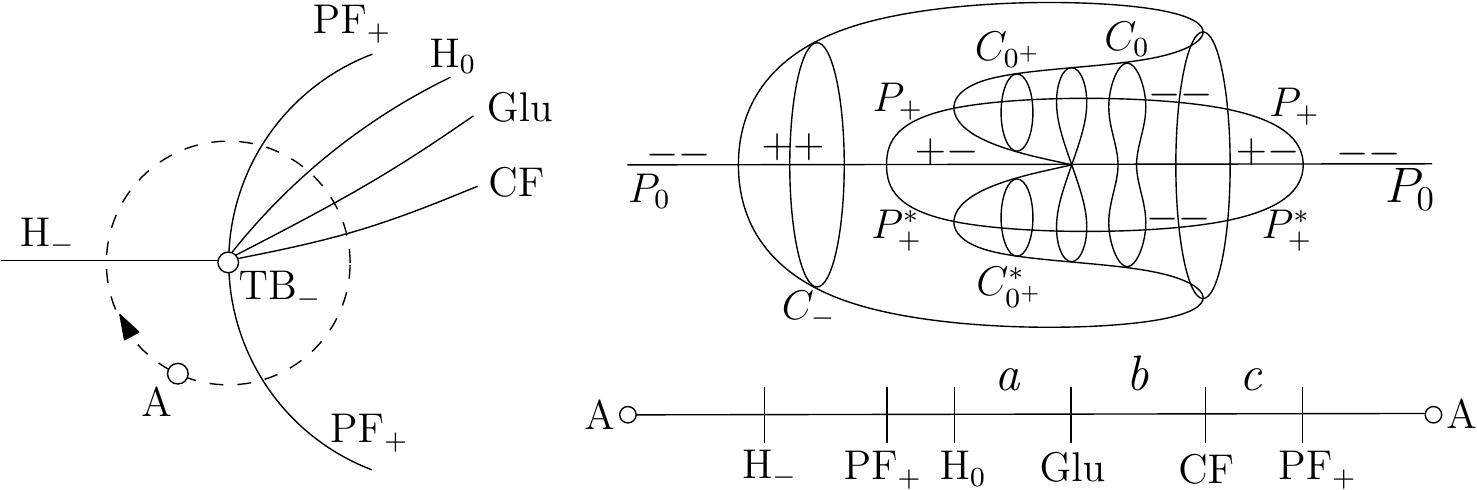} \\[20pt]
    \begin{tabular}{ccc}
      $a$: $\nu=-0.68$ & $b$: $\nu=-0.66$ & $c$: $\nu=-0.65$ \\
      \includegraphics[width=\Figsize]{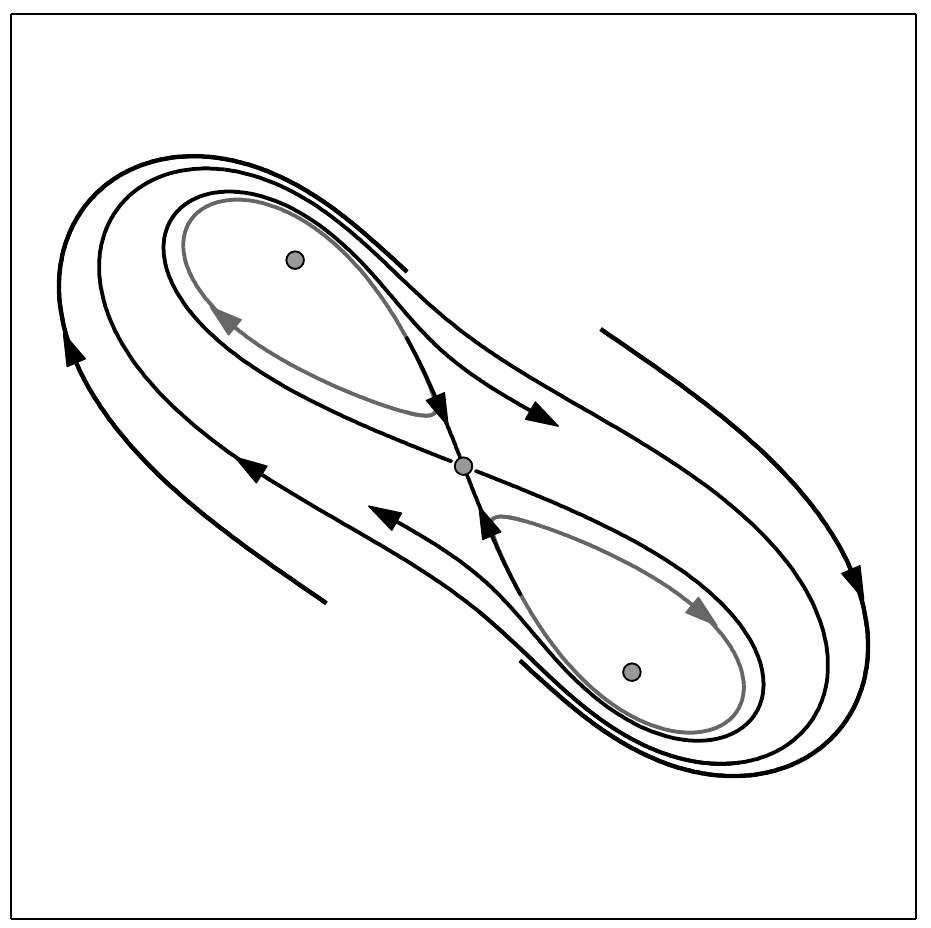} &
      \includegraphics[width=\Figsize]{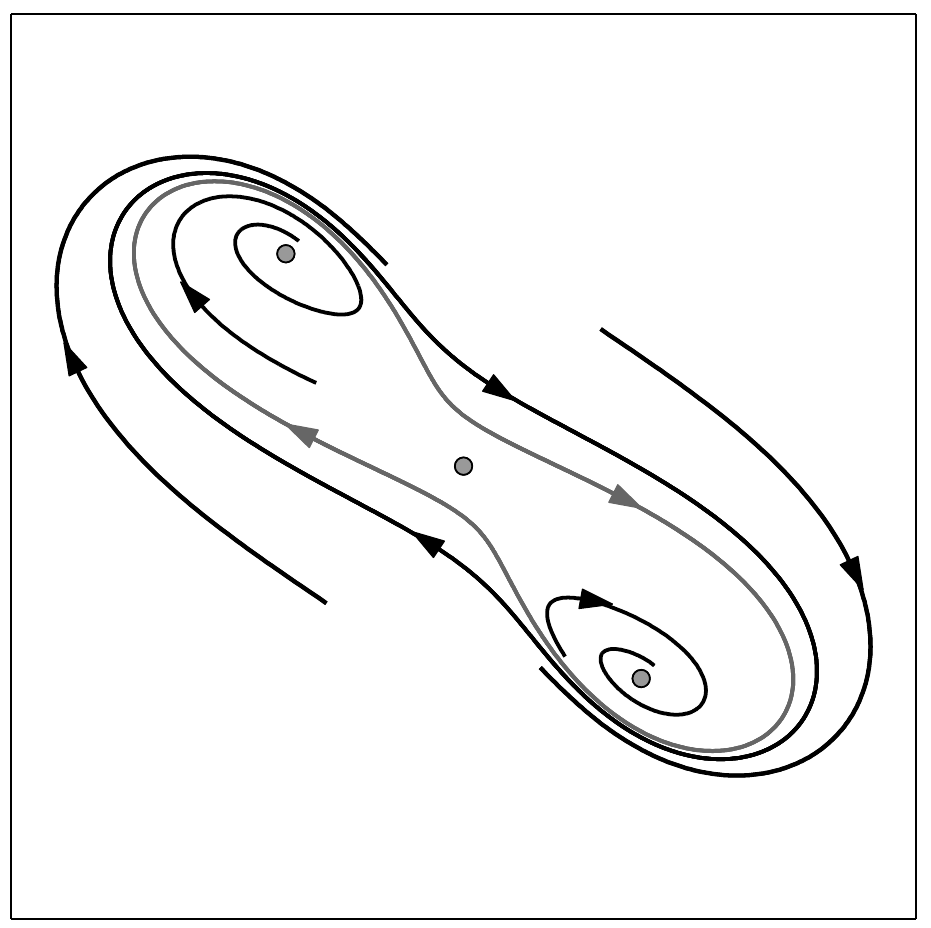} &
      \includegraphics[width=\Figsize]{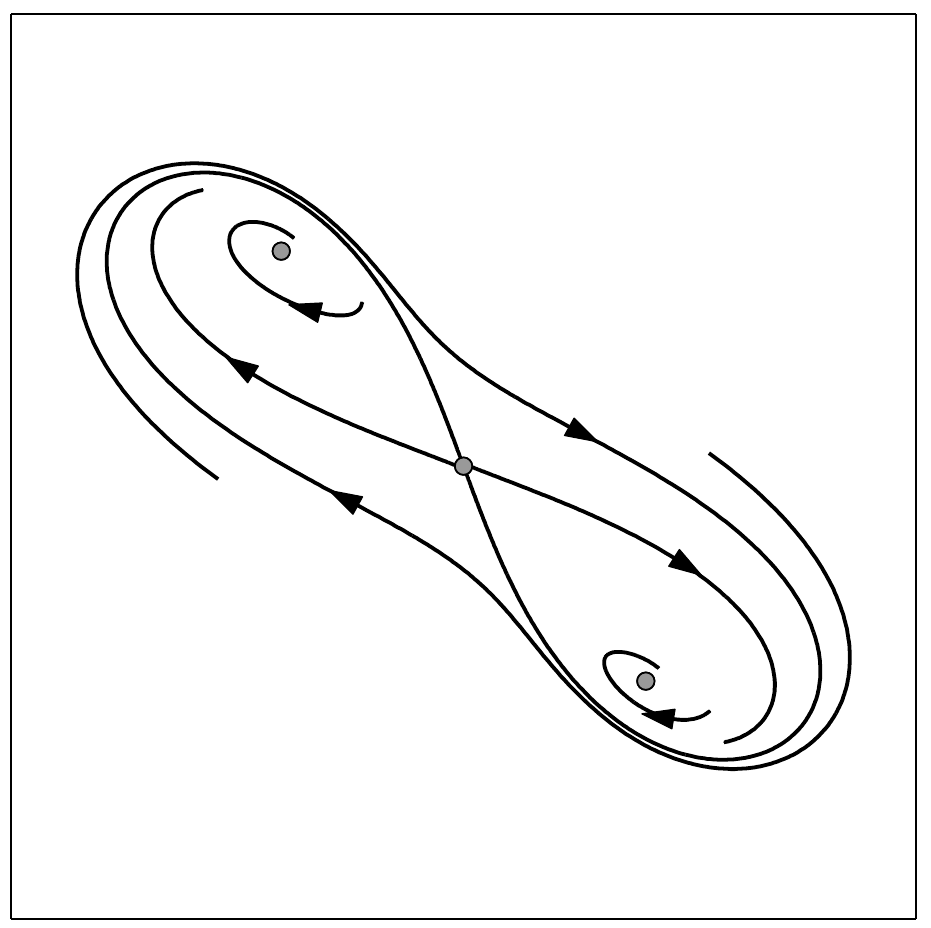} \\
    \end{tabular}
  \end{center}
  \caption{Takens--Bogdanov bifurcation with $Z_2$ symmetry, TB$_-$. The
    top left shows bifurcation curves emanating from TB$_-$ in parameter
    space, with a closed one-dimensional path. The top right shows a
    schematic of the bifurcations along the closed path, starting and
    ending at A.  The fixed point curves are labeled with the signs of
    their eigenvalues. $C_-$ is the periodic solution born at the curve
    H$_-$; $C_{0^+}$ and $C^*_{0^+}$ are the unstable cycles born
    simultaneously at the Hopf bifurcation H$_0$ around the fixed points
    $P_+$ and $P^*_+$; and $C_0$ is the cycle around both fixed points
    that remains after the gluing bifurcation. The bottom panels show
    three numerically computed phase portraits, at points labeled $a$,
    $b$ and $c$, for $\mu=0.5$ and $\nu$ as indicated, illustrating the
    gluing and cyclic fold bifurcations.}
  \label{TB_sym_m}
\end{figure}

The scenario TB$_-$, is depicted in figure~\ref{TB_sym_m}, showing the
bifurcation diagram and also the bifurcations along a closed
one-dimensional path around the codimension-two point. Two curves of
global bifurcations are born at TB$_-$. One corresponds to a gluing
bifurcation Glu, when the two unstable limit cycles $C_{0^+}$ and
$C^*_{0^+}$, born at the Hopf bifurcation H$_0$ around $P_+$ and
$P^*_+$, simultaneously collide with the saddle $P_0$; after the
collision a large cycle $C_0$ results, surrounding the three fixed
points $P_0$, $P_+$ and $P^*_+$. The second global bifurcation curve
corresponds to a saddle-node of cycles, where $C_0$ and $C_-$ collide
and disappear.

\begin{figure}
  \begin{center}
    \includegraphics[width=0.8\linewidth]{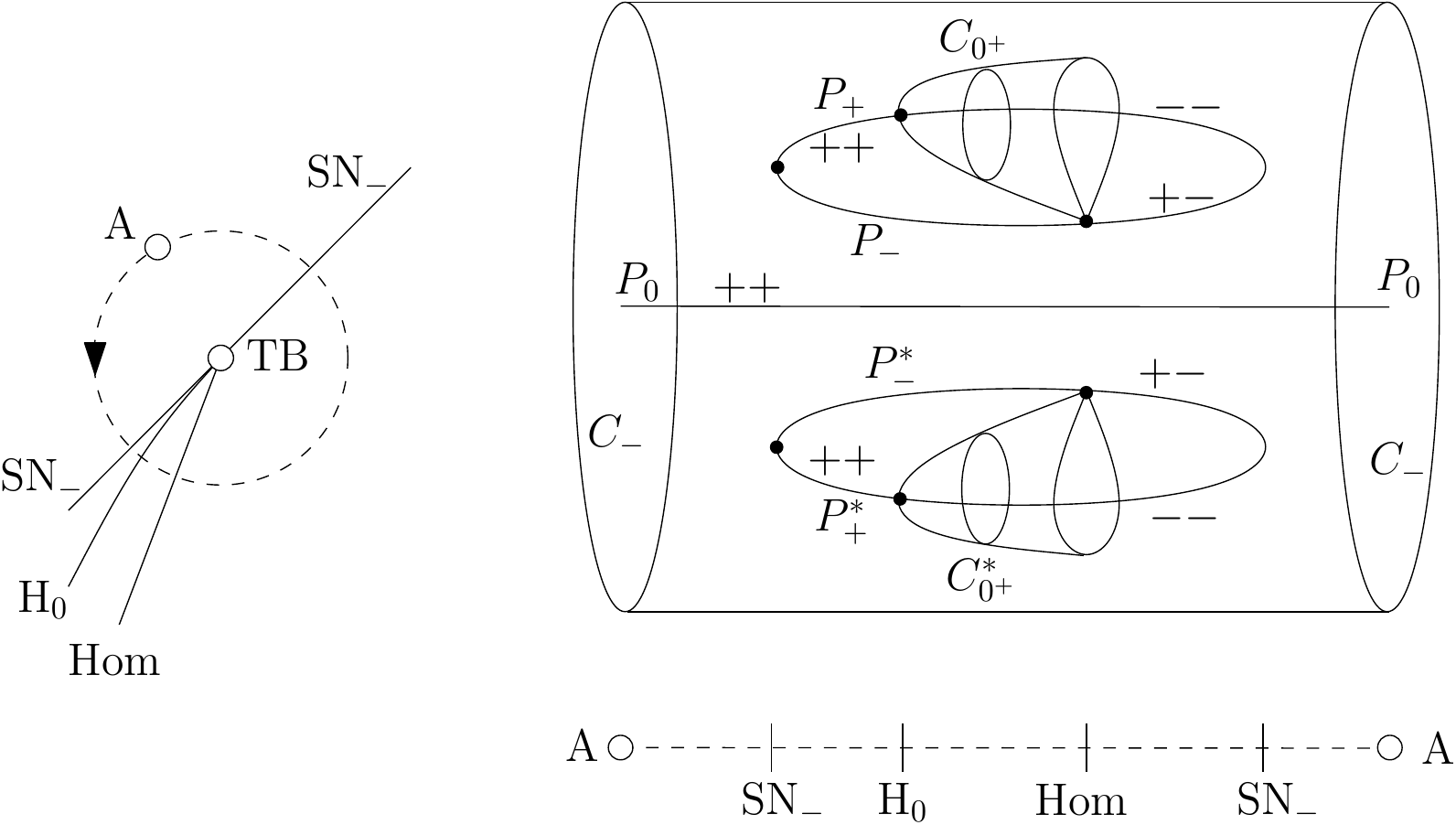}
  \end{center}
  \caption{Double Takens--Bogdanov bifurcation TB: left, bifurcation
    curves emanating from TB in parameter space, with a closed
    one-dimensional path; right, schematics of the bifurcations
    along the closed path, starting and ending at A.}
  \label{TB_double}
\end{figure}

A generic Takens--Bogdanov bifurcation (without symmetry) takes place
at the TB point on the SN$_-$ curve. At the same point in parameter
space, but separate in phase space, two $Z_2$ symmetrically related
Takens--Bogdanov bifurcations take place, with $P_+$ and $P^*_+$ being
the bifurcating states. A schematic of the bifurcations along a
one-dimensional path in parameter space around the TB point is shown
in figure~\ref{TB_double}. Apart from the states
locally connected to both TB bifurcation, there also exist the base
state $P_0$ and the limit cycle $C_-$.

\subsection{Global bifurcations}\label{glob_bif_subsect}

\begin{figure}
  \begin{center}
    \includegraphics[width=0.6\linewidth]{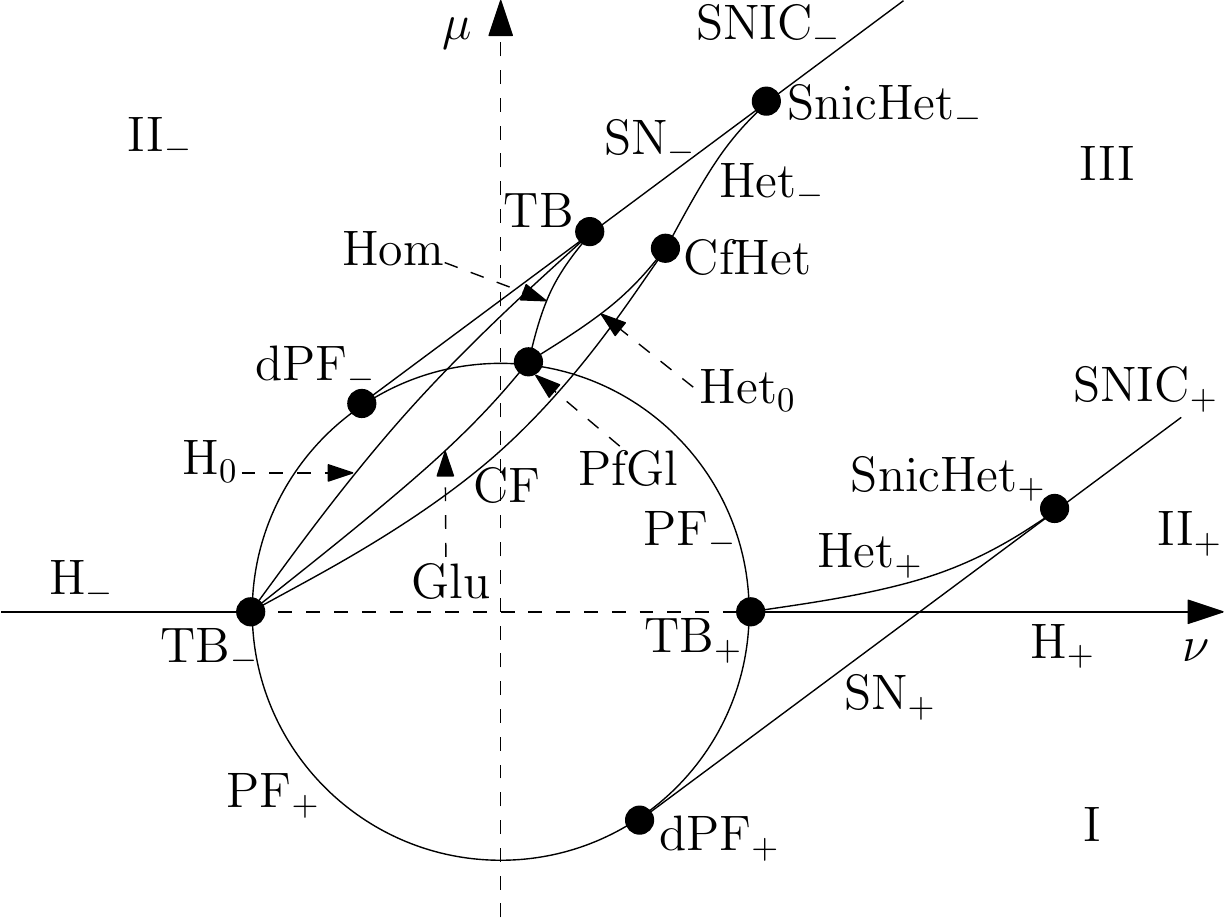}
  \end{center}
  \caption{Global bifurcations in the symmetry breaking of $SO(2)$ to
    $Z_2$ case. Codimension-one bifurcation curves: Gluing Glu,
    cyclic-fold CF, homoclinic collision Hom and heteroclinic loops
    Het$_0$, Het$_\pm$. Codimension-two bifurcation points:
    pitchfork-gluing bifurcation PfGl, cyclic-fold heteroclinic
    bifurcation CfHet and two SNIC-heteroclinic bifurcations
    SnicHet$_\pm$.}
  \label{impHopf_regions}
\end{figure}

In the analysis of the local bifurcations of fixed points we have
found three curves of global bifurcations, a gluing curve Glu and a
saddle-node of cycles CF emerging from TB$_-$, a
heteroclinic loop born at TB$_+$ (Het$_+$), and a homoclinic loop emerging
from TB (Hom). One wonders about the fate of these global bifurcation
curves, and about possible additional global bifurcations. Numerical
simulations of the solutions of the normal form ODE system
\eqref{complexNF}, or equivalently \eqref{cartesianNF}, together with
dynamical systems theory considerations have been used to answer these
questions, and a schematic of all local and global bifurcation curves
is shown in figure~\ref{impHopf_regions}.

\begin{figure}\setlength{\Figsize}{0.4\linewidth}
\begin{center}\footnotesize
\includegraphics[width=\linewidth]{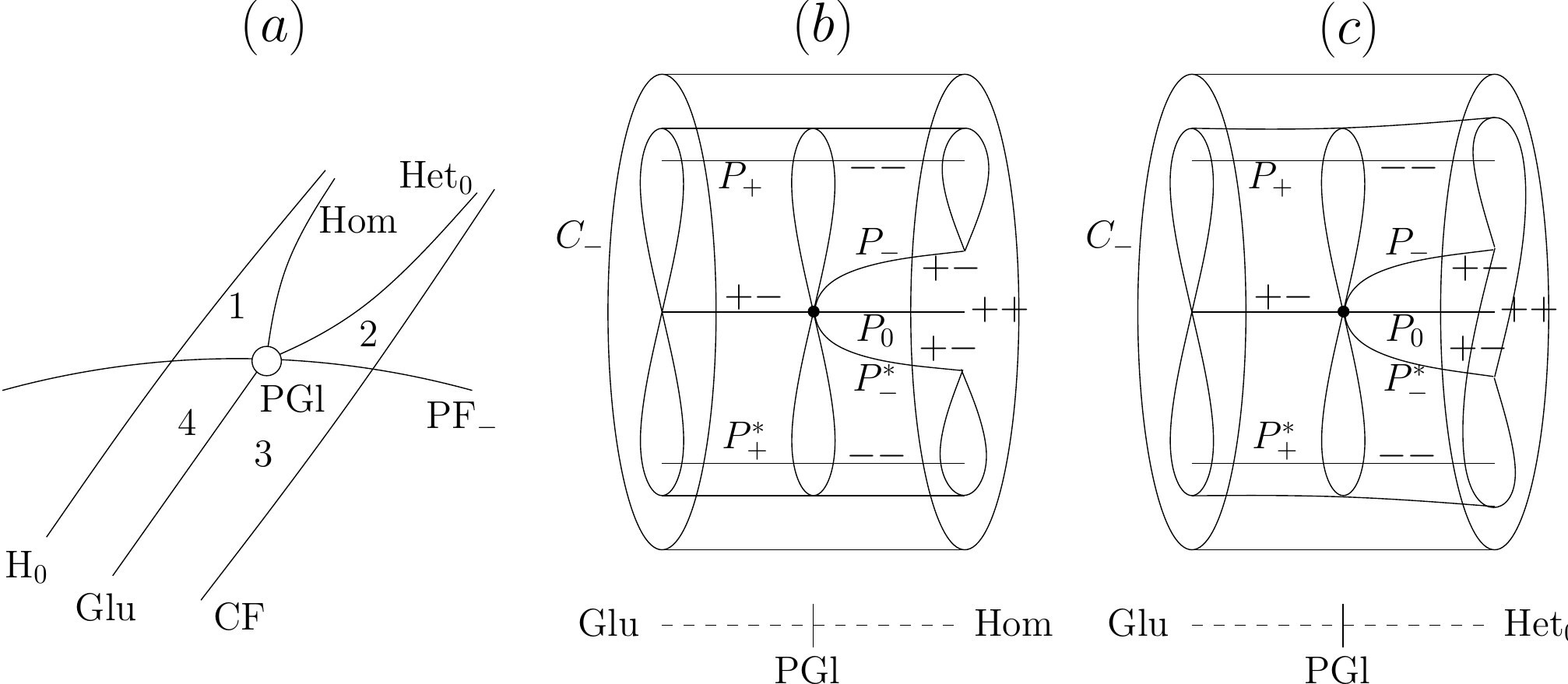} \\[20pt]
\begin{tabular}{@{}l@{\hspace{2pt}}c@{\hspace{5pt}}c@{\hspace{5pt}}c@{}}
$(d)$ 
& 1: $(\mu,\nu)=(0.94,-0.35)$ & 2: $(\mu,\nu)=(0.987,-0.2985)$\\
& \includegraphics[width=\Figsize]{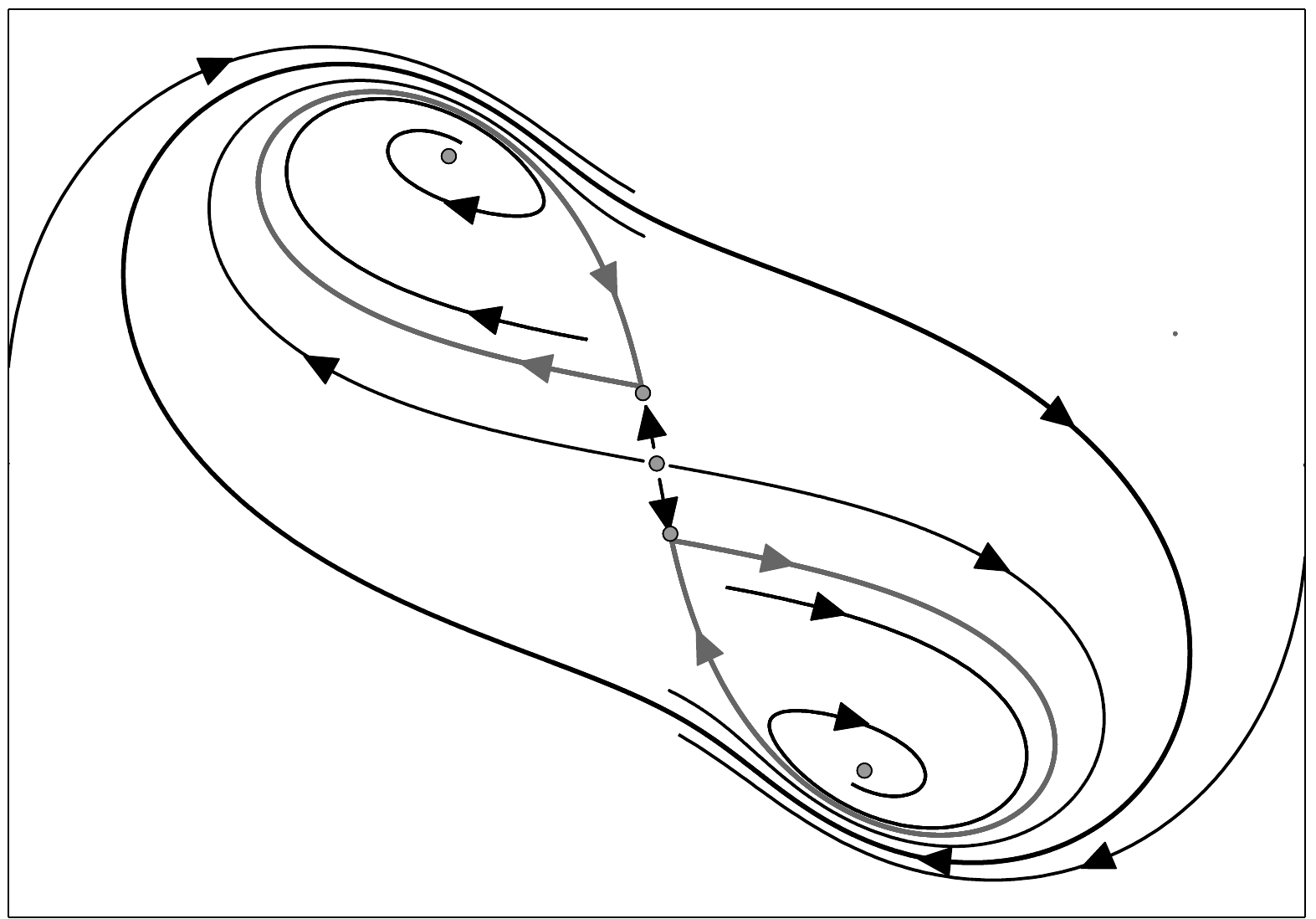} &
  \includegraphics[width=\Figsize]{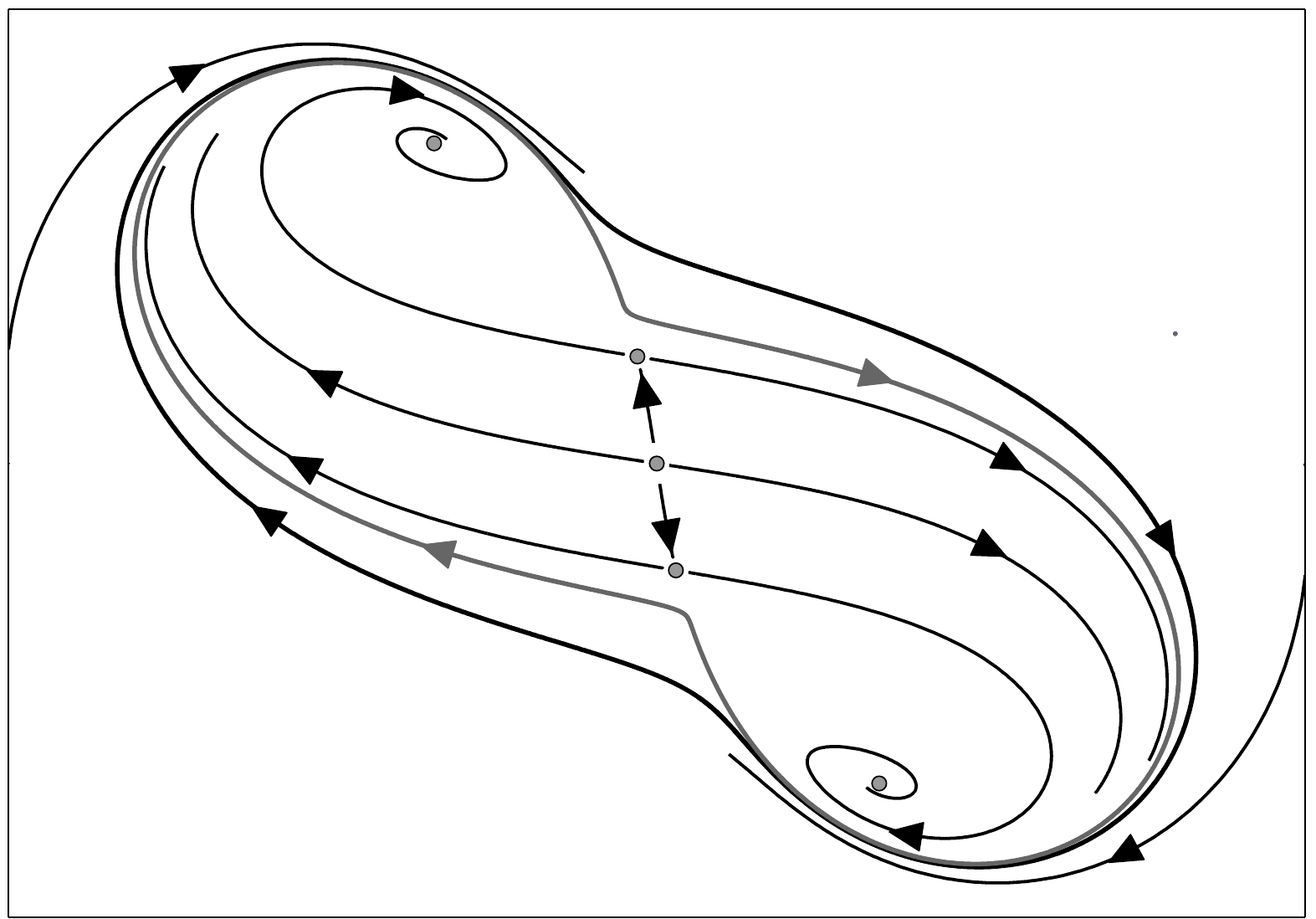} \\
& 3: $(\mu,\nu)=(0.936,-0.34)$ & 4: $(\mu,\nu)=(0.936,-0.352)$\\
& \includegraphics[width=\Figsize]{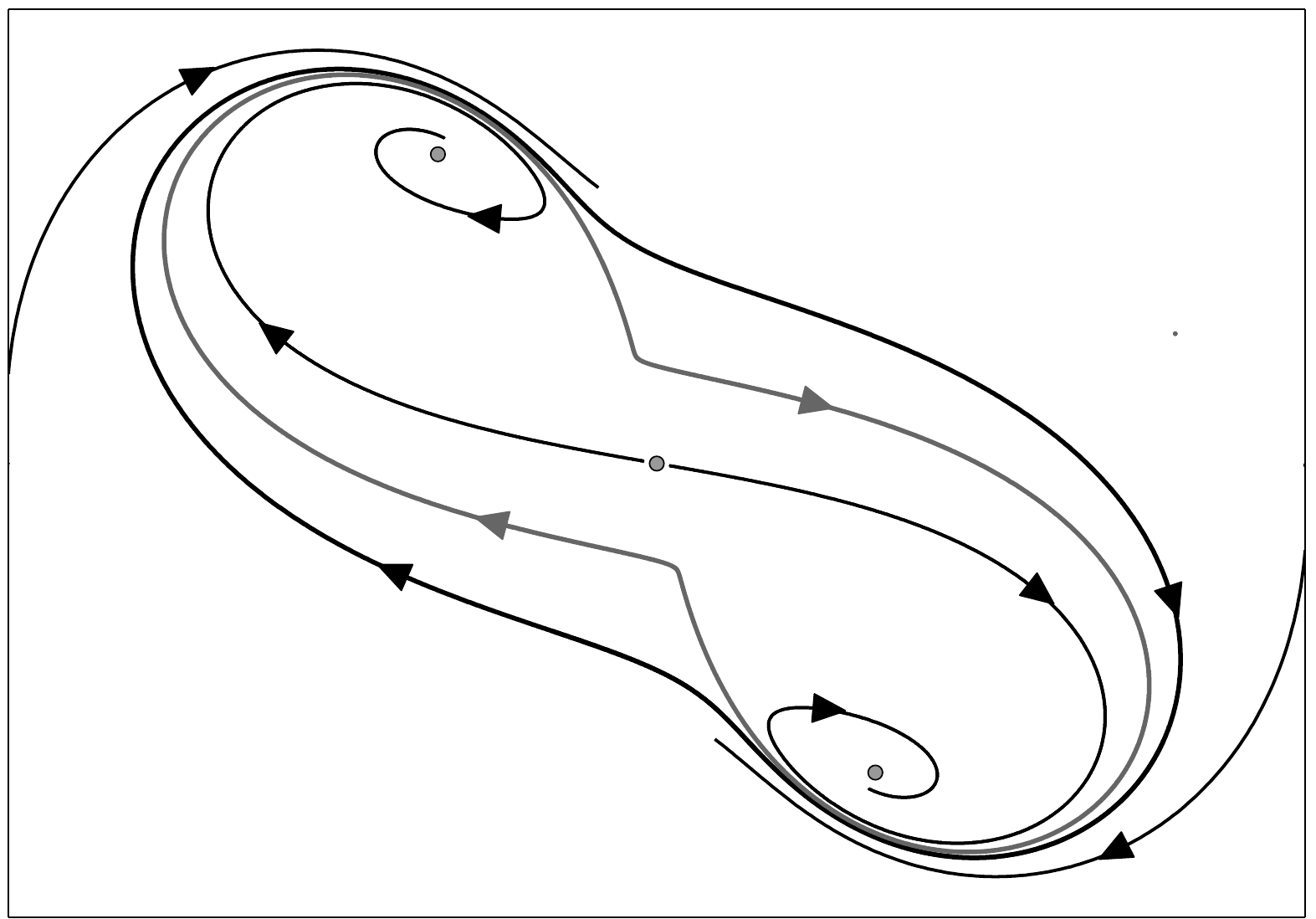} &
  \includegraphics[width=\Figsize]{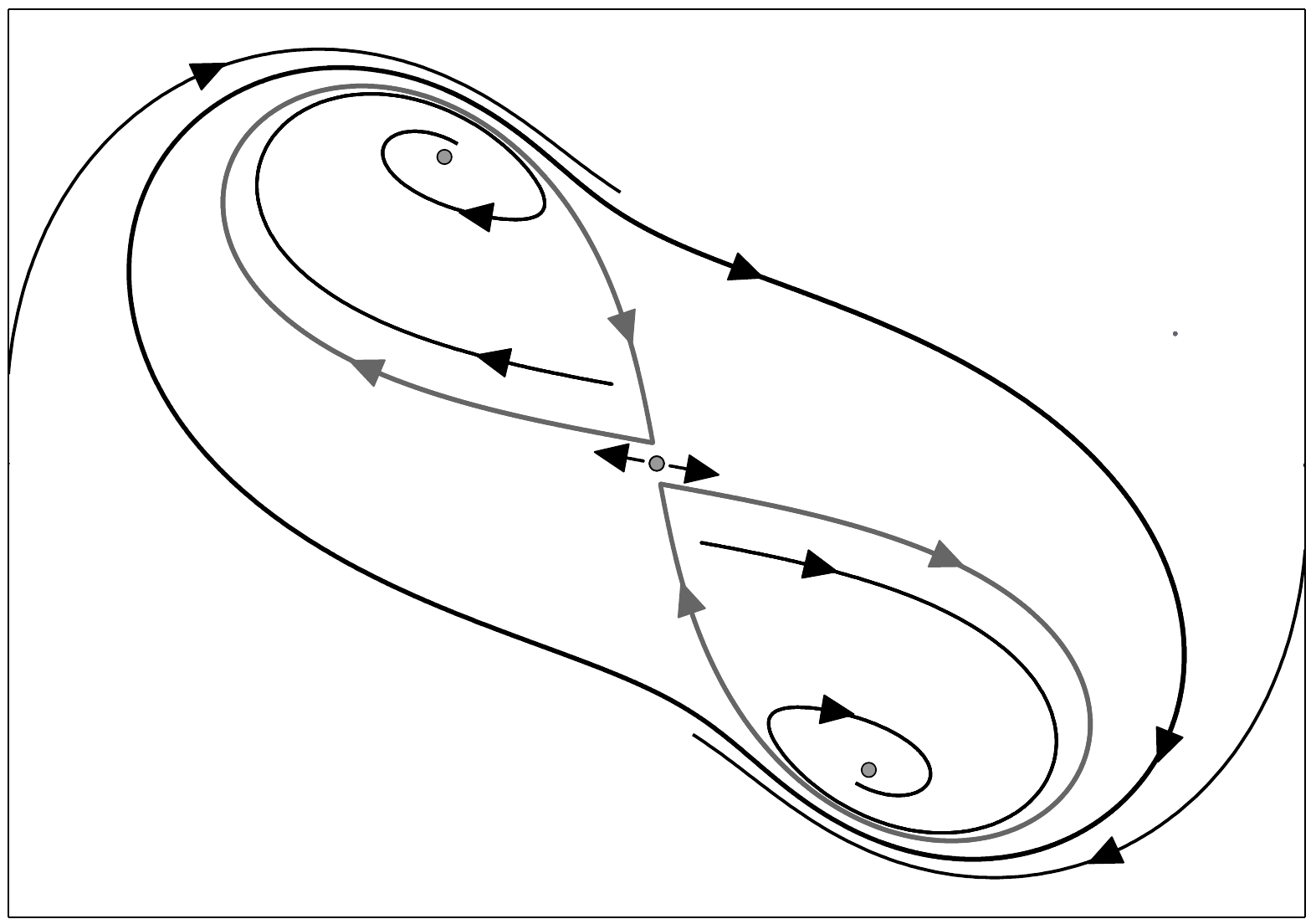} \\
\end{tabular}\vspace*{-6pt}
\end{center}
\caption{Pitchfork-gluing bifurcation PfGl. Top left, bifurcation
  curves emanating from PfGl in parameter space; top right, schematics
  of the bifurcations along the line Glu--Hom. The bottom row shows
  four numerically computed phase portraits, at points labeled 1,
  2, 3 and 4, for the specified parameter values.}
\label{PfGl}
\end{figure}

The gluing bifurcation Glu born at TB$_-$ and the two homoclinic loops
emerging from the two Takens--Bogdanov bifurcations TB (bifurcations
of the symmetric fixed points $P_+$ and $P^*_+$) meet at the point PfGl
on the circle PF$_-$, where the base state $P_0$ undergoes a pitchfork
bifurcation (see figure~\ref{PfGl}$a$). At that point, the two
homoclinic loops of the gluing bifurcation, both homoclinic at the
same point on the stable $P_0$ branch, split when the two fixed points
$P_-$ and $P^*_-$ bifurcate from $P_0$ (see figure~\ref{PfGl}$b$). The
two homoclinic loops are then attached to the bifurcated points and
separate along the curve Hom. The large unstable limit cycle $C_0$,
after the pitchfork bifurcation PF$_-$, collides simultaneously with
both $P_-$ and $P^*_-$, forming a heteroclinic loop along the curve
Het$_0$ (see figure~\ref{PfGl}$c$ and $d$2). Both curves Hom and
Het$_0$ are born at PfGl and separate, leaving a region in between
where none of the cycles $C_{0^+}$, $C^*_{0^+}$ and $C_0$ exist.  The
unstable periodic solution $C_0$ merges with the stable periodic
solution $C_-$ that was born in H$_-$ and existed in region I,
resulting in a cyclic-fold bifurcation of periodic solutions CF. Phase
portraits around PfGl are shown in figure~\ref{PfGl}$(d)$.

The curve Hom born at PfGl ends at the double Takens--Bogdanov point
TB. Locally, around both Takens--Bogdanov bifurcations at TB, after
crossing the homoclinic curve the limit cycles $C_{0^+}$ and
$C^*_{0^+}$ disappear, and no cycles remain. The formation of a large
cycle $C_0$ at Het$_0$ surrounding both fixed points $P_+$ and $P^*_+$ is
a global bifurcation involving simultaneously both $P_-$ and $P^*_-$
unstable points. It is the re-injection induced by the presence of the
$Z_2$ symmetry that is responsible for this global phenomenon
\citep{Adl46,Adl73,Str94,Kuz04}. The two global bifurcation curves
Het$_0$ and CF become very close when leaving the PfGl neighborhood, and
merge at some point in a CfHom (Cyclic-Fold--Heteroclinic collision)
codimension-two global bifurcation. After CfHom, the stable limit cycle
$C_-$, instead of undergoing a cyclic-fold bifurcation, directly
collides with the saddle points $P_-$ and $P^*_-$ along the Het$_-$
bifurcation curve (see figure~\ref{impHopf_regions}). In fact Het$_0$
and Het$_-$ are collisions of a limit cycle with $P_-$ and $P^*_-$, but
the limit cycle is on a different branch of the saddle-node of cycles
CF on each side of CfHom. The two limit cycles born at CF are extremelly
close together in the neighborhood of CfHom, and it is impossible to
see them in a phase portrait, except with a very large zoom around $P_-$
or $P^*_-$.

\begin{figure}
\begin{center}
\begin{tabular}{cc}
 $(a)$ $\mu = 1.6$, $\nu = 0.19$ & $(b)$ $\mu = 1.6$, $\nu = 0.213$ \\
\includegraphics[width=0.45\linewidth]{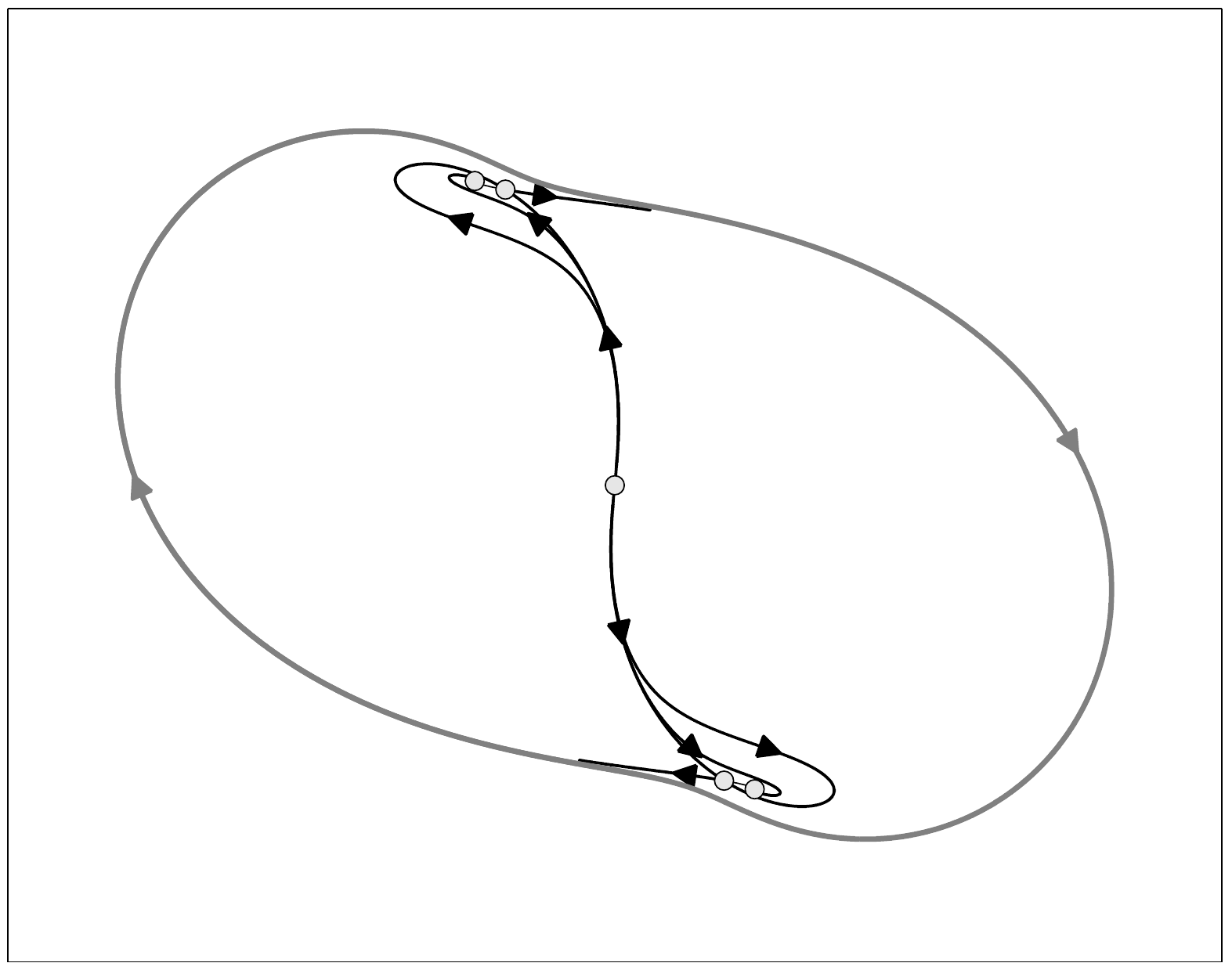} &
\includegraphics[width=0.45\linewidth]{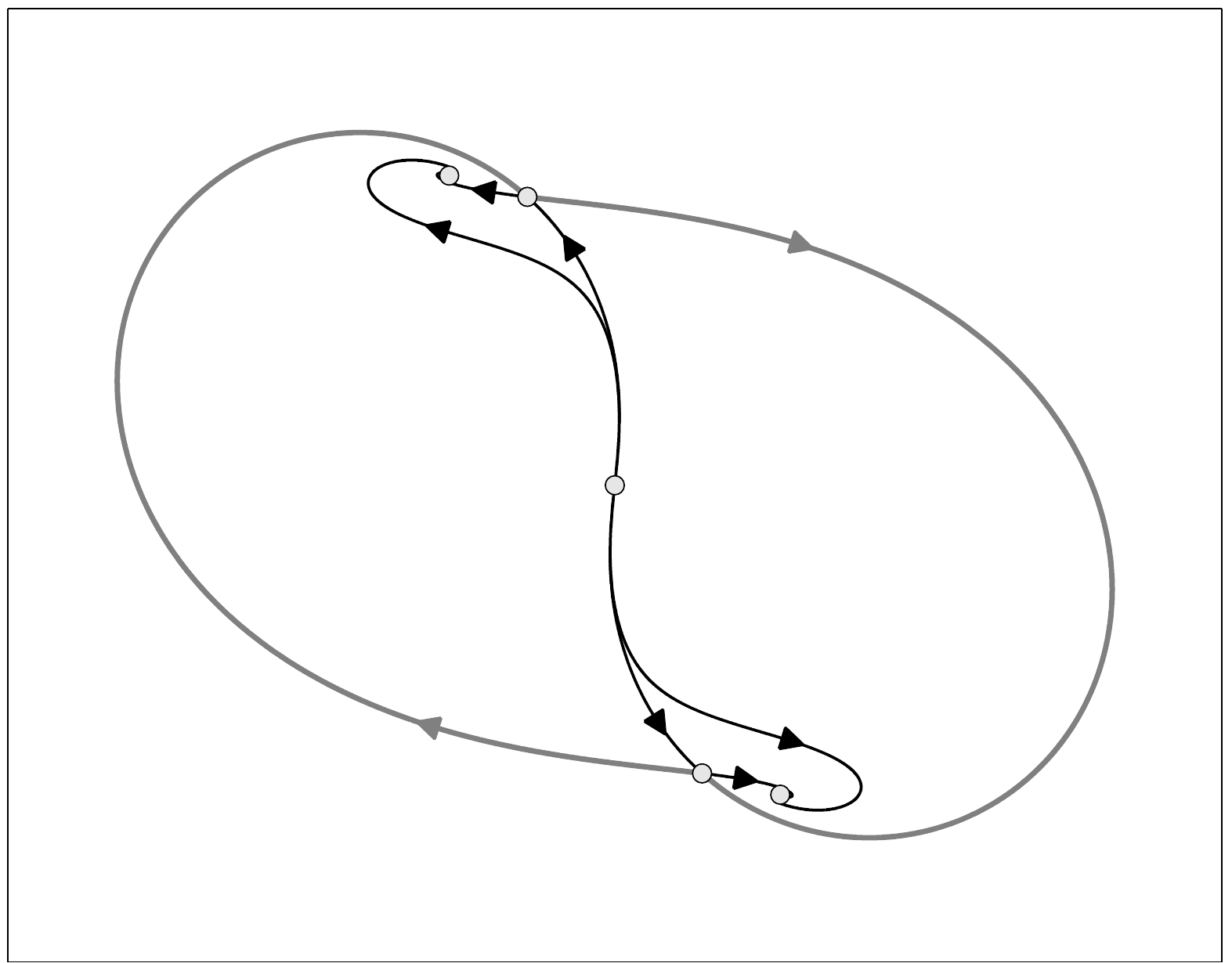} \\
 $(c)$ $\mu = 3.0$, $\nu = 0.19$ & $(d)$ $\mu = 3.0$, $\nu =  1.5858$ \\
\includegraphics[width=0.45\linewidth]{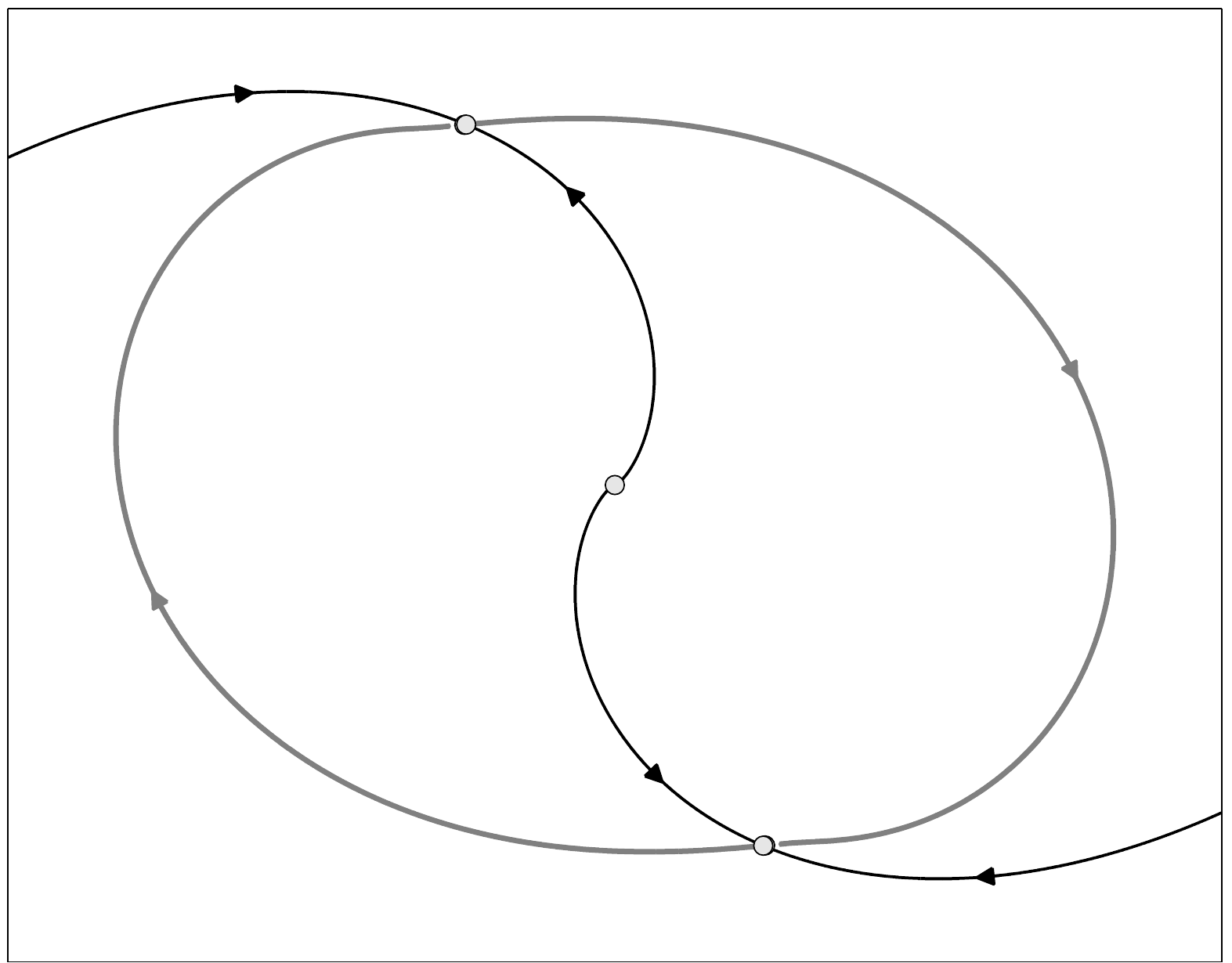} &
\includegraphics[width=0.45\linewidth]{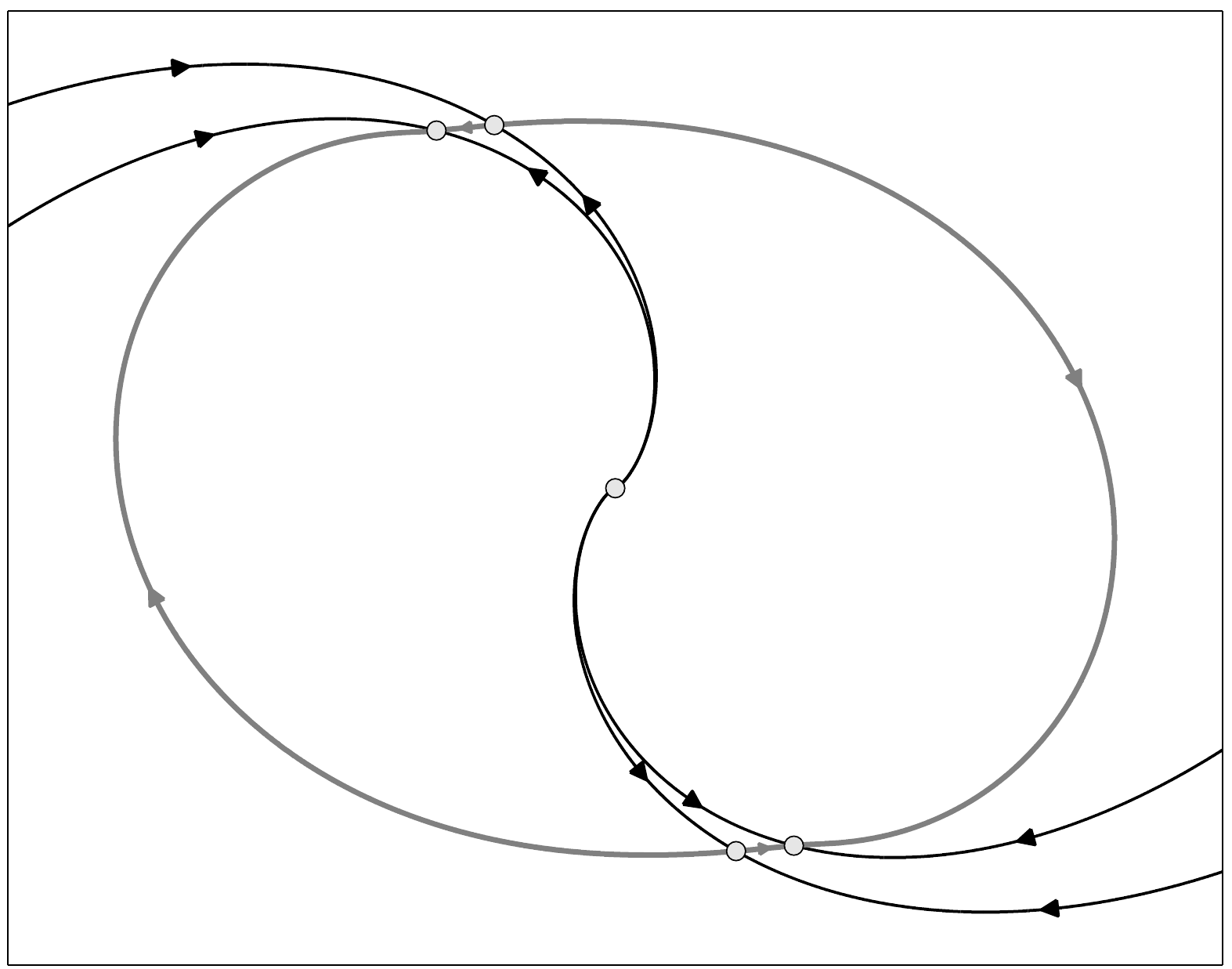}
\end{tabular}
\end{center}
\caption{Numerically computed phase portraits in the $\epsilon\bar z$
  case, for $\alpha_0=45^\circ$ and $\mu$ and $\nu$ as indicated;
  cases $(a)$ and $(b)$ are below the SnicHet$_-$ point and cases $(c)$
  and $(d)$ are above the SnicHet$_-$ point.}
\label{SNIC_Het}
\end{figure}

When increasing $\mu^2+\nu^2$, the heteroclinic loop Het$_-$ born at CfHet
intersects the SN$_-$ curve at a codimension-two global bifurcation point
SnicHet'. When they intersect, the saddle-node appears precisely on
the limit cycle, resulting in a SNIC bifurcation (a saddle-node on an
invariant circle bifurcation). At the SnicHet$_-$ point, the saddle-node
and homoclinic bifurcation curves become tangents, and the saddle-node
curve becomes a SNIC bifurcation curve for larger values of
$\mu^2+\nu^2$. Figure~\ref{SNIC_Het} shows numerically computed phase
portraits for $\alpha_0=45^\circ$, below and above the SnicHet$_-$,
located at $\nu\approx 0.1290$, $\mu\approx 1.505$. In the first case
(figure~\ref{SNIC_Het}$a$) the saddle-node bifurcation SN$_-$ takes place
in the interior of $C_-$, while in the second case
(figure~\ref{SNIC_Het}$b$) SN$_-$ happens precisely on top of $C_-$,
resulting in a SNIC bifurcation.

There remains a global bifurcation curve to be analyzed, the
heteroclinic loop Het born at TB$+$. As shown in
figure~\ref{impHopf_regions}, the curve Het intersects the SN$_+$ curve
tangentially at a codimension-two global bifurcation point
SnicHet. Beyond this point, the SN$_+$ curve becomes a line of SNIC
bifurcations, where the double saddle-node bifurcations appear on the
stable limit cycle $C_+$, which disappears on entering the pinning
region III, exactly in the same way as has been discussed for the
SnicHet' bifurcation.

We can estimate the width of the pinning region as a function of the
magnitude of the imperfection $\epsilon$ and the distance $d$ to the
bifurcation point, $w(d,\epsilon)$. The distance $d$ will be measured
along the line L, and the width $w(d)$ will be the width of the
pinning region measured transversally to L at a distance $d$ from the
origin.  It is convenient to use the coordinates $(u,v)$, so that the
parameter $u$ along L is precisely the distance $d$.  If the pinning
region is delimited by a curve of equation $v=\pm h(u)$, then
$w=2h(u)=2h(d)$. With an imperfection of the form $\epsilon\bar z$,
the case analyzed in this section, the pinning region is of constant
width $w=2$. By restoring the dependence on $\epsilon$, we obtain a
width of value $w(d,\epsilon)=2\epsilon$, independent of the distance
to the bifurcation point; the width of the pinning region is
proportional to $\epsilon$.

\section{Symmetry breaking of $SO(2)$ with quadratics terms}
\label{Sec_quadratic}

\subsection{The $\epsilon\bar z^2$ case}\label{Sec_Z3}

\begin{figure}
\begin{center}
\begin{tabular}{@{}c@{\hspace{4pt}}c@{}}
 $(a)$ & $(b)$  \\
 \includegraphics[width=0.48\linewidth]{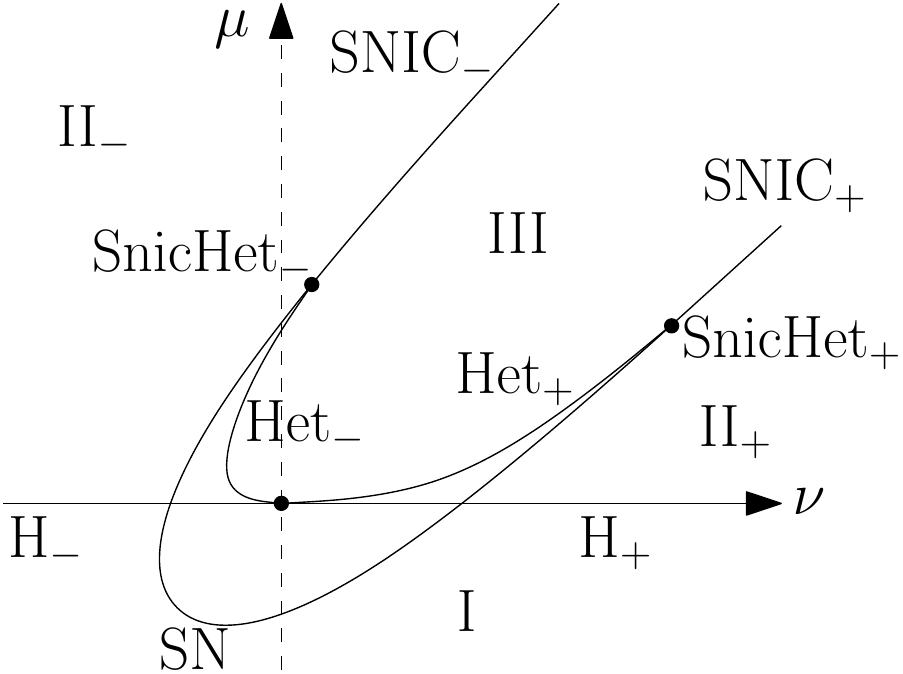} &
 \includegraphics[width=0.48\linewidth]{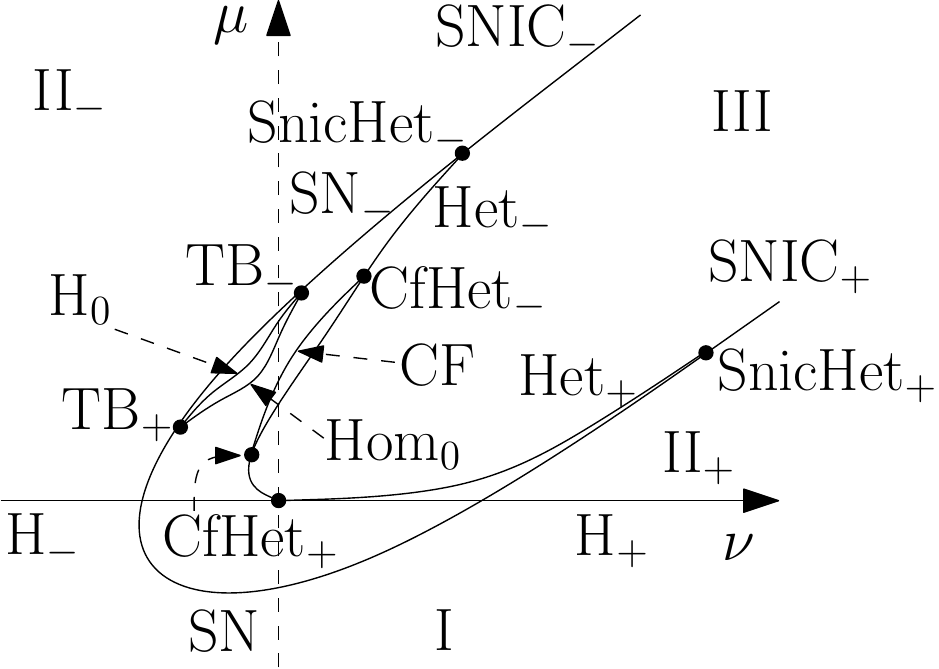} 
\end{tabular}
\end{center}
\caption{Bifurcation curves corresponding to the normal forms with
  quadratic terms \eqref{complexNFbarzz} in the $\epsilon\bar z^2$
  case. $(a)$ $\alpha_0>30^\circ$, $(b)$ $\alpha_0<30^\circ$. H$_0$ is
  tangent to the parabola at the Takens bogdanov points TB$_\pm$, and
  H$_0$, Hom$_0$ almost coincide with SN$_-$; in the figure the distances
  have been exagerated for clarity.}
\label{quadratic_barzz}
\end{figure}

The normal form to be analyzed in this case is
\begin{equation}\label{complexNFbarzz}
 \dot z=z(\mu+i\nu-c|z|^2)+\bar z^2,
\end{equation}
or in terms of the modulus and phase $z=r\ce^{\ci\phi}$,
\begin{equation}\label{realNFbarzz}
\begin{aligned}
 & \dot r=r(\mu-ar^2)+r^2\cos3\phi,\\
 & \dot\phi=\nu-br^2-r\sin3\phi.
\end{aligned}
\end{equation}
The fixed points are the origin $P_0$ ($r=0$) and the solutions of the
same biquadratic equation as in the two previous cases. However, there
is an important difference: due to the factor 3 inside the
trigonometric functions in \eqref{realNFbarzz}, the $P^i_\pm$ points
come in triplets ($i=1,2,3$), each triplet has the same radius $r$ but
their phases differ by 120$^\circ$. This is a consequence of the
invariance of the governing equation \eqref{complexNFbarzz} to the
$Z_3$ symmetry group generated by rotations of 120$^\circ$ around the
origin. This invariance was not present in the two previous cases. The
bifurcation curves of the fixed points (excluding the Hopf
bifurcations of $P^i_\pm$) are still given by
figure~\ref{quadratic_bifs}$(a)$, but now a triplet of symmetric
saddle-node bifurcations take place simultaneously on the SN$_+$ and SN$_-$
curves.

It can be seen (details in \ref{Appendix_quadratic}) that $P^i_-$ are
saddles in the whole of region III, while $P^i_+$ are stable, except
for small angles $\alpha_0<\pi/6$ in a narrow region close to SN$_-$
where they are unstable. For $\alpha_0>\pi/6$, the bifurcation diagram
is exactly the same as in the $z\bar z$ case
(figure~\ref{quadratic_zbarz}$a$), except that the homoclinic curves are
now heteroclinic cycles between the triplets of saddles $P^i_-$; this
case is illustrated in figure~\ref{quadratic_barzz}$(a)$. For
$\alpha_0<\pi/6$, two Takens--Bogdanov bifurcation points appear at
the tangency points between the SN$_-$ curve and the arc H$_0$ of the
ellipse $(b\mu-2a\nu)^2+(a\mu-1)^2=1$, as shown in
figure~\ref{quadratic_barzz}$(b)$. The arc H$_0$ is a Hopf bifurcation
curve of $P^i_+$: three unstable limit cycles $C^i_0$ are born when
the triplet $P^i_+$ becomes stable. These unstable limit cycles
disappear upon colliding with the saddles $P^i_-$ on a curve of
homoclinic collisions, Hom$_0$, that ends at the two Takens--Bogdanov
points TB$_+$ and TB$_-$. This situation is very similar to what
happens in the $\epsilon\bar z$ case analyzed in section
\S\ref{Sec_bar_z}, where a Hopf curve H$_0$ appeared close to SN$_-$
joining two Takens--Bogdanov points. In both cases, the $SO(2)$
symmetry is not completely broken, but a $Z_m$ symmetry remains.  For
$\alpha_0=\pi/6$, the ellipse H$_0$ becomes tangent to SN$_-$ and the two
Takens--Bogdanov points coalesce, disappearing for $\alpha_0>\pi/6$.

Finally, the stable limit cycles $C_-$ and $C_+$, born at the Hopf
bifurcations H$_-$ and H$_+$, upon entering region III collide
simultaneously with the saddles $P^i_+$, $i=1$, 2 and 3, and disappear
along two heteroclinic bifurcation curves Het and Het' for small
values of $\mu$. These curves collide with the parabola at the
codimension-two bifurcation points SnicHet and SnicHet'. For larger
values of $\mu$ the limit cycles $C_-$ and $C_+$ undergo SNIC
bifurcations on the parabola. The two curves Hom$_\pm$ emerge from
the origin, as in the previous cases.

In the three quadratic cases, the pinning region is delimited by the
same parabola $u=v^2-1/4$, using the $(u,v)$ coordinates introduced in
\eqref{uv_munu} (see also figure~\ref{uv_and_epsilon_SN}$(a)$). The width of
the pinning region is easy to compute, and is given by
$w=2v=2\sqrt{u+1/4}\sim2\sqrt{d}$. By using \eqref{epsilon_scaling},
the dependence on $\epsilon$ is restored, resulting in
$w(d,\epsilon)=2\epsilon\sqrt{d}$. The width of the pinning region
increases with the distance $d$ to the bifurcation point, and it is
proportional to the amplitude of the imperfection $\epsilon$.

\subsection{The $\epsilon z\bar z$ case}\label{Sec_zbarz}

The normal form to be analyzed in this case is
\begin{equation}\label{complexNFzbarz}
 \dot z=z(\mu+i\nu-c|z|^2)+z\bar z,
\end{equation}
or in terms of the modulus and phase $z=r\ce^{\ci\phi}$,
\begin{equation}\label{realNFzbarz}
\begin{aligned}
 & \dot r=r(\mu-ar^2)+r^2\cos\phi,\\
 & \dot\phi=\nu-br^2-r\sin\phi.
\end{aligned}
\end{equation}
There are three fixed points: the origin $P_0$ ($r=0$) and the two
solutions $P_\pm$ of the biquadratic equation
$r^4-2(a\mu+b\nu+1/2)r^2+\mu^2+\nu^2=0$, given by
\begin{equation}
  r^2_\pm=a\mu+b\nu+\frac{1}{2}\pm\sqrt{a\mu+b\nu+\frac{1}{4}-(a\nu-b\mu)^2}.
\end{equation}
These solutions are born at the parabola
$a\mu+b\nu+1/4=(a\nu-b\mu)^2$, and exist only in its interior, which
is the pinning region III in figure~\ref{quadratic_zbarz}. The
parabola is a curve of saddle-node bifurcations. It can be seen
(details in \ref{Appendix_quadratic}) that $P_+$ is stable while $P_-$
is a saddle in the whole of region III, so there are no additional
bifurcations of fixed points in the $z\bar z$ case, in contrast to the
$z^2$ case.

\begin{figure}\setlength{\Figsize}{0.28\linewidth}
\begin{center}
\begin{tabular}{m{15pt}m{0.48\linewidth}}
 $(a)$ & \includegraphics[width=\linewidth]{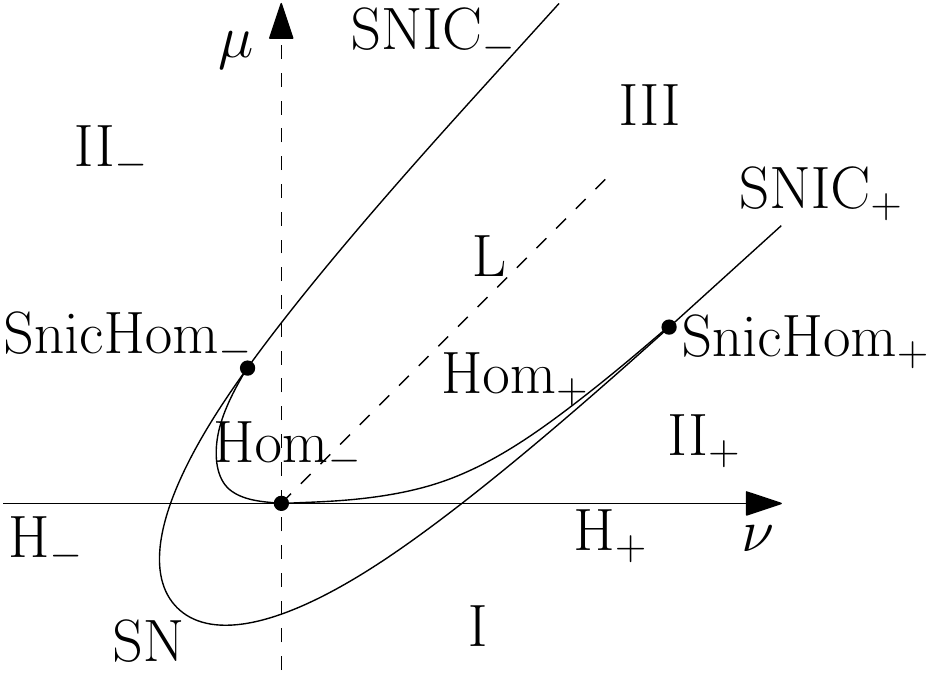}
\end{tabular}\bigskip

\begin{tabular}{@{}m{15pt}@{\hspace{5pt}}m{\Figsize}@{\hspace{5pt}}m{\Figsize}
   @{\hspace{5pt}}m{\Figsize}@{}c@{}}
 $(b)$ & \includegraphics[width=\linewidth,clip=]{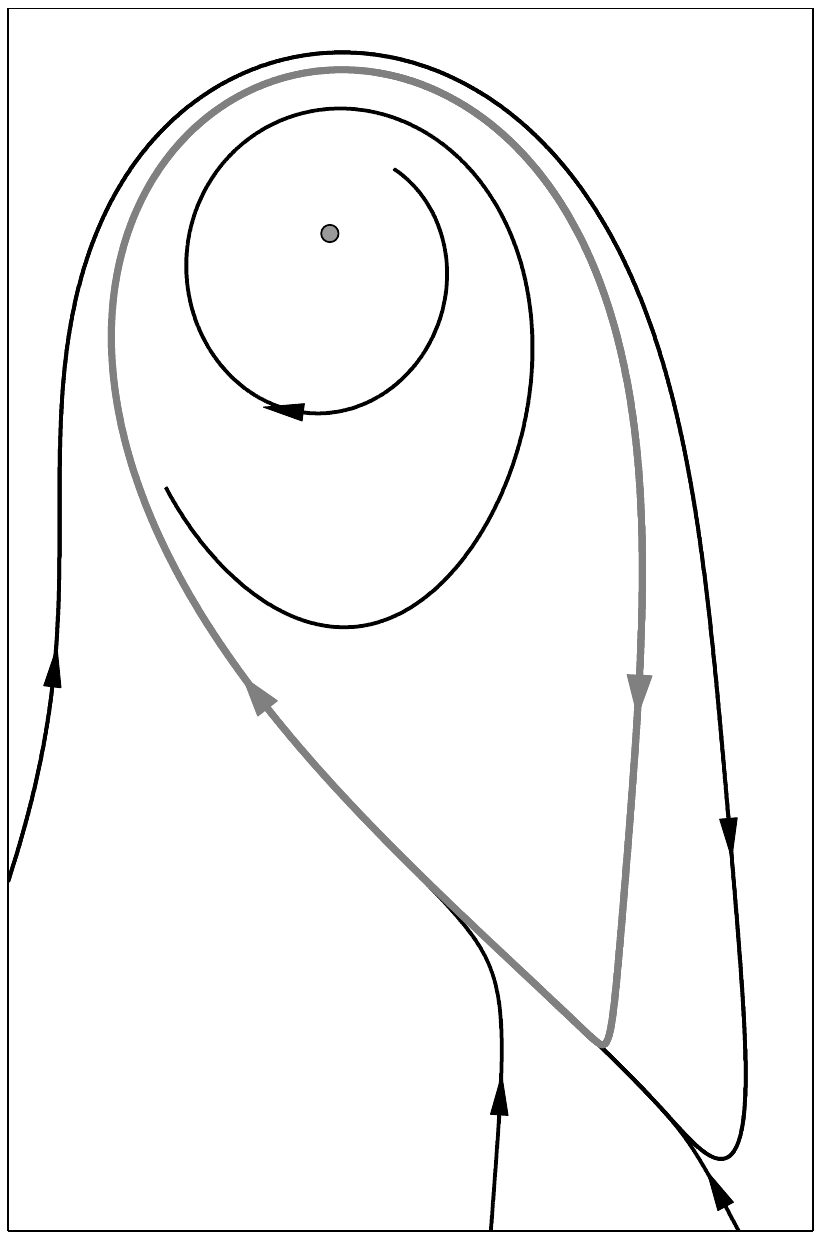} &
   \includegraphics[width=\linewidth,clip=]{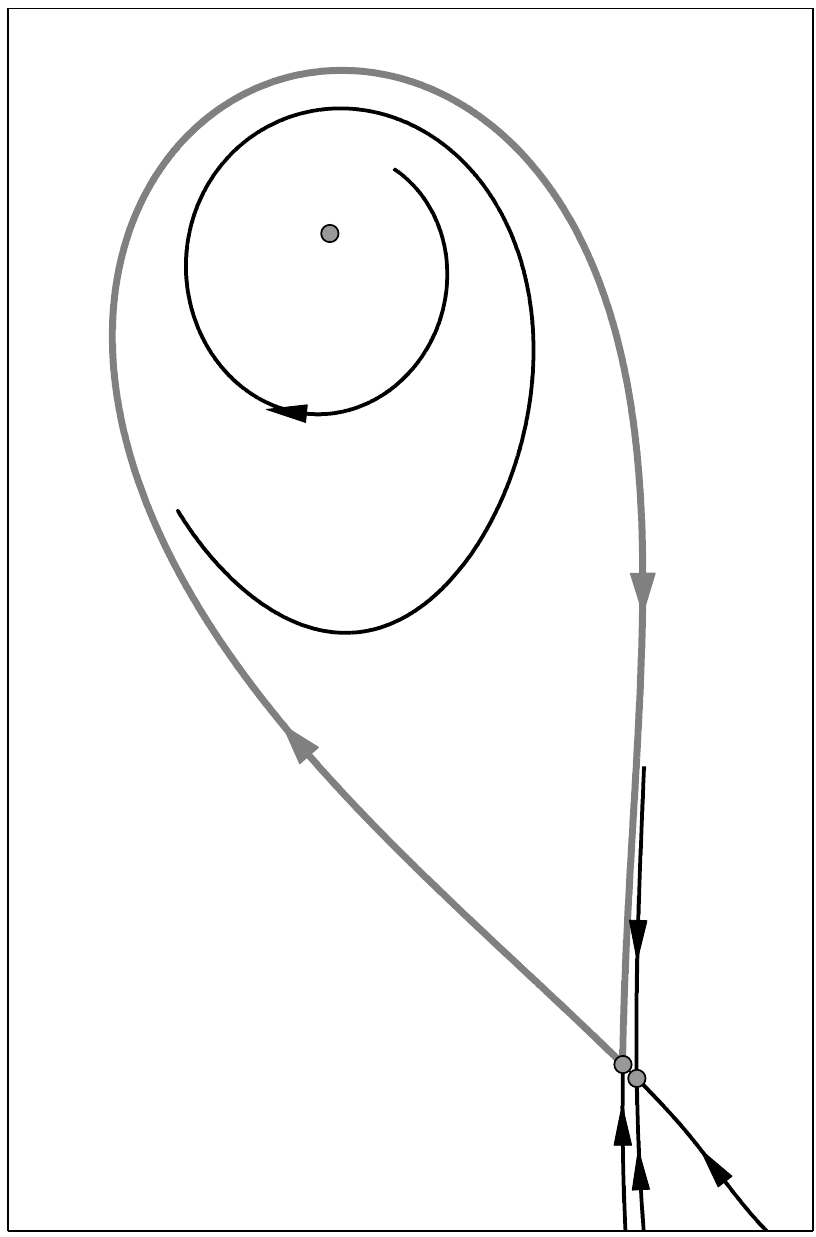} &
   \includegraphics[width=\linewidth,clip=]{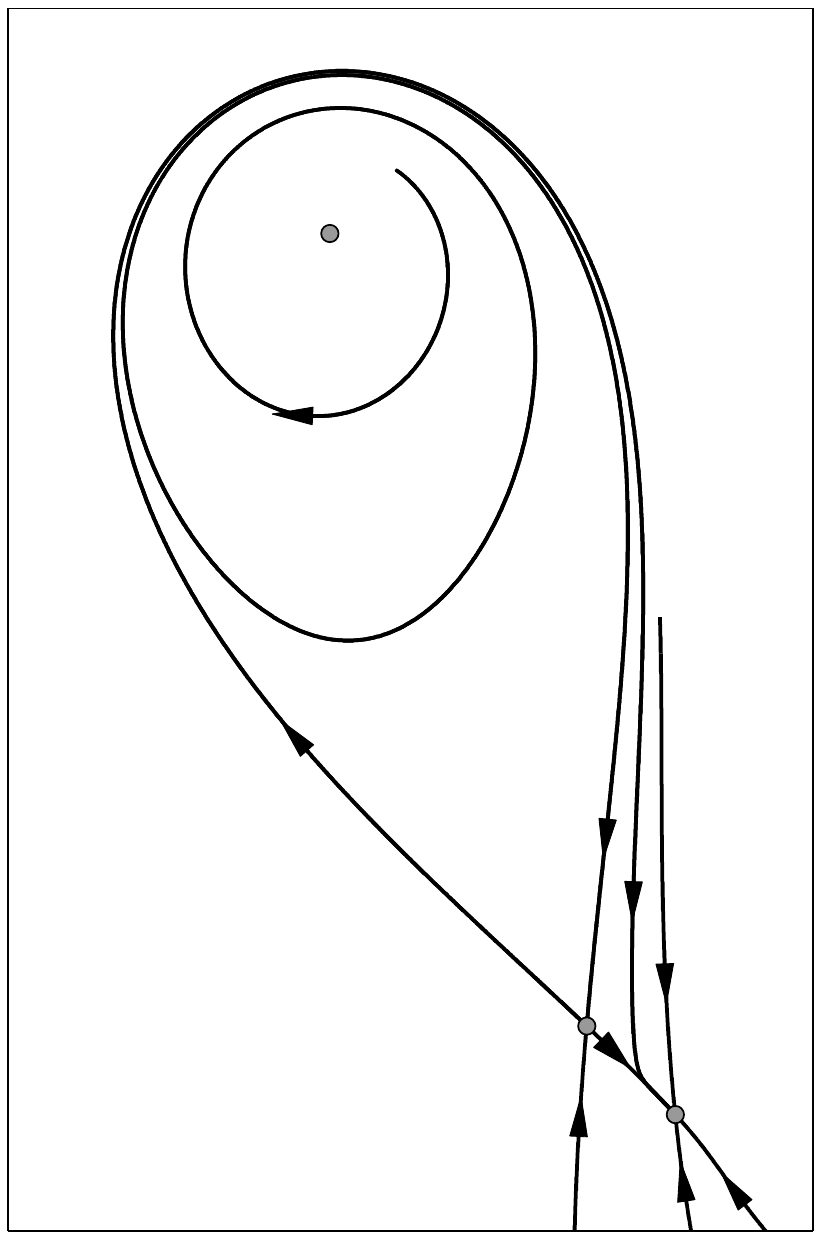} & \\
 & \centering $\nu=-0.306$ & \centering $\nu=-0.30442$ &
   \centering $\nu=-0.303$& \\
 $(c)$ & \includegraphics[width=\linewidth,clip=]{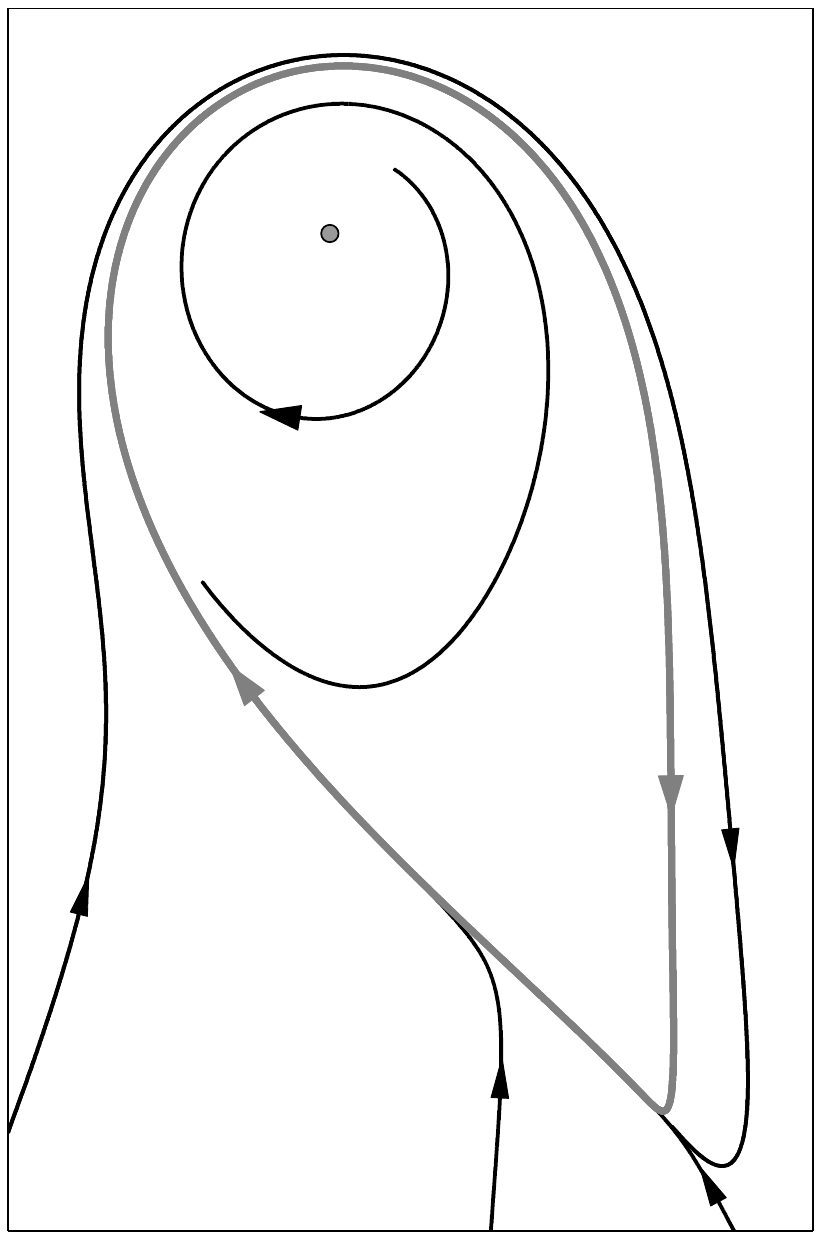} &
   \includegraphics[width=\linewidth,clip=]{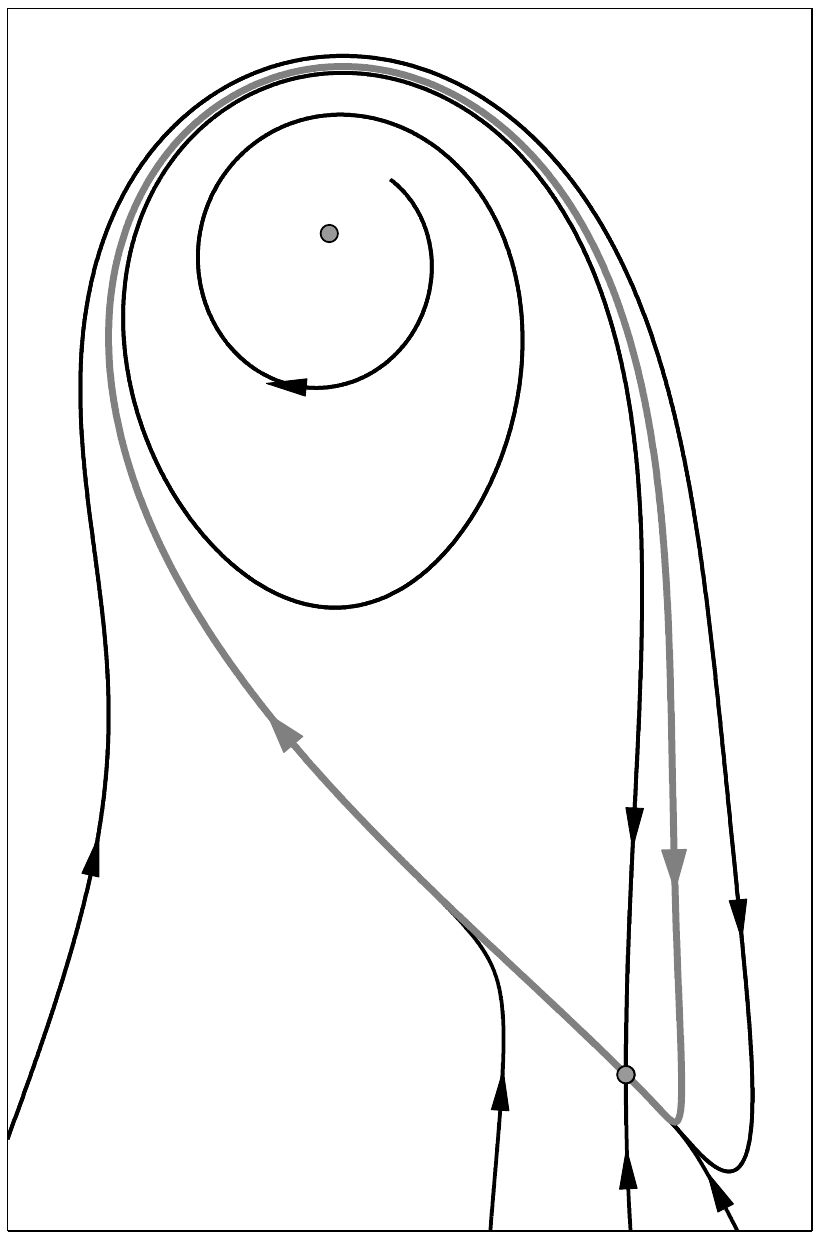} &
   \includegraphics[width=\linewidth,clip=]{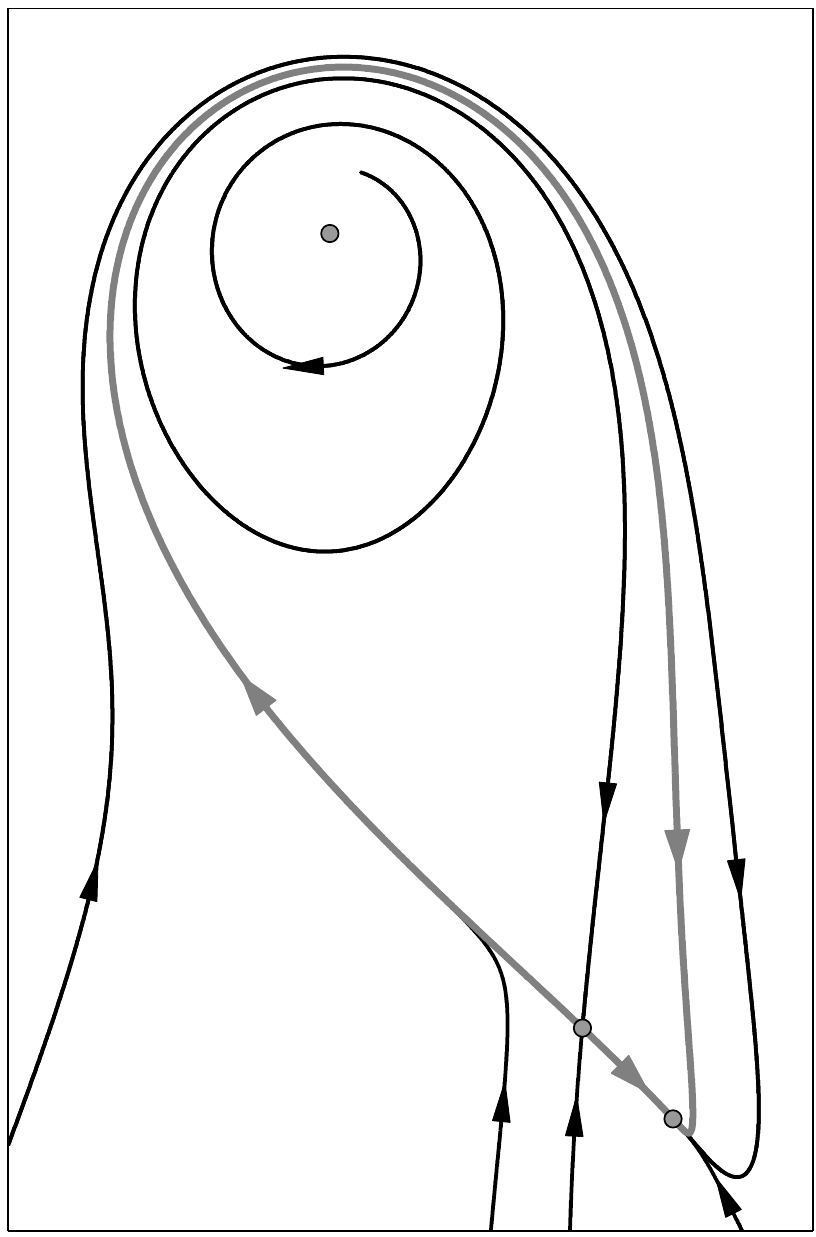} & \\
 & \centering $\nu=-0.307$ & \centering $\nu=-0.30552$ &
   \centering $\nu=-0.304$& 
\end{tabular}
\end{center}
\caption{$(a)$ Schematic of bifurcation curves corresponding to the
  normal form with quadratic terms in the $\epsilon z\bar z$
  case. Phase portraits $(b)$ crossing the Hom$_-$ curve at $\mu=0.03$,
  for $\nu$ values as specified, and $(c)$ crossing the SNIC$_-$ curve at
  $\mu=0.033$. Thick grey lines correspond to the periodic orbit and
  the homoclinic and heteroclinic loops.}
\label{quadratic_zbarz}
\end{figure}

As the perturbation is of second order, the Jacobian at $P_0$ is the
same as in the unperturbed case, and the Hopf bifurcations of $P_0$
take place along the horizontal axis $\mu=0$. As in the unperturbed
case, the Hopf frequency is negative for $\nu<0$ (H$_-$), it is zero
at the origin, and becomes positive for $\nu>0$ (H$_+$). The Hopf
curves H$_-$ and H$_+$ extend in this case up to the origin, in
contrast with the zero and first order cases examined previously,
where the Hopf curves ended in Takens--Bogdanov bifurcations without
reaching the origin.

The stable limit cycles $C_-$ and $C_+$, born at the Hopf bifurcations
H$_-$ and H$_+$, upon entering region III collide with the saddle
$P_+$ and disappear along two homoclinic bifurcation curves Hom$_\pm$
for small values of $\mu$. For larger values of $\mu$, the curves Hom$_\pm$
collide with the parabola at the codimension-two bifurcation
points SnicHom$_\pm$, and for larger values of $\mu$, the
saddle-node bifurcations take place on the parabola and the limit
cycles $C_-$ and $C_+$ undergo SNIC bifurcations, as in the previous
$\epsilon$ and $\epsilon\bar z$
cases. Figure~\ref{quadratic_zbarz}$(a)$ summarizes all the local and
global bifurcation curves in the $z\bar z$ case, and shows numerically
computed phase portraits around the SnicHom$_-$
point. Figure~\ref{quadratic_zbarz}$(b)$ shows the SN$_+$ and Hom$_-$
bifurcations before SnicHom$_-$ (at $\mu=0.03$), and
figure~\ref{quadratic_zbarz}$(c)$ illustrates the SNIC$_-$ bifurcation
after SnicHom$_-$ (at $\mu=0.033$).

\subsection{The $\epsilon z^2$ case}\label{Sec_z2}

\begin{figure}
\begin{center}
\begin{tabular}{cc}
 $(a)$ & $(b)$  \\
\includegraphics[width=0.48\linewidth]{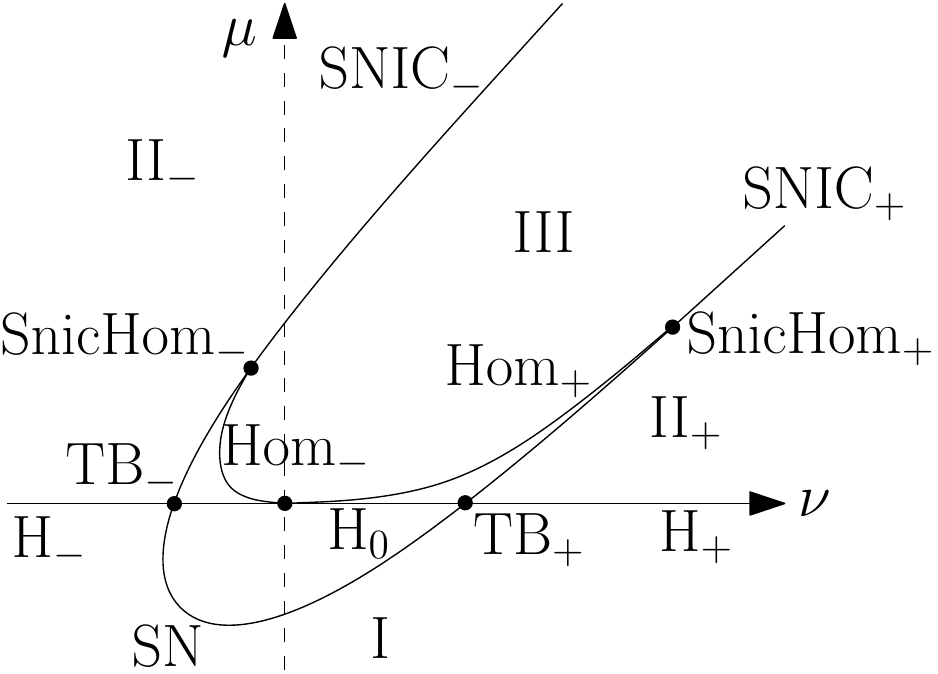} &
\includegraphics[width=0.48\linewidth]{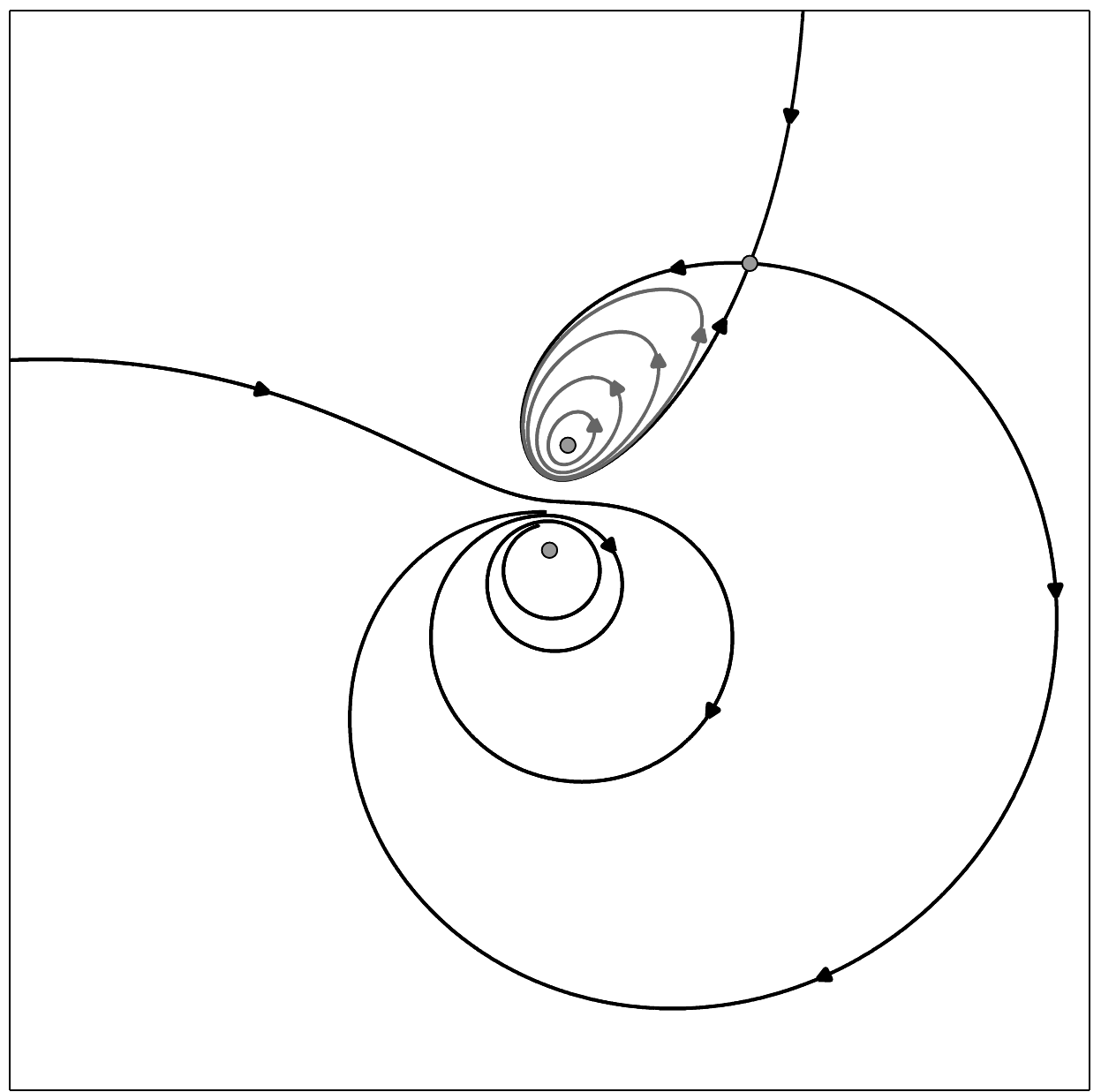}
\end{tabular}
\end{center}
\caption{$(a)$ Bifurcation curves corresponding to the normal form with
  quadratic terms in the $\epsilon z^2$ case. $(b)$ Phase portrait on
  the curve H$_0$ showing the degenerate Hopf bifurcation and the
  associated homoclinic loop.}
\label{quadratic_bifs}
\end{figure}

The normal form to be analyzed in this case is \eqref{NFsinglepert}
with $p=q=2$ and $\epsilon=1$:
\begin{equation}\label{complexNFzz}
 \dot z=z(\mu+i\nu-c|z|^2)+z^2.
\end{equation}
The normal form \eqref{complexNFzz}, in terms of the modulus and phase
$z=r\ce^{\ci\phi}$, reads
\begin{equation}\label{realNFzz}
\begin{aligned}
 & \dot r=r(\mu-ar^2)+r^2\cos\phi,\\
 & \dot\phi=\nu-br^2+r\sin\phi.
\end{aligned}
\end{equation}
There are three fixed points: the origin $P_0$ ($r=0$) and the two
solutions $P_\pm$ of the biquadratic equation
$r^4-2(a\mu+b\nu+1/2)r^2+\mu^2+\nu^2=0$, which are the same as in the
$\epsilon z\bar z$ case (\S\ref{Sec_zbarz}). In fact, the fixed points
in these two cases have the same modulus $r$ and their phases have
opposite sign; changing $\phi\to-\phi$ in \eqref{realNFzbarz} results
in \eqref{realNFzz}. Therefore, the bifurcation curves of the fixed
points (excluding Hopf bifurcations of $P_\pm$) in this case are also
given by figure~\ref{quadratic_zbarz}$(a)$.

It can be seen (details in \ref{Appendix_quadratic}) that $P_+$ is a
saddle in the whole of region III, while $P_-$ is stable for $\mu>0$,
unstable for $\mu<0$, and undergoes a Hopf bifurcation H$_0$ along the
segment of $\mu=0$ delimited by the parabola of saddle-node
bifurcations. The points TB$_\pm$ where H$_0$ meets the parabola are
Takens--Bogdanov codimension-two
bifurcations. Figure~\ref{quadratic_bifs} summarize all the local and
global bifurcation curves just described.

In the present case, the two Takens--Bogdanov bifurcations and the
Hopf bifurcations along H$_0$ are degenerate, as shown in
\ref{Appendix_quadratic}. Detailed analysis and numerical simulations
show that the Hopf and homoclinic bifurcation curves emerging from the
Takens--Bogdanov point are both coincident with the H$_0$ curve
previously mentioned. Moreover, the interior of the homoclinic loop is
filled with periodic orbits, and no limit cycle exists on either side
of H$_0$. This situation is illustrated in the phase portrait in
figure~\ref{quadratic_bifs}$(b)$. This highly degenerate situation
will be broken by the presence of additional terms in the normal form,
and of the continuous family of periodic orbits, only a few will
remain. \citet{drs87}, who have analyzed in detail the unfolding of
such a degenerate case, find that at most two of the periodic orbits
survive.

Finally, and exactly in the same way as in the $z\bar z$ case examined
in the previous subsection, the stable limit cycles $C_-$ and $C_+$
born respectively at the Hopf bifurcations H$_-$ and H$_+$ and
existing in regions I$_-$ and I$_+$ (figure~\ref{quadratic_bifs}$a$),
on entering region III collide with the saddle $P_+$ and disappear
along two homoclinic bifurcation curves Hom$_\pm$ for small values
of $\mu$. These curves collide with the parabola at the
codimension-two bifurcation points SnHom$_\pm$ and for larger values
of $\mu$, the limit cycles $C_-$ and $C_+$ undergo SNIC bifurcations
on the parabola. The two curves Hom$_\pm$ emerge from the origin,
as in the previous case.

\section{Symmetry breaking $SO(2)\to Z_m$, $m\ge4$}\label{Large_Zm_symmetry}

For completeness, and also for intrinsic interest, we will explore the
breaking of the $SO(2)$ symmetry to $Z_m$, so that the imperfections
added to the normal form \eqref{zeroHopfNF} preserve the $Z_m$ subgroup
of $SO(2)$ generated by rotations of $2\pi/m$ (also called $C_m$). The
lowest order monomial in $(z,\bar z)$ not of the form $z|z|^{2p}$ and
equivariant under $Z_m$, is $\bar z^{m-1}$, resulting in the normal
form
\begin{equation}\label{HopfZm}
 \dot z=z(\mu+i\nu-c|z|^2)+\epsilon\bar z^{m-1}.
\end{equation}
In terms of the modulus and phase of the complex amplitude
$z=r\ce^{\ci\phi}$, the normal form becomes
\begin{equation}\label{HopfZm2}
\begin{aligned}
 & \dot r=r(\mu-ar^2)+\epsilon r^{m-1}\cos m\phi,\\
 & \dot\phi=\nu-br^2-\epsilon r^{m-2}\sin m\phi.
\end{aligned}
\end{equation}
The cases $m=2$ and $m=3$ have already been examined in
\S\ref{Sec_bar_z} and \S\ref{Sec_Z3} respectively. When $m\ge 4$, the
term $\epsilon\bar z^{m-1}$ is smaller than the remaining terms in
\eqref{HopfZm}, so the effect of the symmetry breaking is going to be
small compared with the other cases analyzed in this paper. The
fixed point solutions of \eqref{HopfZm2}, apart from the trivial
solution $P_0$ ($r=0$), are very close the zero-frequency line L
in the perfect system. Using the coordinates $(u,v)$ along and
orthogonal to L \eqref{uv_munu}, the nontrivial fixed points of
\eqref{HopfZm2} satisfy
\begin{equation}\label{fix_points_Zm}
  (r^2-u)^2=\epsilon^2r^{2m-4}-v^2.
\end{equation}

\begin{figure}
\begin{center}
\begin{tabular}{c@{\hspace{20pt}}c}
 $(a)$ $m=4$ & $(b)$ $m>4$  \\[5pt]
 \includegraphics[width=0.36\linewidth]{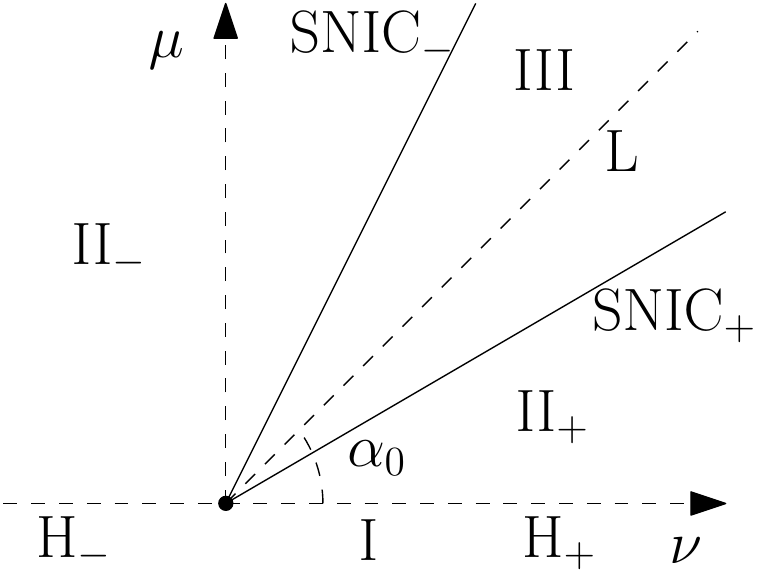} &
 \includegraphics[width=0.36\linewidth]{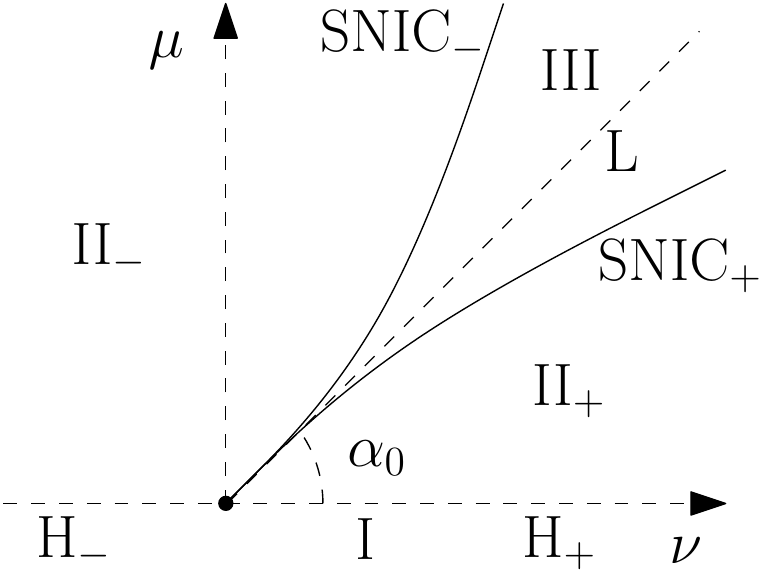} 
\end{tabular}\vspace*{-8pt}
\end{center}
\caption{Bifurcation curves corresponding to the normal forms
  retaining a $Z_m$ symmetry \eqref{HopfZm2}; $(a)$ corresponds to the
  $m=4$ case and $(b)$ corresponds to the $m>4$ cases.}
\label{hornsZm}
\end{figure}

\begin{figure}
\begin{center}
\begin{tabular}{m{12pt}m{0.9\linewidth}}
 $(a)$ & \includegraphics[width=\linewidth]{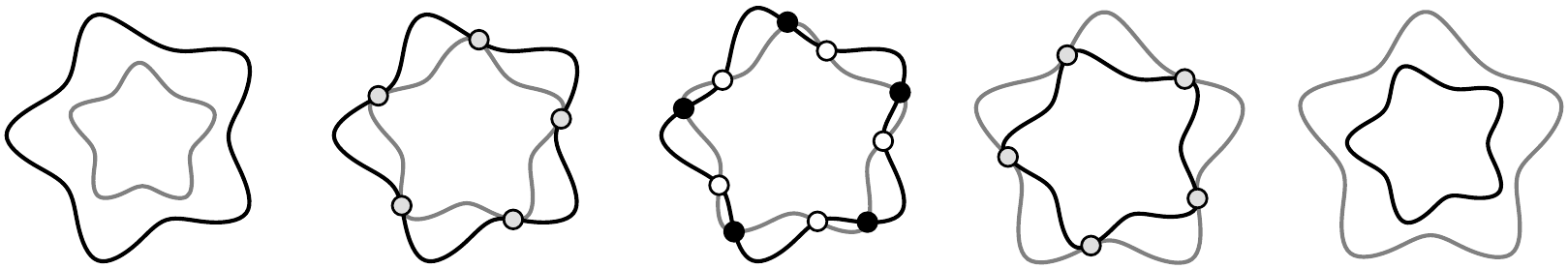} \\
 $(b)$ &  \includegraphics[width=\linewidth]{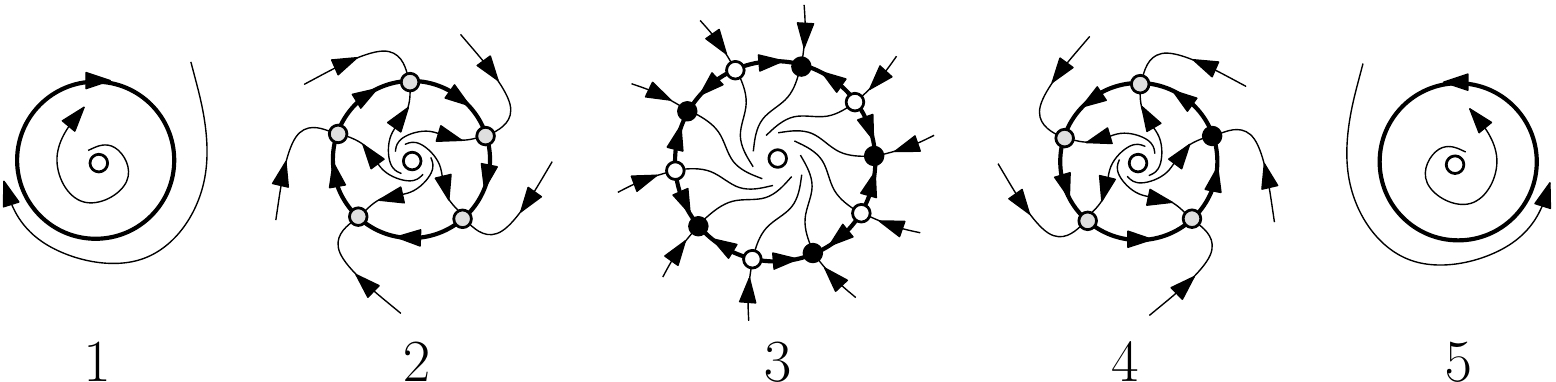} 
\end{tabular}
\end{center}
\caption{Crossing the horn for $m=5$. $(a)$ Fixed point solutions of
  \eqref{HopfZm2} at the five points in figure~\ref{hornsZm}$(b)$;
  grey points in 2 and 4 are the saddle-node points, that split in a
  stable point (black) and saddle (white). $(b)$ Phase portraits
  corresponding to the five cases in $(a)$.}
\label{SNICprocess}
\end{figure}

On L, $r^2=u$; close to L, the fixed points $P_\pm$ are given by
$r^2\sim u\pm\sqrt{\epsilon^2u^{m-2}-v^2}$. The pinning region, at
dominant order in $\epsilon$, is $v=\pm\epsilon u^{(m-2)/2}$. This
gives a wedge-shaped region around L for $m=4$ and a horn for $m>4$,
as illustrated in figure~\ref{hornsZm}. On the boundaries of the
pinning region, the fixed points merge in saddle-node bifurcations
that take place on the limit cycles $C_\pm$. These are curves of SNIC
bifurcations, exactly the same phenomena that is observed in
Neimark-Sacker bifurcations \citep{arpl90}, and that we have
encountered also in the previous cases analyzed in the present
study. Due to the symmetry $Z_m$, from \eqref{HopfZm2} we see that at
the boundaries of the wedge, $m$ simultaneous saddle-node bifurcations
take place. Figure~\ref{SNICprocess} shows how the fixed points
appear and disappear in saddle-node bifurcations on the limit cycle
$C_\pm$ when crossing the horn for the $m=5$ case, at the five points
in parameter space indicated in figure~\ref{hornsZm}$(b)$. The
nontrivial fixed points, from \eqref{HopfZm2} and
$r^2\sim\mu/a\sim\nu/b$, satisfy
\begin{equation}\label{fixpointsZm2}
\begin{aligned}
 & r^2=\mu/a+\epsilon(\mu/a)^{(m-1)/2}\cos m\phi,\\
 & r^2=\nu/b-\epsilon(\nu/b)^{(m-1)/2}\sin m\phi,
\end{aligned}
\end{equation}
and the intersection of these circles modulated by the $\sin m\phi$ and
$\cos m\phi$ terms is illustrated in
figure~\ref{SNICprocess}$(a)$. Phase portraits corresponding to the
five points in parameter space are also schematically shown in
figure~\ref{SNICprocess}$(b)$.

The width of the pinning region in the symmetry breaking $SO(2)\to
Z_m$ case is obtained from the shape of the horn region and is given
by $w(d,\epsilon)=2v=2\epsilon d^{(m-2)/2}$. Again, as in all the
preceding cases, the width of the pinning region is proportional to
the amplitude of the imperfection $\epsilon$.

\section{Common features in the different ways to break $SO(2)$ symmetry}
\label{Sec_general}

\begin{figure}
\begin{center}
\includegraphics[width=0.65\linewidth]{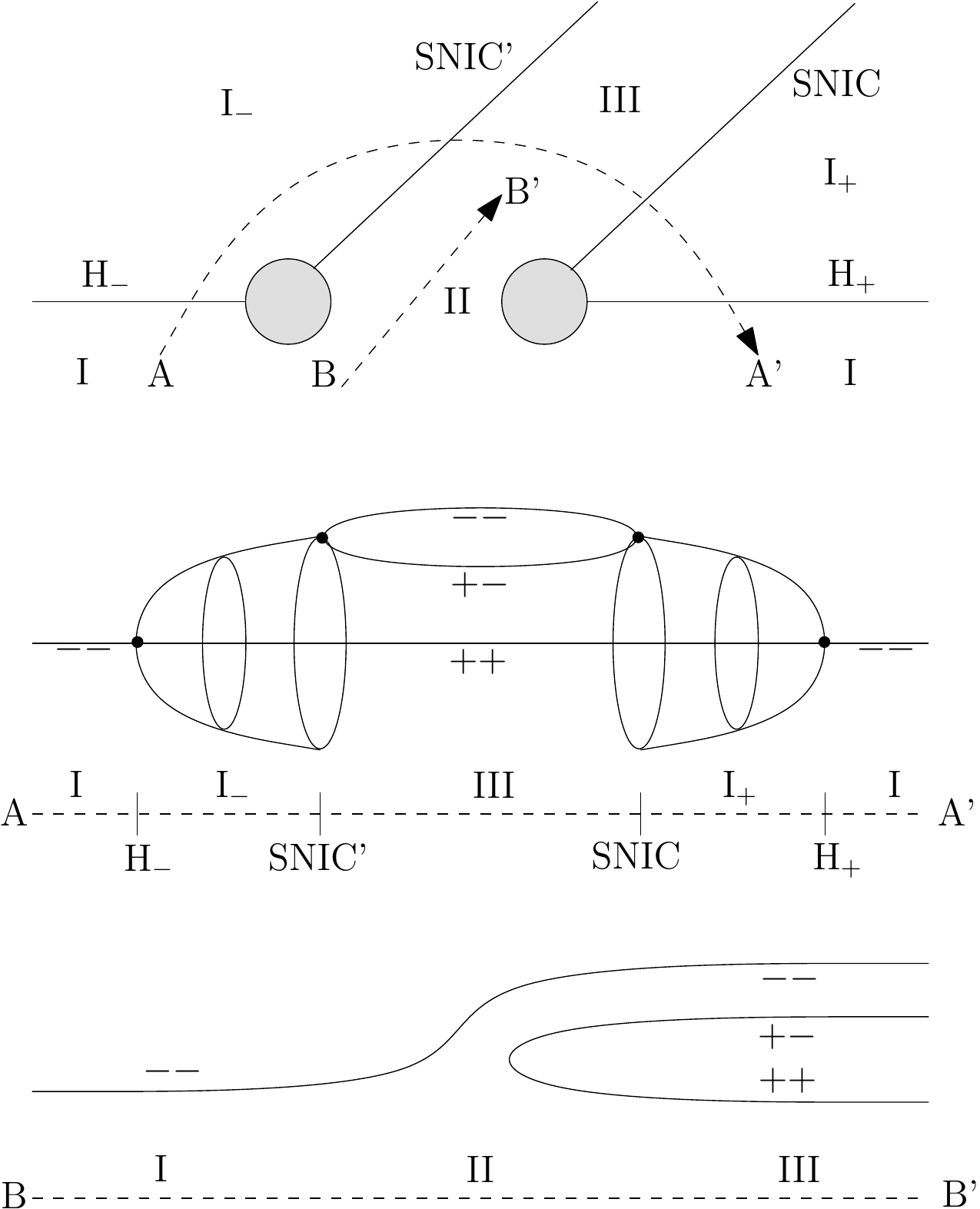}
\end{center}
\caption{Imperfect Hopf under general perturbations: $(a)$ regions in
  parameter space; $(b$) and $(c)$, bifurcation diagrams along the two
  one-dimensional paths, (1) and (2), respectively. The signs ($++$,
  $--$\ldots) indicate the sign of the real part of the two
  eigenvalues of the solution branch considered.}
\label{impHopf_general_min}
\end{figure}

Here we summarize the features that are common to the different
perturbations analyzed in the previous sections. The most important
feature is that the curve of zero frequency splits into two curves
with a region of zero-frequency solutions appearing in between (the
so-called pinning region). Of the infinite number of steady solutions
that exist along the zero-frequency curve in the \emph{perfect} system
with $SO(2)$ symmetry, only a small finite number remain. These steady
solutions correspond to the pinned solutions observed in experiments
and in numerical simulations, like the ones to be described in
\S\ref{Sec_experiments}. The number of remaining steady solutions
depends on the details of the symmetry-breaking imperfections, but
when $SO(2)$ is completely broken and no discrete symmetries remain,
there are three steady solutions in the pinning region III (see
figure~\ref{impHopf_general_min}$a$). One corresponds to the base
state, now unstable with eigenvalues $(+,+)$. The other two are born
on the SNIC curves delimiting region III away from the origin. Of
these two solutions, one is stable (the only observable state in
region III) and the other is a saddle (see
figure~\ref{impHopf_general_min}$b$ and $c$). There are also the two
Hopf bifurcation curves H$_-$ and H$_+$. The regions where the Hopf
bifurcations meet the infinite-period bifurcations cannot be described
in general, and as has been shown in the examples in the previous
sections, will depend on the specifics of how the $SO(2)$ symmetry is
broken, i.e.\ on the specifics of the imperfections present in the
problem considered. These regions contain complex bifurcational
processes, and are represented as grey disks in
figure~\ref{impHopf_general_min}$(a)$. The stable limit cycle
existing outside III, in regions I$_\pm$, undergoes a SNIC bifurcation
and disappears upon entering region III (see
figure~\ref{impHopf_general_min}$b$). When the SNIC bifurcation curves
approach the Hopf bifurcation curves (i.e.\ enter the grey disks
regions), the saddle-node bifurcations do not occur on the stable
limit cycle but very close, and the limit cycle disappears in a
saddle-loop homoclinic collision that occurs very close to the
saddle-node bifurcations. These homoclinic collisions behave like a
SNIC bifurcation, except in a very narrow region in parameter space
around the saddle-node curves, as has been discussed in
\S\ref{glob_bif_subsect}.

In all cases considered, the width of the pinning region scales
linearly with the strength of the symmetry breaking $\epsilon$. In all
cases, we have found $w(d,\epsilon)=2\epsilon d^{(p-1)/2}$, where $p$
is the order of the symmetry breaking considered. For lower order
terms, the width decreases ($\epsilon$ case, order zero) or remains
constant ($\epsilon\bar z$ case, order one) with increasing distance
from the bifurcation point. For quadratic and higher order terms, the
width increases with the distance. When arbitrary perturbations are
included, we expect a behavior of the form $w(d,\epsilon)=\epsilon
f(d)$, where the function $f$ will depend on the details of the
symmetry-breaking terms involved. The size of the regions containing
complex bifurcational processes (the grey disks in
figure~\ref{impHopf_general_min}$a$) is of order $\epsilon$ or
smaller, as we have seen in all cases considered. Therefore these
regions are comparable in size or smaller than the width of the
pinning region.

\begin{figure}
\begin{center}
\begin{tabular}{cc}
 $(a)$ & $(b)$ \\
\includegraphics[width=0.47\linewidth]{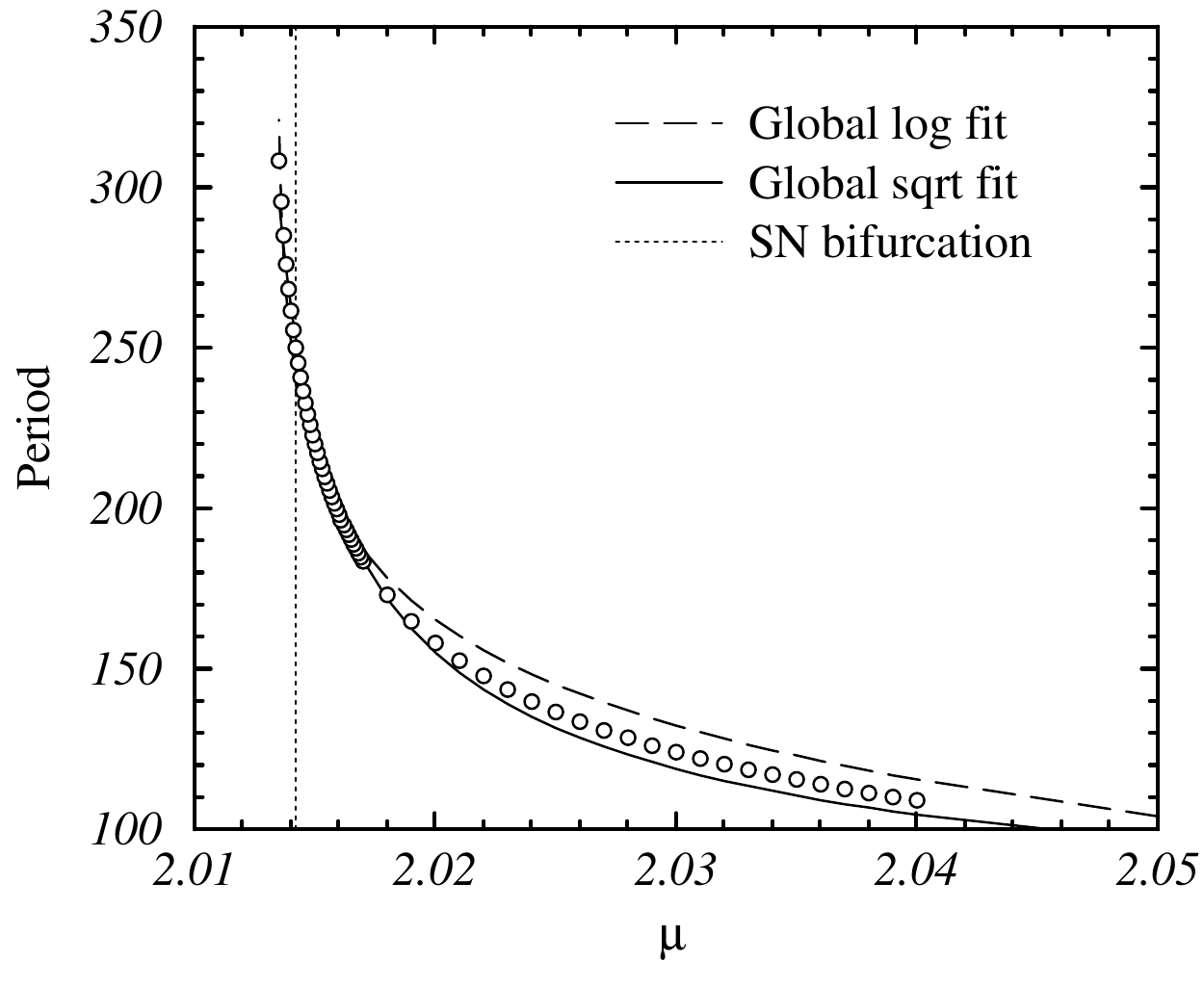} &
\includegraphics[width=0.47\linewidth]{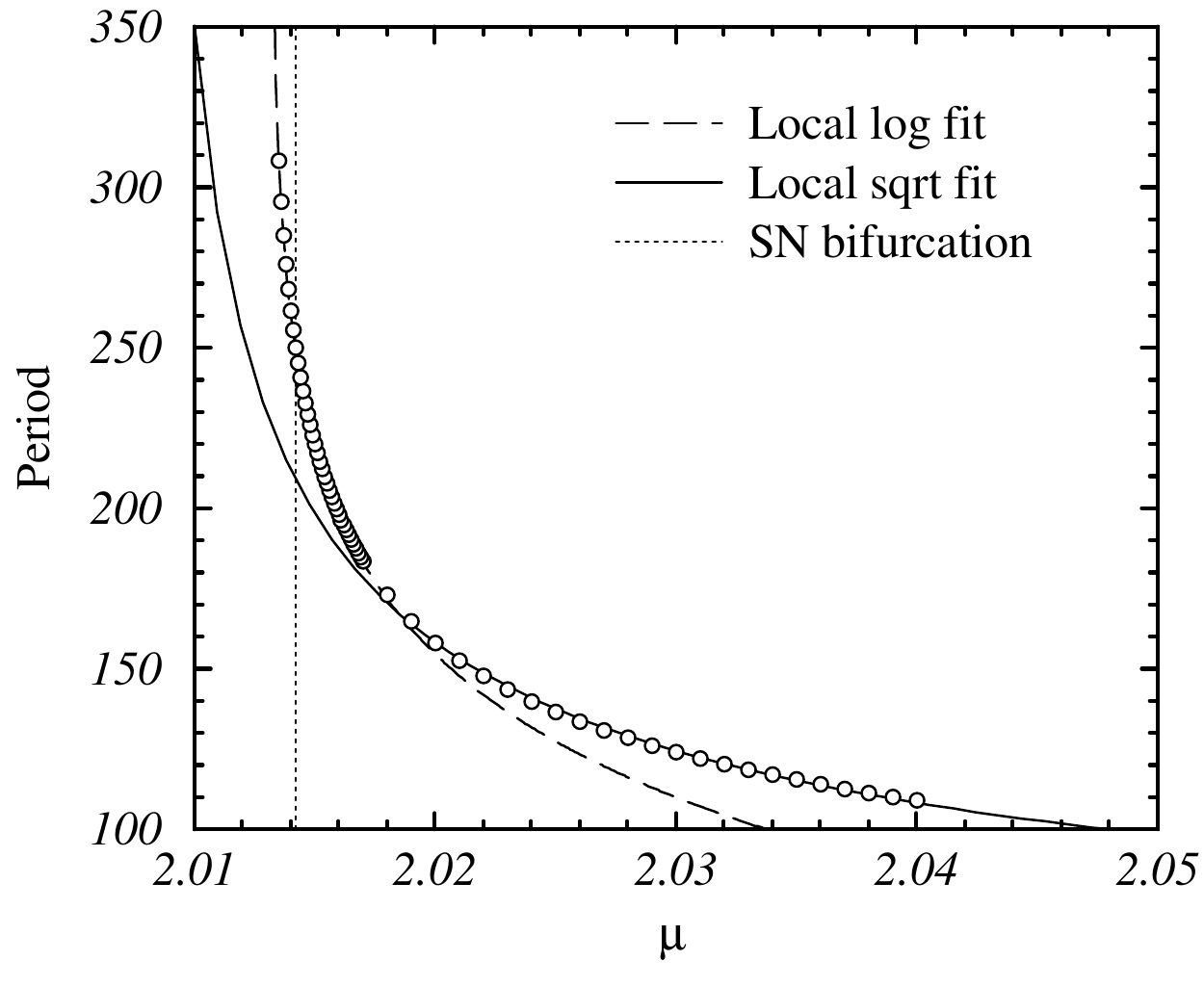}
\end{tabular}
\end{center}
\caption{Square root and logarithmic fits to the periods of $C_-$
  approaching the Hom$_-$ curve at $\nu=0.6$: $(a)$ fitted curves in the
  complete range $\mu\in[2.0135,2.2]$, and $(b)$ fits in the range
  $\mu\in[2.0135,2.027]$ for the log fit and $\mu\in[2.018,2.2]$ for
  the square root fit.}
\label{SNIC_versus_Hom}
\end{figure}

\subsection{SNIC versus homoclinic collision: scaling of the period}
\label{SNICvsHom}

Another feature we have found in the five scenarios discussed above is
that the SNIC$_\pm$ bifurcations, where the stable limit cycle
disappears on entering the pinning region away from the origin
($\mu=\nu=0$), become homoclinic (or heteroclinic if the $SO(2)$
symmetry is not completely broken) collisions between the stable limit
cycle and the saddle point that is born in the saddle-node
bifurcations SN$_\pm$ that now take place not on the limit cycle but
very close to it. There is a codimension-two global bifurcation,
termed SnicHom in the bifurcation diagrams discussed above, where the
SNIC curve, the saddle-node curve and the homoclinic collision curve
meet.

The scaling laws of the periods when approaching a
homoclinic or a SNIC bifurcation are different, having logarithmic
or square root profiles:
\begin{equation}\label{per_fits}
  T_\text{Het}=\frac{1}{\lambda}\ln\frac{1}{\mu-\mu_c}+O(1),\quad
  T_\text{SNIC}=\frac{k}{\sqrt{\mu-\mu_c}}+O(1),
\end{equation}
where $\lambda$ is the positive eigenvalue of the saddle, and $k$ a
constant. So the question is what happens with the
scaling of the period close to the SnicHom bifurcation?

We have analized in detail the $\epsilon\bar z$ case, being the other
cases very similar. We refer to the figures~\ref{impHopf_regions} and
\ref{SNIC_Het} From a practical point of view, close to but before the
SnicHet$_-$ point, the saddle-node bifurcation SN$_-$ is very closely
followed by the heteroclinic collision of the limit cycle $C_-$ with
the saddles $P_-$ and $P^*_-$, and it becomes almost indistinguishable
from the SNIC$_-$ bifurcation.  We have numerically computed the
period of $C_-$ at $\nu=0.6$ for decreasing $\mu$ values approaching
SN$_-$ in the range $\mu\in[2.0135,2.2]$, for $\alpha_0=45^\circ$.
Figure~\ref{SNIC_versus_Hom}$(a)$ shows both fits using the values of
the period over the whole computed range. The log fit overestimates
the period while the square-root fit underestimates it, and this
underestimate gets larger as the heteroclinic collision is
approached. Figure~\ref{SNIC_versus_Hom}$(b)$ again shows both fits,
but now using values close to the collision for the log fit and values
far away from the collision for the square-root fit. Both fits are now
very good approximations of the period in their corresponding
intervals, and together cover all the values numerically
computed. When the interval between the SN$_-$ bifurcation and the
heteroclinic collision (in figure~\ref{SNIC_versus_Hom},
$\mu_\text{SN}=2.01420$ and $\mu_c=2.01336$, respectively) is very
small, it cannot be resolved experimentally (or even numerically in an
extended systems with millions of degrees of freedom, as is the case
in fluid dynamics governed by the three-dimensional Navier--Stokes
equations). In such a situation the square-root fit looks good enough,
because away from the SN$_-$ point, the dynamical system just feels
the ghost of the about to be formed saddle-node pair and does not
distinguish between whether the saddle-node appears on the limit cycle
or very close to it. However, if we are able to resolve the very
narrow parameter range between the saddle-node formation and the
subsequent collision with the saddle, then the log fit matches the
period in this narrow interval much better.

Due to the presence of two very close bifurcations (Het' and SN$_-$), the
scaling laws become cross-contaminated, and from a practical point of
view the only way to distinguish between a SNIC and a Homoclinic
collision is by computing or measuring periods extremelly close to the
infinite-period bifurcation point. We can also see this from the log
fit equation in \eqref{per_fits}; when both bifurcations are very
close, $\lambda$, the positive eigenvalue of the saddle, goes to zero
(it is exactly zero at the saddle-node point), so the log fit becomes
useless, except when the periods are very large.

When the interval between
the SN$_-$ bifurcation and the homoclinic collision is very small, it
cannot be resolved experimentally (or even numerically in an extended
systems with millions of degrees of freedom, as is the case in fluid
dynamics governed by the three-dimensional Navier--Stokes
equations). In such a situation, the square-root fit appears to be
good enough, because away from the SN$_-$ point, the dynamical system
just feels the ghost of the about-to-be-formed saddle-node pair and
does not distinguish between whether the saddle-node appears on the
limit cycle or just very close to it. However, if we are able to
resolve the very narrow parameter range between the saddle-node
formation and the subsequent collision with the saddle, then the log
fit matches the period in this narrow interval much better.

\subsection{Codimension-two bifurcations of limit cycles}\label{Sec_codtwo_LC}

The bifurcations that a limit cycle can undergo have been an active
subject of research since dynamical systems theory was born. Even in
the case of isolated codimension-one bifurcations, a complete
classification was not completed until fifteen years ago
\citep{TuSh95}, when the blue-sky catastrophe was found. The seven
possible bifurcations are: the Hopf bifurcation, where a limit cycle
shrinks to a fixed point, and the length of the limit cycle reduces to
zero. Three bifurcations where both the length and period of the limit
cycle remain finite: the saddle node of cycles (or cyclic fold), the
period doubling and the Neimark-Sacker bifurcations. Two bifurcations
where the length remains finite but the period goes to infinity: the
collision of the limit cycle with an external saddle forming a
homoclinic loop, and the appearance of a saddle-node of fixed points
on the limit cycle (the SNIC bifurcation). Finally, there is the
blue-sky bifurcation where both the length and period go to infinity,
corresponding to the appearance of a saddle-node of limit cycles
transverse to the given limit cycle. The seven bifurcations are
described in many books on dynamical systems,
e.g.\ \citet{SSTC01,Kuz04}; they are also described on the web page
{\tt http://www.scholarpedia.org/article/Blue-sky\_catastrophe}
maintained by A.\ Shil'nikov and D.\ Turaev.

Of the seven bifurcations, only four (Hopf, cyclic fold, homoclinic
collision and SNIC) are possible in planar systems, as is the case in
the present study, and we have found the four of them in the
different scenarios explored. We have also found a number of
codimension-two bifurcations of limit cycles. For these bifurcations a
complete classification is still lacking, and it is interesting to
list them because some of the bifurcations obtained are not very
common.  The codimension-two bifurcations of limit cycles associated
to codimension-two bifurcations of fixed points can be found in many
dynamical systems texbooks, and include Takens--Bogdanov bifurcations
(present in almost all cases considered here) and the Bautin bifurcation
(in the $\epsilon$ case). Codimension-two bifurcations of limit cycles
associated only to global bifurcations are not so common. We have
obtained five of them, that we briefly summarize here.
\begin{description}

\item[PfGl] A gluing bifurcation with the saddle point undergoing a
  pitchfork bifurcation; it may happen in systems with $Z_2$ symmetry;
  see \S\ref{Sec_bar_z}.

\item[CfHom] A cyclic-fold and a homoclinic collision occurring
  simultaneously; see \S\ref{Sec_bar_z}.

\item[CfHet] A cyclic-fold and a heteroclinic collision occurring
  simultaneously; see \S\ref{Sec_epsilon}.

\item[SnicHom] A SNIC bifurcation and a homoclinic collision occurring
  simultaneously; see \S\ref{Sec_quadratic}.

\item[SnicHet] A double SnicHom bifurcation  occurring
  simultaneously in $Z_2$ symmetric systems; see \S\ref{Sec_bar_z}.

\end{description}
The SnicHom bifurcation is particularly important in our problem because it
separates the two possible scenarios upon entering the pinning region:
the stable limit cycle outside may disappear in a homoclinic collision
or a SNIC bifurcation.

\section{Fluid dynamics examples of pinning due to breaking the
  $SO(2)$ symmetry}\label{Sec_experiments}

\subsection{Pinning in small aspect ratio Taylor-Couette flow}
\label{TCpinning}

Experiments in small aspect-ratio Taylor--Couette flows have reported
the presence of a band in parameter space where rotating waves become
steady non-axisymmetric solutions (a pinning effect) via
infinite-period bifurcations \citep{PSL91}. Previous numerical
simulations, assuming $SO(2)$ symmetry of the apparatus, were unable
to reproduce these observations \citep{MaLo06}. Recent additional
experiments suggest that the pinning effect is not intrinsic to the
dynamics of the problem, but rather is an extrinsic response induced
by the presence of imperfections that break the $SO(2)$ symmetry of
the ideal problem. Additional controlled symmetry-breaking
perturbations were introduced into the experiment by tilting one of
the endwalls \citep{AHHP08}. \citet{PLM11} conducted direct numerical
simulations of the Navier--Stokes equations including the tilt of one
endwall by a very small angle. Those simulations agree very well with
the experiments, and the normal form theory developed in this paper
provides a theoretical framework for understanding the observations. A
brief summary of those results follows.

\begin{figure}
  \begin{center}
    \includegraphics[width=0.6\linewidth]{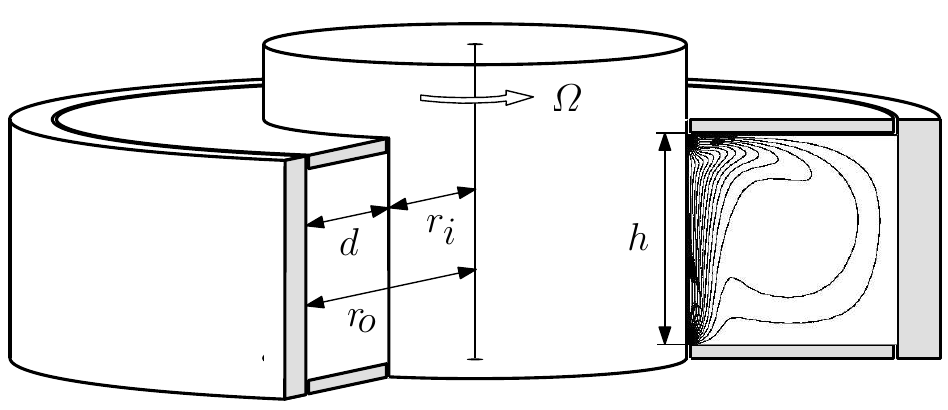}\vspace*{-12pt}
  \end{center}
  \caption{Schematic of the Taylor-Couette apparatus.}
  \label{schemTC}
\end{figure}

Taylor--Couette flow consists of a fluid confined in an annular region
with inner radius $r_i$ and outer radius $r_o$, capped by endwalls a
distance $h$ apart. The endwalls and the outer cylinder are
stationary, and the flow is driven by the rotation of the inner
cylinder at constant angular speed $\varOmega$ (see
figure~\ref{schemTC} for a schematic). The system is governed by three
parameters: 
\begin{equation}\begin{aligned}
 & \text{the Reynolds number} & & Re=\varOmega r_i (r_o-r_i)/\nu, \\
 & \text{the aspect ratio}    & & \Gamma=h/(r_o-r_i), \\
 & \text{the radius ratio}    & & \eta=r_i/r_o, 
\end{aligned}\end{equation}
where $\nu$ is the kinematic viscosity of the fluid. The system is
invariant to arbitrary rotations about the axis, $SO(2)$ symmetry, and
to reflections about the mid-height, a $Z_2$ symmetry that commutes
with $SO(2)$. In both the experiments and the numerical simulations,
the radius ratio was kept fixed at $\eta=0.5$. $Re$ and $\Gamma$ were
varied, and these correspond to the parameters $\mu$ and $\nu$ in the
normal forms studied here.

\begin{figure}
  \begin{center}
    \begin{tabular}{m{0.02\linewidth}@{\hspace{12pt}}m{0.45\linewidth}
        @{\hspace{4mm}}m{0.02\linewidth}@{\hspace{12pt}}m{0.4\linewidth}}
      $(a)$ & \includegraphics[width=\linewidth]{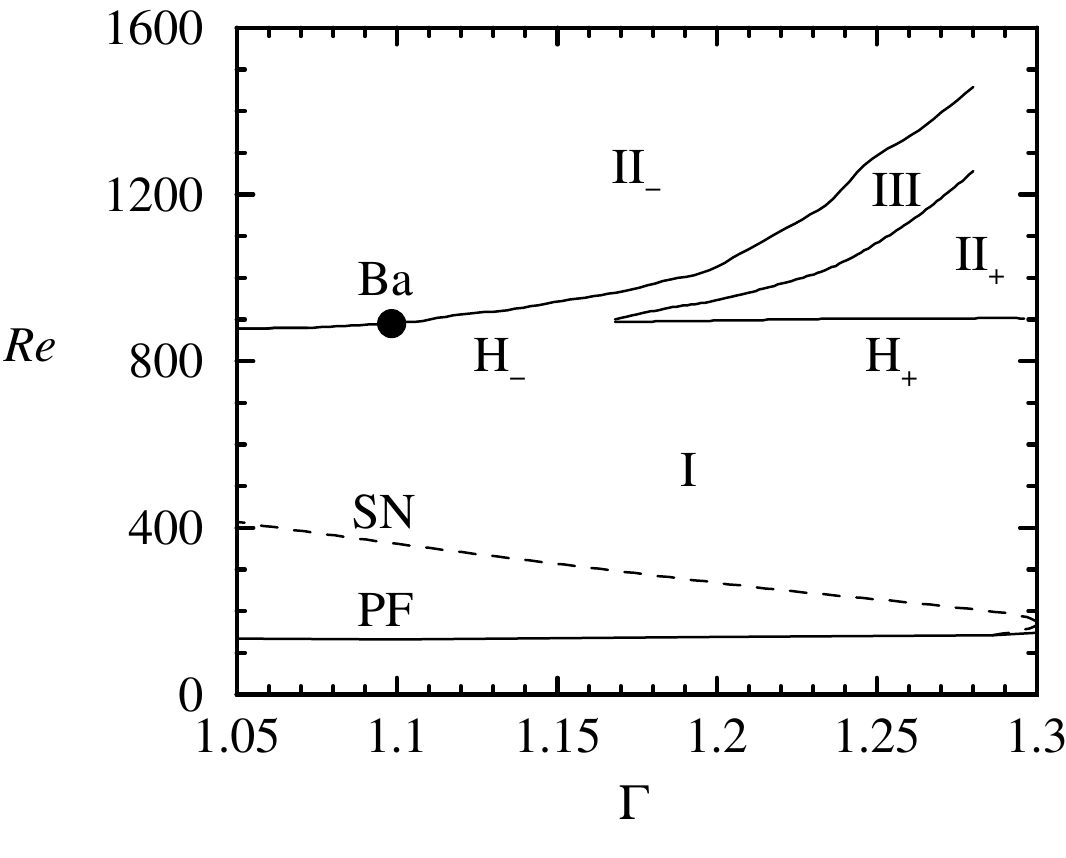} &
      $(b)$ & \includegraphics[width=\linewidth]{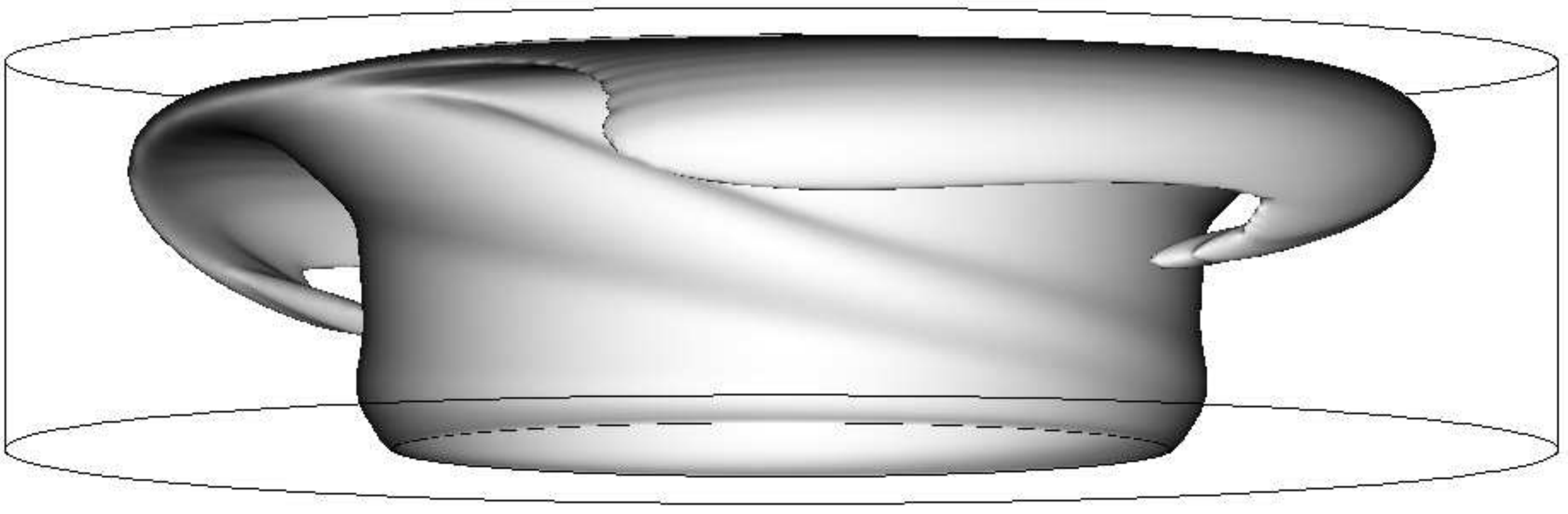}
    \end{tabular}\vspace*{-20pt}
  \end{center}
  \caption{$(a)$ Experimental regimes found in the small aspect-ratio
    Taylor-Couette problem, with a pinning region, adapted from
    \citet{PSL91}. $(b)$ Numerically computed rotating wave
    from \citet{MaLo06}, found in regions II$_\pm$.}
  \label{pinning_TC}
\end{figure}

For small $Re$, below the curve PF in figure~\ref{pinning_TC}$(a)$,
the flow is steady, axisymmetric and reflection symmetric, consisting
of two Taylor vortices \citep{MaLo06}. The $Z_2$ reflection symmetry
is broken in a pitchfork bifurcation along the curve PF, and a pair of
steady axisymmetric one-vortex states that have a jet of angular
momentum emerging from the inner cylinder boundary layer near one or
other of the endwalls is born. Both are stable, and which is realized
depends on initial conditions. The only symmetry of these
symmetrically-related solutions is $SO(2)$. The inset in
figure~\ref{schemTC} shows the azimuthal velocity associated with the
state with the jet near the top. These steady axisymmetric one-vortex
states are stable in region I. There are other flow states that are
stable in this same region. For example, above the dashed curve SN in
figure~\ref{pinning_TC}$(a)$, the two-Taylor-vortex state becomes
stable and coexists with the one-vortex states. However, the
two-vortex and the one-vortex states are well separated in phase space
and the experiments and numerics we describe below are focused on the
one-vortex state. On increasing $Re$, the one-vortex state suffers a
Hopf bifurcation that breaks the $SO(2)$ symmetry and a rotating wave
state emerges with azimuthal wave number
$m=2$. Figure~\ref{pinning_TC}$(b)$ shows an isosurface of axial
angular momentum, illustrating the three-dimensional structure of the
rotating wave. For slight variations in aspect ratio, the rotating
wave may precess either prograde (in region II$_+$ above H$_+$) or
retrograde (in region II$_-$ above H$_-$) with the inner cylinder, and in
between a pinning region III is observed. This is observed even with a
nominally perfect experimental system, i.e.\ with the $SO(2)$ symmetry
to within the tolerances in constructing the apparatus. The Hopf
bifurcation is supercritical around the region where the precession
frequency changes sign. However, for smaller aspect ratios the Hopf
bifurcation becomes subcritical at the Bautin point Ba in
figure~\ref{pinning_TC}$(a)$.

\begin{figure}
  \begin{center}
    \begin{tabular}{m{0.02\linewidth}@{\hspace{6pt}}m{0.43\linewidth}
        @{\hspace{4mm}}m{0.02\linewidth}@{\hspace{6pt}}m{0.43\linewidth}}
      $(a)$ & \includegraphics[width=\linewidth]{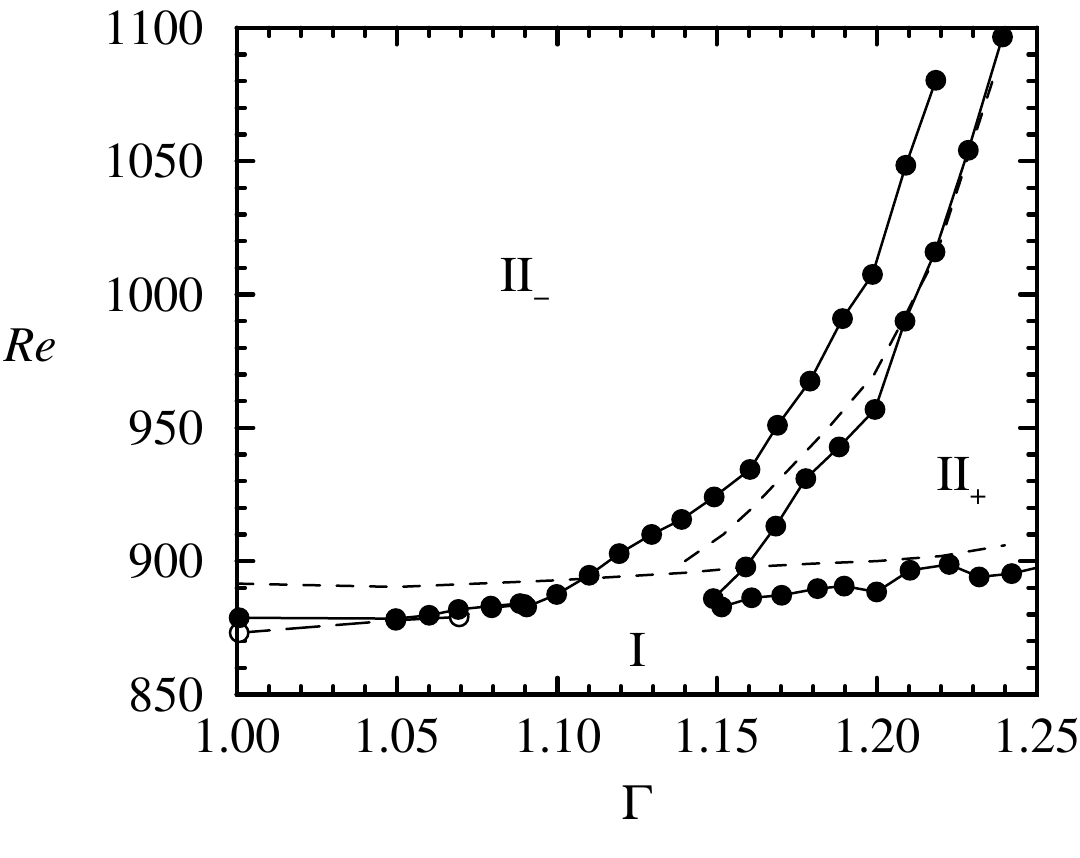} &
      $(b)$ & \includegraphics[width=\linewidth]{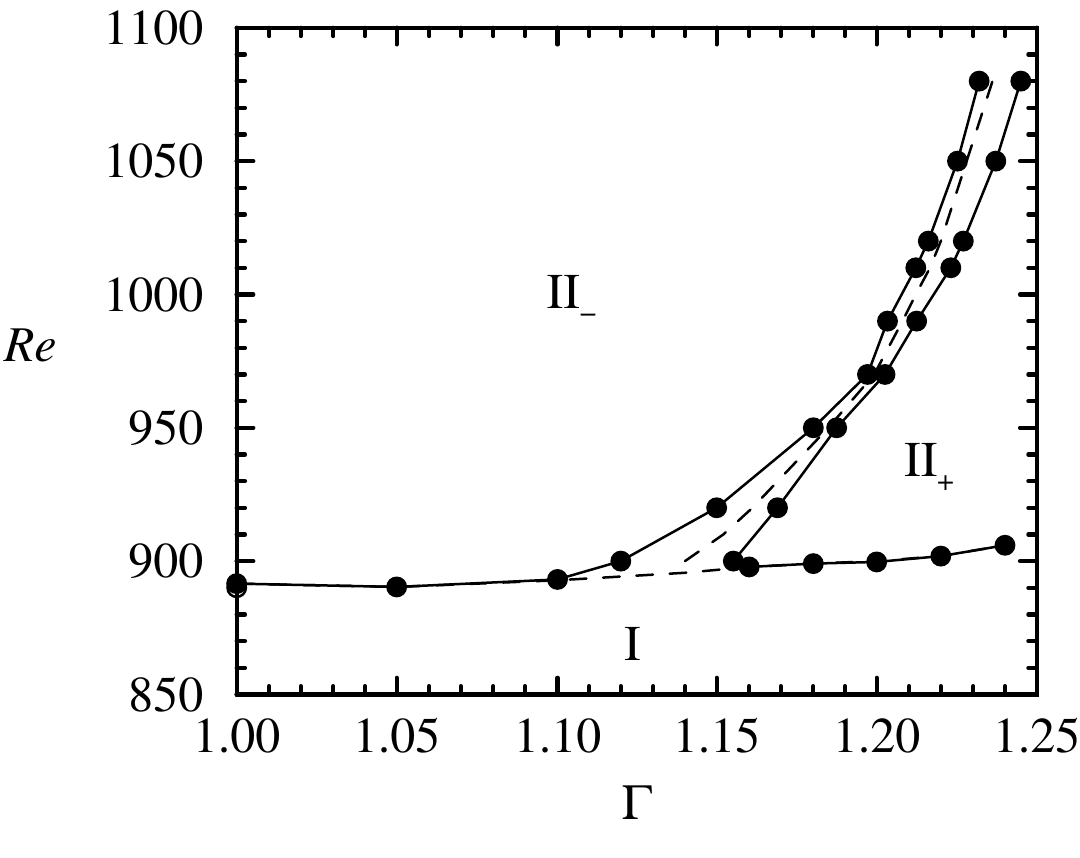}
    \end{tabular}\vspace*{-20pt}
  \end{center}
  \caption{Bifurcation diagrams for the one-cell state from $(a)$ the
    experimental results of \citet{AHHP08} with the natural
    imperfections of their system, and $(b)$ the numerical results of
    \citet{PLM11} with a tilt of $0.1^\circ$ on the upper lid. The
    dotted curve in both is the numerically determined Hopf curve with
    zero tilt.}
  \label{bif_diagram_exp}
\end{figure}

Various different experiments in this regime have been conducted in
the nominally perfect system, as well as with a small tilt of an
endwall \citep{PSCM88,PSL91,PBE92,AHHP08}.
Figure~\ref{bif_diagram_exp}$(a)$ shows a bifurcation diagram from the
laboratory experiments of \citet{AHHP08}. These experiments show that
without an imposed tilt, the natural imperfections of the system
produce a measurable pinning region, and that the additional tilting
of one endwall increases the extent of the pinning region. Tilt angles
of the order of $0.1^\circ$ are necessary for the tilt to dominate
over the natural imperfections. Figure~\ref{bif_diagram_exp}$(b)$
shows a bifurcation diagram from the numerical results of
\citet{PLM11}, including a tilt of one of the endwalls of about
$0.1^\circ$, showing very good agreement with the experimental
results. Included in figures~\ref{bif_diagram_exp}$(a)$ and $(b)$ is
the numerically computed bifurcation diagram in the perfect system,
shown as dotted curves. The effects of imperfections are seen to be
only important in the parameter range where the Hopf frequency is
close to zero and a pinning region appears. It is bounded by
infinite-period bifurcations of limit cycles. The correspondence
between these results and the normal form theory described in this
paper is excellent, strongly suggesting that the general remarks on
pinning extracted from the analysis of the five particular cases are
indeed realized both experimentally and numerically. These two studies
\citep{AHHP08,PLM11} are the only cases we know of where quantitative
data about the pinning region are available. Yet, even in these cases
the dynamics close to the intersection of the Hopf curve with the
pinning region, that according to our analysis should include
complicated bifurcational processes, has not been explored either
numerically or experimentally. This is a very interesting problem that
deserves further exploration.

\subsection{Pinning in rotating Rayleigh-B\'enard convection}
\label{conrot_pinning}

\begin{figure}
  \begin{center}
    \begin{tabular}{m{0.02\linewidth}@{\hspace{5mm}}m{0.36\linewidth}
        @{\hspace{8mm}}m{0.02\linewidth}@{\hspace{5mm}}m{0.36\linewidth}}
      $(a)$ & \includegraphics[width=\linewidth]{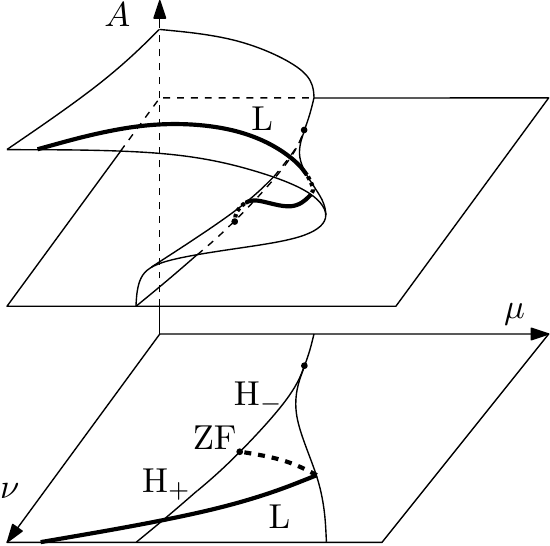}&
      $(b)$ & \includegraphics[width=\linewidth]{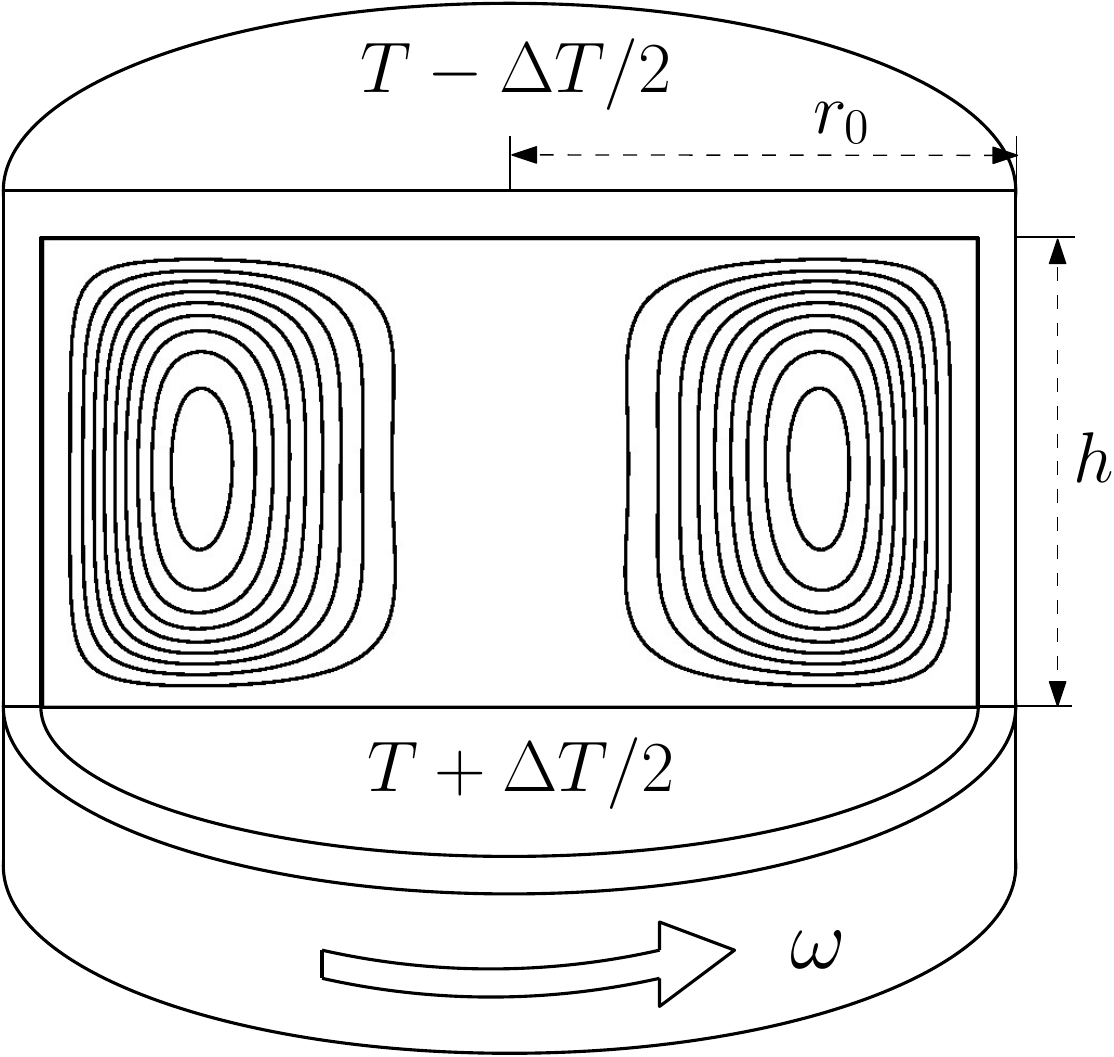}\\
    \end{tabular}\vspace*{-12pt}
  \end{center}
  \caption{$(a)$ Schematic of the Bautin bifurcation, including the
    path L of frequency zero bifurcated states in the $SO(2)$ perfect
    system, shown in a two parameter space $(\mu, \nu)$ with $A$ a
    global measure of the solution. Also shown is the projection of
    the saddle-node surface on parameter space. $(b)$ Schematic of the
    rotating convection apparatus, with the streamlines of the basic
    state shown in the inset.}
  \label{schemRC}
\end{figure}

Up to now, we have considered the zero-frequency Hopf problem in the
context of a supercritical Hopf bifurcation. However, in the
Taylor--Couette example discussed in the previous section, the zero
frequency occurs quite close to a Bautin bifurcation, at which the
Hopf bifurcation switches from being supercritical to subcritical, and
a natural question is what are the consequences of the zero-frequency
occurring on a subcritical Hopf bifurcation. The normal form theory
for the behavior local to the Hopf bifurcation carries over by
changing the direction of time and the sign of the parameters $\mu$
and $\nu$ as discussed before, but then both the limit cycle and the
pinned state are unstable and not observable in a physical experiment
or direct numerical simulation. The limit cycle becomes observable as
it undergoes a saddle-node of limit cycles (a cyclic fold) bifurcation
at the fold associated with the Bautin bifurcation (see a schematic in
figure~\ref{schemRC}$a$), and we expect that the pinned state does
likewise with a saddle-node of fixed points bifurcation along the same
fold. In this subsection, we identify a rotating convection problem
where precisely this occurs \citep{MMBL07,LoMa09}, and conduct new
numerical simulations by introducing an $SO(2)$ symmetry-breaking
bifurcation that produces a pinning region on the upper branch of the
subcritical Hopf bifurcation.

The rotating convection problem consists of the flow in a circular
cylinder of radius $r_0$ and height $h$, rotating at a constant rate
$\omega$~rad/s. The cold top and hot bottom endwalls are maintained at
constant temperatures $T_0\mp0.5\Delta T$, where $T_0$ is the mean
temperature and $\Delta T$ is the temperature difference between the
endwalls. The sidewall has zero heat flux. Figure~\ref{schemRC}$(b)$
shows a schematic of the flow configuration.

Using the Boussinesq approximation that all fluid properties are
constant except for the density in the gravitational and centrifugal
buoyancy terms, and using $h$ as the length scale, $h^2/\kappa$ as the
time scale, and $\Delta T$ as the temperature scale, the governing
equations written in the rotating frame of reference are:
\begin{align}
 & (\partial_t + \bm{u}\cdot\nabla)\bm{u}= -\nabla p +
   \sigma\nabla^2\bm{u} +\sigma Ra\Theta\hat{z}
   +2\sigma\varOmega\bm{u}\times\hat{z} -\frac{\sigma Fr Ra}{\gamma}
   (\Theta-z)\bm{r},\label{goveq1}\\
 & (\partial_t +\bm{u}\cdot\nabla)\Theta=w+\nabla^2\Theta,\qquad\quad
   \nabla\cdot\bm{u}=0\,,
\end{align}
where $\bm{u}=(u,v,w)$ is the velocity field in cylindrical
coordinates $(r,\theta,z)$ in the rotating frame, $p$ is the kinematic
pressure (including gravitational and centrifugal contributions),
$\hat{z}$ the unit vector in the vertical direction $z$, and
$\bm{r}$ is the radial vector in cylindrical coordinates. Instead
of the non-dimensional temperature $T$, we have used the temperature
deviation $\Theta$ with respect to the conductive profile,
$T=T_0/\Delta T -z+\Theta$, as is customary in many thermal convection
studies.

There are five non-dimensional independent parameters:
\begin{equation}\begin{aligned}
 & \text{the Rayleigh number} & & Ra=\alpha g h^3 \Delta T/(\kappa\nu),\\
 & \text{the Froude number}   & & Fr=\omega^2 r_0/g,\\
 & \text{the Coriolis number} & & \varOmega=\omega h^2/\nu,\\
 & \text{the Prandtl number}  & & \sigma=\nu/\kappa,\\
 & \text{the aspect ratio}    & & \gamma=r_0/h,
\end{aligned}\end{equation}
where $\alpha$ is the coefficient of volume expansion, $g$ is the
gravitational acceleration, $\kappa$ is the thermal diffusivity, and
$\nu$ is the kinematic viscosity. 

The boundary conditions for $\bm{u}$ and $\Theta$ are:
\begin{alignat}{2}
 & r=\gamma: && \Theta_r=u=v=w=0, \\
 & z=\pm1/2: &\quad& \Theta=u=v=w=0.
\end{alignat}
For any $Fr\ne0$, the system is not invariant to the so-called
Boussinesq symmetry corresponding to invariance to a reflection $K_z$
about the half-height $z=0$, whose action is
$K_z(u,v,w,\Theta,p)(r,\theta,z)=(u,v,-w,-\Theta,p)(r,\theta,-z)$.
The system is only invariant under rotations about the axis of the
cylinder, the $SO(2)$ symmetry.

The governing equations have been solved using a second-order
time-splitting method combined with a pseudo-spectral method for the
spatial discretization, utilizing a Galerkin--Fourier expansion in the
azimuthal coordinate $\theta$ and Chebyshev collocation in $r$ and
$z$. The details are presented in \citet{MBA10}. We have used
$n_r=36$, $n_{\theta}=40$ and $n_z=64$ spectral modes in $r$, $\theta$ and $z$
and a time-step $dt=2\times 10^{-5}$ thermal time units in
all computations. We have checked the spectral convergence of the code
using the infinity norm of the spectral coefficients of the computed
solutions. The trailing coefficients of the spectral expansions are at
least five orders of magnitude smaller than the leading
coefficients. In order to compute the zero-frequency line L in the
subcritical region of the Bautin bifurcation, where the fixed points
and limit cycles involved are unstable, we have used arclength
continuation methods for fixed points and for rotating waves adapted
to our spectral codes \citep{SML02,MBA06}.

\begin{figure}
  \begin{center}
    \begin{tabular}{m{0.02\linewidth}@{\hspace{5mm}}m{0.38\linewidth}
        @{\hspace{8mm}}m{0.02\linewidth}@{\hspace{5mm}}m{0.38\linewidth}}
      $(a)$ & \includegraphics[width=\linewidth]{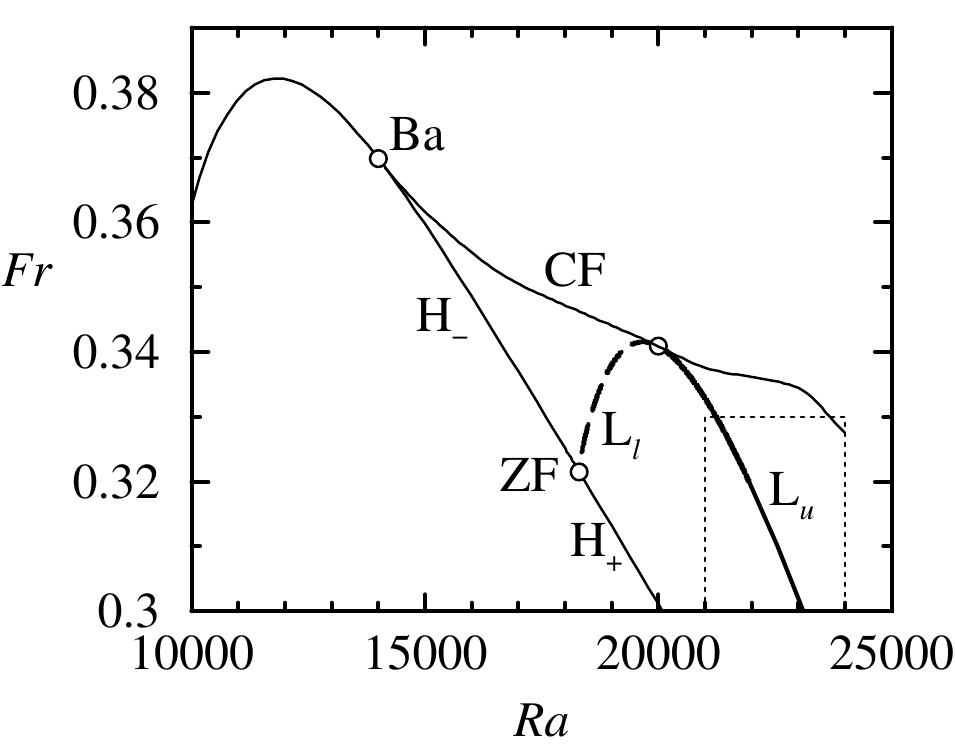}&
      $(b)$ & \includegraphics[width=\linewidth]{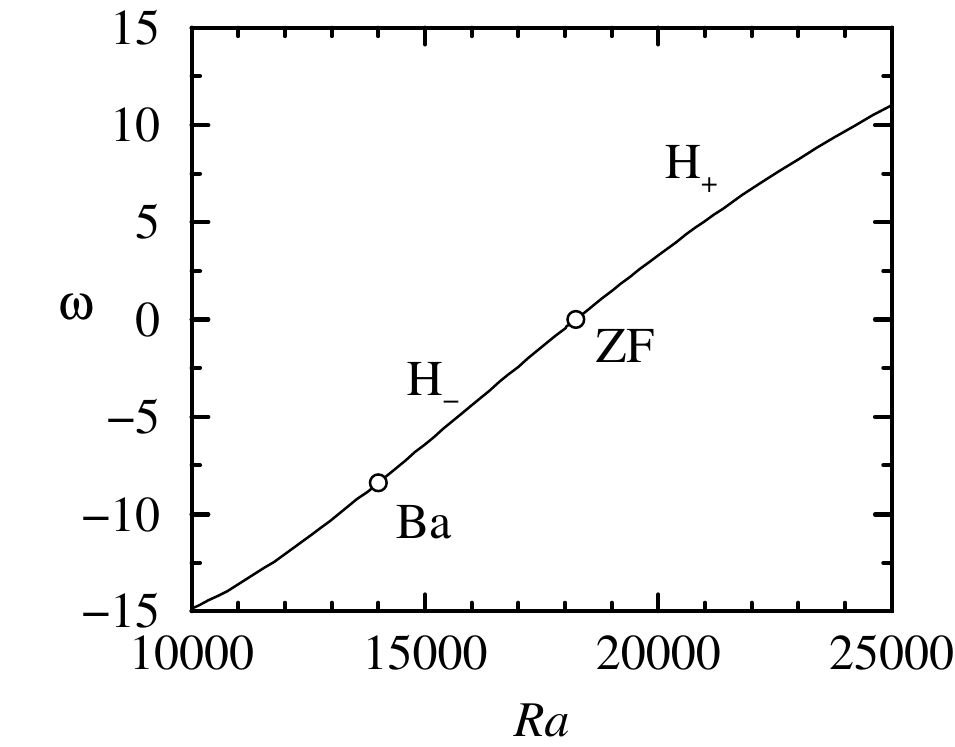}
    \end{tabular}\vspace*{-20pt}
  \end{center}
  \caption{$(a)$ Bifurcation curves for $\varOmega=100$, $\sigma=7$
    and $\gamma=1$, where H$_{\mp}$ are the segments of the Hopf
    bifurcation with negative and positive frequency, the Bautin point
    Ba is where the Hopf bifurcation switches from super- to
    subcritical and the cyclic-fold bifurcation curve CF emerges. The
    lines L$_l$ and L$_u$ are the loci where the rotating wave has
    zero frequency on the lower and upper branches of the cyclic
    fold. The rectangle around L$_u$ corresponds to
    figure~\ref{RCpinning}$a$. $(b)$ The frequency along the Hopf
    bifurcation H$_{\mp}$ with the Bautin point and the point ZF, where
    the sense of precession changes, marked as open symbols.}
  \label{RCbifs}
\end{figure}

Figure~\ref{RCbifs}$(a)$ shows the parameter region of interest in
this convection problem. In the region of high Froude number (region
I) we have a stable steady solution, consisting of a single
axisymmetric convective roll where the warm fluid moves upwards close
to the axis (due to the rotation of the container), and returns
along the sidewall, as illustrated in figure~\ref{schemRC}$(b)$; the inset
shows streamlines of the flow in this base state, that is
$SO(2)$-equivariant with respect to rotations about the cylinder axis.
This base state looses stability when the Froude number $Fr$
decreases, in a Hopf bifurcation along the curves H$_\pm$. The
bifurcation is supercritical for $Ra<14\,157$ and subcritical for higher
$Ra$; the change from supercritical to subcritical happens at the
codimension-two Bautin bifurcation point Ba, at
$(Ra,Fr)\approx(14\,157,0.3684)$. The bifurcated limit cycle, a rotating
wave with azimuthal wave number $m=3$, is unstable, but becomes stable
at the cyclic fold curve CF (a saddle-node bifurcation of limit
cycles). This curve CF originates at the Bautin point Ba. There are
other flow states that are stable in this same region \citep{LoMa09}; these
additional states are well separated in phase space and the numerics
we describe below are focused on the base state and the $m=3$
bifurcated rotating wave.

Figure~\ref{RCbifs}$(b)$ shows the computed frequency of the limit
cycle along the Hopf bifurcation curve. This frequency is negative
along H$_-$ and positive along H$_+$, and is zero at the ZF (zero
frequency) point. At this point we have precisely the scenario
discussed in the present paper: a flow (the base state) with $SO(2)$
symmetry undergoing a Hopf bifurcation that has zero frequency at that
point. Figure~\ref{RCbifs}$(a)$ also includes the line L in parameter
space where the frequency of the bifurcated states is zero. This curve
has been computed using continuation methods since the
zero-frequency state is unstable in the lower part (L$_l$) of the
saddle-node CF, and therefore cannot be obtained via time evolution.
The zero frequency state becomes stable upon crossing the saddle-node
curve CF and moving to the upper part L$_u$ of the saddle-node CF, and
becomes observable both experimentally and by numerical simulations
advancing the Navier-Stokes equations in time.

\begin{figure} 
  \begin{center}
    \begin{tabular}{m{0.02\linewidth}@{\hspace{5mm}}m{0.4\linewidth}
        @{\hspace{6mm}}m{0.02\linewidth}@{\hspace{5mm}}m{0.4\linewidth}}
      $(a)$ & \includegraphics[width=\linewidth]{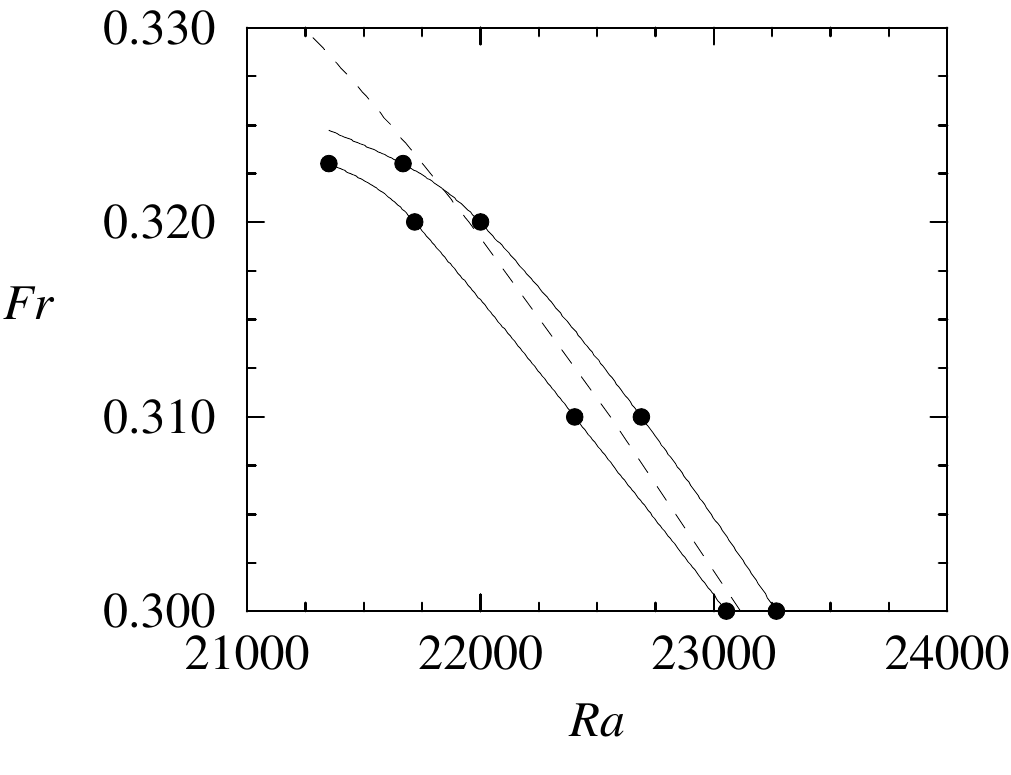} &
      $(b)$ & \includegraphics[width=\linewidth]{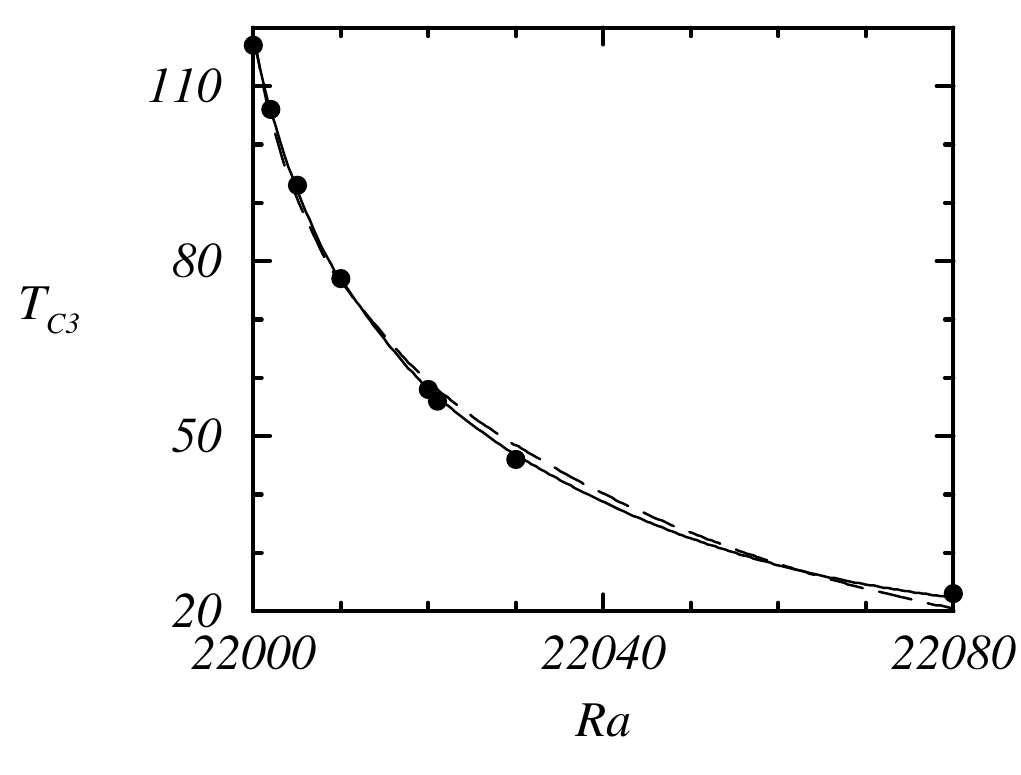} \\
    \end{tabular}\vspace*{-20pt}
  \end{center}
  \caption{$(a)$ SNIC bifurcation curves in $(Ra,Fr)$ space bounding
    the pinning region for $\varOmega=100$, $\sigma=7$ and
    $\gamma=1$. The region shown is the rectangle in
    figure~\ref{RCbifs}$a$, and the dashed line is the L$_u$ curve in
    the perfect case. $(b)$ The period of the rotating wave C3 as it
    approaches the SNIC bifurcations at $Fr=0.32$ and other parameters
    as in $(a)$.}
  \label{RCpinning}
\end{figure}

In order to break the $SO(2)$ symmetry and see if a pinning region
appears, an imperfection has been introduced, in the form of an
imposed linear profile of temperature at the top lid: 
\begin{equation}
  \Theta(r,\theta,z)=\epsilon\,r\cos\theta \quad\text{at}\quad z=1/2,
\end{equation}
where $\epsilon$ is a measure of the symmetry breaking. This term
completely breaks the rotational symmetry of the governing equations,
and no symmetry remains. Figure~\ref{RCpinning}$(a)$ shows that the
line L becomes a band of pinned solutions, steady solutions with
frequency zero, as predicted by the normal form theory presented in
this paper. We can also check the nature of the bifurcation taking
place at the boundary of the pinning
region. Figure~\ref{RCpinning}$(b)$ shows the variation of the period
of the limit cycle approaching the pinning region. It is an infinite
period bifurcation, and the square root fit (shown in the figure)
works better than the logarithmic fit. We estimate that the
bifurcation is a SNIC bifurcation, as the normal form theory presented
predicts it should be sufficiently far from the zero frequency point
ZF.

\begin{figure} 
  \begin{center}
    \begin{tabular}{m{0.02\linewidth}@{\hspace{5mm}}m{0.4\linewidth}
        @{\hspace{6mm}}m{0.02\linewidth}@{\hspace{5mm}}m{0.4\linewidth}}
    $(a)$ & \includegraphics[width=\linewidth]{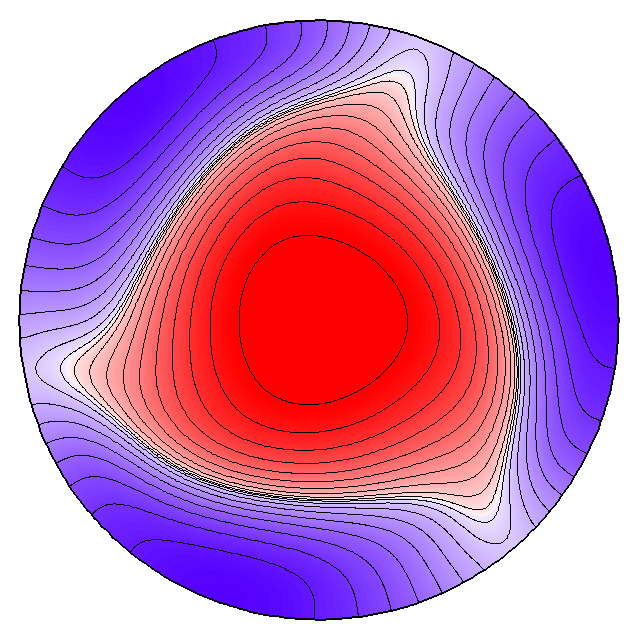} &
    $(b)$ & \includegraphics[width=\linewidth]{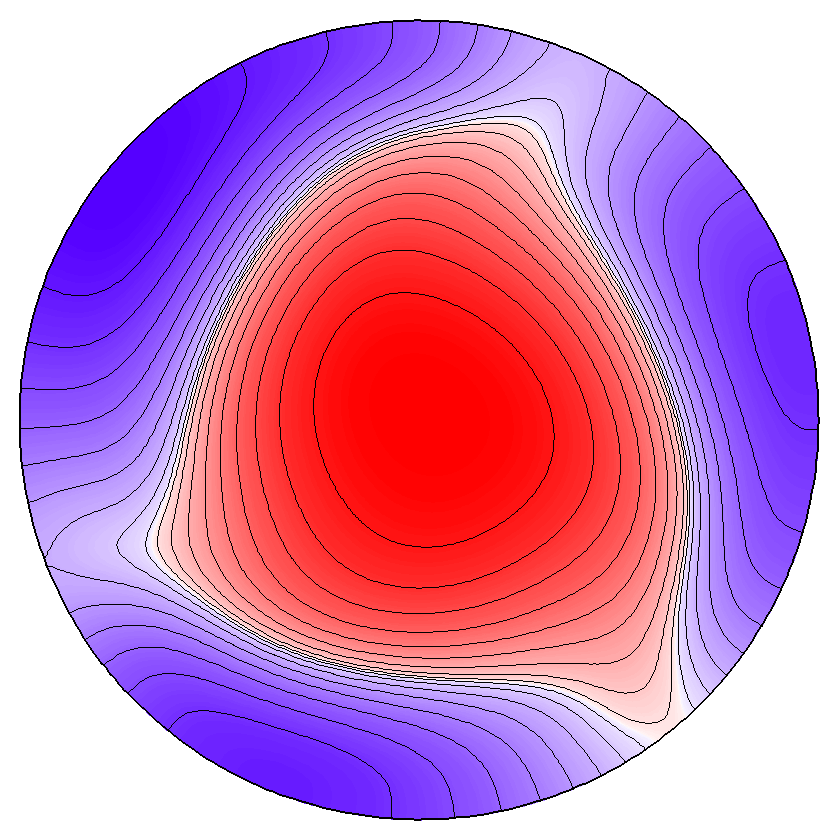}\\
    \end{tabular}\vspace*{-12pt}
  \end{center}
  \caption{Temperature contours at mid height ($z=0$) for
    $\varOmega=100$, $\sigma=7$, $\gamma=1$, $Ra=21950$ and $Fr=0.32$.
    $(a)$ is the symmetric solution without imperfection
    ($\epsilon=0$), and $(b)$ is a pinned solution with an
    imperfection $\epsilon=0.05$. There are 20 quadratically spaced
    contours in the interval $T\in[-0.31,0.31]$, with blue (red) for the
    cold (warm) fluid.}
  \label{temp_contours}
\end{figure}

Figure~\ref{temp_contours} shows contours of the temperature at a
horitzontal section at mid height ($z=0$) for the symmetric system
($\epsilon=0$) in figure~\ref{temp_contours}$(a)$, and for the system
with an imperfection $\epsilon=0.05$, corresponding to a maximum
variation of temperature of 5\%$\Delta T$ at the top lid, in
figure~\ref{temp_contours}$(b)$. The parameter values for both
snapshots are $Ra=21950$ and $Fr=0.32$, inside the pinning region in
figure~\ref{RCpinning}$(a)$. The pinned solution
\ref{temp_contours}$(b)$ has broken the $SO(2)$ symmetries, is a
steady solution, and we can see that one of the three arms of the
solution is closer to the wall than the other two. The attachment of
the solution to the side wall, due to the imperfection at the top lid,
results in the pinning phenomenon.

\section{Summary and conclusions}\label{conclusions}

The aim of this paper has been to provide a general dynamical systems
description of the pinning phenomenon which is observed in systems
possessing two ingredients: slowly traveling or rotating waves and
imperfections. The description boils down to the unfolding of a Hopf
bifurcation in an $SO(2)$ equivariant system about the point where the
Hopf frequency is zero. This turns out to be a very complicated
problem due to the degeneracies involved, but by considering all of
the low-order ways in which $SO(2)$ symmetry may be broken near a zero
frequency Hopf bifurcation, we can identify  a number of general
features which are common to all scenarios, and hence can be expected
to be found in practice. These are that the curve of zero frequency
splits into a region in parameter space of finite width that scales
with the strength of the imperfection, and this region is delimited by
SNIC bifurcations. In the very small neighborhood of the zero
frequency Hopf bifurcation point, where the SNIC curves and the Hopf
curve approach each other, the dynamics is extremely complicated,
consisting in a multitude of codimension-two local bifurcations and
global bifurcations. The details depend on the particulars of the
imperfection, but all of these complications are very localized and
are not resolvable in any practical sense. We provide two examples in
canonical fluid dynamics to illustrate both the pinning phenomenon and
the use of the theory to describe it. These are a Taylor-Couette flow
in which the Hopf bifurcation is supercritical and a rotating
Rayleigh-B\'enard flow where the Hopf bifurcation is subcritical. 

\subsection*{Acknowledgments}
  This work was supported by the National Science Foundation grants
  DMS-05052705 and CBET-0608850, the Spanish Government grants
  FIS2009-08821 and BES-2010-041542, and the Korean Science
  and Engineering Foundation WCU grant R32-2009-000-20021-0.

\appendix
\section{Notation and description of bifurcations}\label{appendix}

\begin{center}\renewcommand{\arraystretch}{1.15}
  \begin{tabular}{l@{\quad}l@{~}}\hline
    \multicolumn{2}{c}{Codimension-one bifurcations}\\\hline
    Name & Description \\\hline
    SN$_{\pm,0}$ & Saddle-node (also called fold) bifurcations \\
    H$_{\pm,0}$ & Hopf bifurcations \\
    PF$_\pm$ & Pitchfork bifurcations \\
    CF & Cyclic fold: two limit cycles are born simultaneously \\
    L, L$_{l,u}$ & Limit cycle becoming a family of fixed points \\
    Hom$_{\pm,0}$, Hom & Homoclinic collision of a limit cycle with a saddle \\
    Het$_{\pm,0}$ & Heteroclinic collision of a limit cycle with saddles \\
    SNIC$_{\pm,0}$ & Saddle-node appearing on a limit cycle \\
    Glu & Gluing bifurcation -- two limit cycles collide with a saddle \\\hline
  \end{tabular}
\end{center}

\begin{center}\renewcommand{\arraystretch}{1.15}
  \begin{tabular}{l@{\quad}l@{~}}\hline
    \multicolumn{2}{c}{Codimension-two bifurcations}\\\hline
    Name & Description \\\hline
    Cusp$_\pm$ & Cusp bifurcations \\
    TB$_\pm$, TB & Takens-Bogdanov bifurcations \\
    dPF$_\pm$ & Degenerate pitchfork -- zero cubic term \\
    Ba & Bautin bifurcation -- degenerate Hopf with zero cubic term \\
    PfGl & Simultaneous gluing Gl and pitchfork PF bifurcations \\
    CfHom & Simultaneous cyclic-fold CF and homoclinic collision Hom \\
    CfHet$_\pm$ & Simultaneous cyclic-fold CF and heteroclinic collision Hom \\
    SnicHom$_{\pm,0}$ & Simultaneous SNIC and homoclinic collision \\
    SnicHet$_{\pm,0}$ & Simultaneous SNIC and heteroclinic collision \\\hline
  \end{tabular}
\end{center}

\section{Symmetry breaking of $SO(2)$ with an $\epsilon$ term: computations}
\label{Appendix_epsilon}

\subsection*{Fixed points}

The fixed points of the normal form \eqref{complexNFepsilon} are
given by $\dot r=\dot\phi=0$, i.e.
\begin{equation}
 \left.\begin{array}{l}
    \cos\phi=r(ar^2-\mu), \\ \sin\phi=r(\nu-br^2), \end{array}\right\}\quad
 \Rightarrow\quad r^2[(\mu-ar^2)^2+(\nu-br^2)^2]=1,
\end{equation}
resulting in the cubic equation $f(\rho)=\rho^3-2(a\mu+b\nu)\rho^2+
(\mu^2+\nu^2)\rho-1=0$, where $\rho=r^2$. This equation always has a
real solution with $\rho>0$, and in some regions in parameter space
may have three solutions. The curve separating these behaviors is a
curve of saddle-node bifurcations, where a couple of additional fixed
points are born. This saddle-node curve is given by
$f(\rho)=f'(\rho)=0$; from these equations we can obtain $(\mu,\nu)$
as a function of $\rho$. In order to describe the curve it is better
to use the rotated reference frame $(u,v)$, where the $u$-axis
coincides with the line L, introduced in \eqref{uv_munu} (see also
figure~\ref{uv_and_epsilon_SN}$(a)$). The saddle-node curve is given by
\begin{equation}
  (u,v)=\frac{1}{2\rho^2_2}\Big(1+2\rho^3_2,\pm\sqrt{4\rho^3_2-1}\,\Big),
  \quad\rho_2\in(2^{-2/3},+\infty),
\end{equation}
where $\rho_2$ is the double root of the cubic equation
$f(\rho)=0$. The third root $\rho_0$ is given by
$\rho_0\rho_2^2=1$. If there is any point where the three roots
coincide (i.e.\ a cusp bifurcation point, where two saddle-node curves
meet), it must satisfy $f(\rho)=f'(\rho)=f''(\rho)=0$. There are two
such points Cusp$_\pm$, given by $\rho_2=1$ and $(u,v)_{\text{Cusp}_\pm}=(3/2,
\pm\sqrt{3}/2)$, dividing the saddle-node curve into three branches:
SN$_0$, joining Cusp$_+$ and Cusp$_-$, and unbounded branches SN$_+$ and SN$_-$
starting at Cusp$_+$ and Cusp$_-$ respectively, and becoming asymptotic to
the line L.  Along SN$_0$, $\rho_2<1<\rho_0$, while along SN$_+$ and SN$_-$,
$\rho_0<1<\rho_2$. At the cusp points, the three roots coincide and
their common value is $+1$.

A better parametrization of the saddle-node curve is obtained by
introducing $s=\pm\sqrt{(4\rho^3_2-1)/3}$, so that now
$s\in(-\infty,+\infty)$, Cusp$_\pm$ corresponds to $s=\pm1$ and 
\begin{equation}\label{epsSNcurve}
 (u,v)=\frac{(3(1+s^2),2\sqrt{3}s)}{[2(1+3s^2)]^{2/3}},\quad
 \rho_0=\Big(\frac{4}{1+3s^2}\Big)^{2/3},\quad
 \rho_2=\Big(\frac{1+3s^2}{4}\Big)^{1/3}.
\end{equation}

\subsection*{Hopf bifurcations of the fixed points}

Using Cartesian coordinates $z=x+\ci y$ in \eqref{complexNFepsilon} we obtain
\begin{equation}\label{cartesianNFepsilon}
  \begin{pmatrix} \dot x \\ \dot y \end{pmatrix}=
  \begin{pmatrix} 1 \\ 0 \end{pmatrix}+
  \begin{pmatrix} \mu & -\nu \\ \nu & \mu \end{pmatrix}
  \begin{pmatrix} x \\ y \end{pmatrix}-
  (x^2+y^2)\begin{pmatrix} ax-by \\ bx+ay \end{pmatrix},
\end{equation}
where we have set $\epsilon=1$. The Jacobian of the right-hand side of
\eqref{cartesianNFepsilon} is given by
\begin{equation}\label{jacobian_epsilon}
  J=\begin{pmatrix} \mu-3ax^2-ay^2+2bxy & -\nu+bx^2+3by^2-2axy \\
             \nu-3bx^2-by^2-2axy & \mu-ax^2-3ay^2-2bxy \end{pmatrix}.
\end{equation}
The invariants of the Jacobian are given by
\begin{equation}
 T=2(\mu-2ar^2),\quad D=\mu^2+\nu^2-4(a\mu+b\nu)r^2+3r^4.
\end{equation}
A Hopf bifurcation takes place iff $T=0$ and $D>0$. When $T=0$,
$4a^2D=4(a\nu-b\mu)^2-\mu^2$. The fixed points satisfying $T=0$ are
given by $f(\rho)=0$ and $\mu=2a\rho$, resulting in the curve T in
parameter space
\begin{equation}
  4a\mu\nu(a\nu-b\mu)=8a^3-\mu^3,
\end{equation}
that can be parametrized as
\begin{equation}
  (\mu,\nu)=a^{1/3}(1-s^2)^{1/3}\bigg(2\,,
  \frac{b}{a}+\frac{s}{\sqrt{1-s^2}}\bigg),\quad s\in(-1,+1).
\end{equation}
For $s\to\pm1$, $\mu=0$ and $\nu\to\pm\infty$, and the curve is
asymptotic to the $\mu=0$ axis, the Hopf curve for $\epsilon=0$. Along
the T curve, the determinant $D$ is given by
\begin{equation}
  D=a^{2/3}(1-s^2)^{-1/3}\Big(2s^2-1-2\frac{b}{a}s\sqrt{1-s^2}\,\Big),
\end{equation}
resulting in two Hopf bifurcation curves (when $D>0$):
\begin{equation}
  {\rm H}_-:\ s\in\Big(-1,-\sqrt{(1-b)/2}\Big),\qquad
  {\rm H}_+:\ s\in\Big(\sqrt{(1+b)/2},+1\Big).
\end{equation}
The end points of these curves have $T=D=0$, and are Takens--Bogdanov
bifurcation points TB$_\pm$. They are precisely on the saddle-node curve
\eqref{epsSNcurve}, where both curves are tangent. The coordinates of
the four points Cusp$_\pm$ and TB$_\pm$ are:
\begin{align}
 &  (\mu,\nu)_{\text{Cusp}_+}=\frac{3}{2}\Big(a-\frac{b}{\sqrt{3}},
    b+\frac{a}{\sqrt{3}}\Big),
 && (\mu,\nu)_{\text{Cusp}_-}=\frac{3}{2}\Big(a+\frac{b}{\sqrt{3}},
    b-\frac{a}{\sqrt{3}}\Big),\\
 &  (\mu,\nu)_{\text{TB}_+}=\frac{(2a,2b+1)}{\big(2(1+b)\big)^{1/3}},\quad
 && (\mu,\nu)_{\text{TB}_-}=\frac{(2a,2b-1)}{\big(2(1-b)\big)^{1/3}}.
\end{align}
As TB$_\pm$ are on the saddle-node curve, its $\tilde s_\pm$ parameter,
according to \eqref{epsSNcurve}, can be computed. The result is 
$\tilde s_\delta=\delta\sqrt{(1-\delta b)/3(1+\delta b)}$, with
$\delta=\pm'$. Therefore TB$_+\in\,$SN$_0$, closer to Cusp$_+$ than to
Cusp$_-$; TB$_-\in\,$SN$_-$ if $\alpha_0<60^\text{o}$, TB$_-\in\,$SN$_0$ when
$\alpha_0>60^\text{o}$ and TB$_-=\,$Cusp$_-$ for $\alpha_0=60^\text{o}$.

\section{Symmetry breaking of $SO(2)$ to $Z_2$: computations}
\label{Appendix_bar_z}

\subsection*{Fixed points}

The fixed points of the normal form \eqref{complexNF} are
given by $\dot r=\dot\phi=0$. One solution is $r=0$ (named $P_0$). The
other fixed points are the solutions of
\begin{equation}
 \left.\begin{array}{l}
    \cos 2\phi=ar^2-\mu, \\ \sin 2\phi=\nu-br^2, \end{array}\right\}\quad
 \Rightarrow\quad (\mu-ar^2)^2+(\nu-br^2)^2=1,
\end{equation}
resulting in the bi-quadratic equation
$r^4-2(a\mu+b\nu)r^2+\mu^2+\nu^2-1=0$, whose solutions are
\begin{gather}
 r^2_\pm=a\mu+b\nu\pm\Delta,\quad
  \ce^{2\ci\phi_\pm}=(a\nu-b\mu\mp\ci\Delta)\ce^{\ci\alpha_0}\\
 \Delta^2=(a\mu+b\nu)^2+1-\mu^2-\nu^2=1-(a\nu-b\mu)^2,
\end{gather}
and every $\phi_\pm$ admits two solutions, differing by $\pi$ (they
are related by the symmetry $Z_2$, $z\to-z$ 
discussed above). Introducing a new phase $\alpha_1$,
\begin{equation}
 a\nu-b\mu-\ci\Delta=\ce^{\ci\alpha_1},
\end{equation}
where $\alpha_1\in[-\pi,0]$ because $\Delta>0$, we immediately obtain:
\begin{equation}
 \ce^{2\ci\phi_\pm}=\ce^{\ci(\alpha_0\pm\alpha_1)}\quad\Rightarrow\quad
 \phi_\pm=(\alpha_0\pm\alpha_1)/2,
\end{equation}
with the other solution being $(\alpha_0\pm\alpha_1)/2+\pi$;
$\alpha_1$ is a function of $(\mu,\nu)$ while $\alpha_0$ is a fixed
constant. We have obtained two pairs of $Z_2$ symmetric points,
$P_+=r_+\ce^{\ci\phi_+}$ and $P^*_+=-r_+\ce^{\ci\phi_+}$, and
$P_-=r_-\ce^{\ci\phi_-}$ and $P^*_-=-r_-\ce^{\ci\phi_-}$.

\subsection*{Hopf bifurcations of fixed points}

The Jacobian of the right-hand side of \eqref{cartesianNF} is given by
\begin{equation}\label{jacobian}
  J=\begin{pmatrix} \mu+1-3ax^2-ay^2+2bxy & -\nu+bx^2+3by^2-2axy \\
             \nu-3bx^2-by^2-2axy & \mu-1-ax^2-3ay^2-2bxy \end{pmatrix}.
\end{equation}
The invariants of the Jacobian are the trace $T$, the
determinant $D$ and the discriminant $Q=T^2-4D$. They are given by
\begin{align}
 & T=2(\mu-2ar^2), \\
 & D=\mu^2+\nu^2-1-4(a\mu+b\nu)r^2+3r^4+2\big(a(x^2-y^2)-2bxy\big), \\
 & Q=4\Big(1-\nu^2+4b\nu r^2+(1-4b^2)r^4-2\big(a(x^2-y^2)-2bxy\big)\Big).
\end{align}
The eigenvalues of the Jacobian matrix \eqref{jacobian} in terms of
the invariants are $\lambda_\pm=\frac{1}{2}(T\pm\sqrt{Q})$. For
example, a Hopf bifurcation takes place iff $T=0$ and $Q<0$. For the
fixed points $P_s$ and $P^*_s$, where $s=\pm$, we obtain
\begin{align}
 & T(P_s)=T(P^*_s)=2\big((b^2-a^2)\mu-2ab\nu-2as\Delta\big), \\
 & D(P_s)=D(P^*_s)=4s\Delta r^2_s, \\
 & Q(P_s)=Q(P^*_s)=4\big((\mu-2ar^2_s)^2-4s\Delta r^2_s\big).
\end{align}
As a result, $Q(P_-)=Q(P^*_-)>0$ and $P_-$ and $P^*_-$ never
experience a Hopf bifurcation. After some computations, $T(P_+)=0$
results in the ellipse
\begin{equation}
  \mu^2-4ab\mu\nu+4a^2\nu^2=4a^2,
\end{equation}
centered at the origin, contained between the straight lines SN$_+$ and
SN$_-$ and passing through the points $(\mu,\nu)=(0,\pm1)$, the ends of
the horizontal diameter of the circle $\mu^2+\nu^2=1$. This ellipse is
tangent to SN$_+$ and SN$_-$ at the points $(\mu,\nu)=\pm(2b,(b^2-a^2)/a)$.
The condition $Q<0$ is only satisfied on the elliptic arc from
$(\mu,\nu)=(0,-1)$ to $(2b,(b^2-a^2)/a)$ with $\mu>0$; along this arc
$P_+$ and $P^*_+$ undergo a Hopf bifurcation. We have assumed that $a$
and $b$ are both positive. The properties of the ellipse are:
\begin{align}
 & \text{major semiaxis}\quad (1-\ell_-)\mu=2ab\nu,
   && \text{length}\quad 2a/\sqrt{\ell_-}, \\
 & \text{minor semiaxis}\quad (1-\ell_-)\nu=-2ab\mu,
   && \text{length}\quad 2a/\sqrt{\ell_+},
\end{align}
where $2\ell_\pm=1+4a^2\pm\sqrt{1+8a^2}$. The eccentricity $e$ is given by
\begin{equation}\label{excentricity}
  \frac{2}{e^2}=1+\frac{1+4a^2}{\sqrt{1+8a^2}}.
\end{equation}

\subsection*{Codimension-two bifurcations of fixed points}

The Jacobian evaluated at the three points TB$_+$, TB$_-$ and TB is:
\begin{align}
 & J(\text{TB}_+)=\begin{pmatrix} 1 & -1 \\ 1 & -1 \end{pmatrix},\\
 & J(\text{TB}_-)=\begin{pmatrix} 1 & 1 \\ -1 & -1 \end{pmatrix},\\
 & J(\text{TB})=\begin{pmatrix} 1 & (1+b)/a \\
      (b-1)/a & -1 \end{pmatrix}.
\end{align}
The three matrices have double-zero eigenvalues and are of rank one,
so the three of them correspond to Takens--Bogdanov bifurcations. The
state that bifurcates at the TB point is $P_+$, without any symmetry,
so that it is an ordinary Takens--Bogdanov bifurcation, although the
$Z_2$ symmetric state $P^*_+$ also bifurcates at the same point in
parameter space (but removed in phase space) at another ordinary
Takens--Bogdanov bifurcation. The state that bifurcates at the
TB$_\pm$ points is $P_0$. This state is $Z_2$ symmetric, and so these
are Takens--Bogdanov bifurcations with $Z_2$ symmetry.

The Jacobian evaluated at the two points dPF$_\pm$ is
\begin{equation}
  J(\text{dPF}_+)=\begin{pmatrix} 1-b & -a \\ a & -1-b \end{pmatrix},\quad
  J(\text{dPF}_-)=\begin{pmatrix} b+1 & a \\ -a & b-1 \end{pmatrix}.
\end{equation}
The corresponding eigenvalues are $\lambda_+=-2b$ and $\lambda_-=0$
for dPF$_+$, and $\lambda_+=2b$ and $\lambda_-=0$ for dPF$_-$. Both
points are pitchfork bifurcations, and in order to determine if they
are degenerate, their normal form needs to be computed in order to
verify that the cubic term is zero. However, since in a degenerate
pitchfork bifurcation a curve of saddle-node bifurcations emerges that
is tangent to the pitchfork bifurcation curve, from
figure~\ref{bif_fixed_points} it is immediate apparent that both
dPF$_\pm$ are degenerate pitchfork bifurcations.

\section{Symmetry breaking of $SO(2)$ with quadratic terms: computations}
\label{Appendix_quadratic}

The three cases \eqref{realNFzbarz}, \eqref{realNFzz} and
\eqref{realNFbarzz} can be dealt with by considering the normal form
\begin{equation}\label{realNFquadratic}
\begin{aligned}
 & \dot r=r(\mu-ar^2)+r^2\cos m\phi,\\
 & \dot\phi=\nu-br^2+r\sin m\phi,
\end{aligned}
\end{equation}
where $m=1$ for the $\epsilon z^2$ case (\S\ref{Sec_z2}), $m=-1$ for
the $\epsilon z\bar z$ case (\S\ref{Sec_zbarz}) and $m=-3$ for the
$\epsilon\bar z^2$ case (\S\ref{Sec_Z3}). The fixed points, other than
the trivial solution $P_0$ ($r=0$), in the three cases are given by
the biquadratic equation $r^4-2(a\mu+b\nu+1/2)r^2+\mu^2+\nu^2=0$, with
solutions $P_\pm$
\begin{equation}
  r^2_\pm=a\mu+b\nu+1/2\pm\big(a\mu+b\nu+1/4-(a\nu-b\mu)^2\big)^{1/2}.
\end{equation}
The phases $\phi$ of the $P_\pm$ fixed points can be recovered from
\begin{equation}
 \cos m\phi=ar-\mu/r,\quad \sin m\phi=br-\nu/r.
\end{equation}
For $m=\pm1$ the solution is unique; for $m=3$ the solutions come in
triples, differing by $2\pi/m$. It is convenient to use the phase
space coordinates adapted to line L, introduced in \eqref{uv_munu}
(see also figure~\ref{uv_and_epsilon_SN}$(a)$). In terms of these coordinates,
$r^2_\pm=u+1/2\pm\sqrt{u+1/4-v^2}$, and the fixed points $P_\pm$ exist
only in the interior of the parabola $u=v^2-1/4$, whose axis is the
line L. On the parabola, these points are born in saddle-node
bifurcations. In order to explore additional bifurcations of these
points, we compute the Jacobian matrix of the normal form
\eqref{realNFquadratic},
\begin{equation}
 J=\begin{pmatrix} \mu-3ar^2+2r\cos m\phi & -mr^2\sin m\phi \\
                   -2br+\sin m\phi & mr\cos m\phi \end{pmatrix},
\end{equation}
whose trace and determinant, for the $P_\pm$ points, are easily computed:
\begin{equation}
\begin{aligned}
 & T(P_\pm)=(m-1)ar^2_\pm-(m+1)\mu, \\
 & D(P_\pm)=-2m\sqrt{u+1/4-v^2}\Big(\sqrt{(u+1/2)^2-u^2-v^2}\pm(u+1/2)\Big).
\end{aligned}
\end{equation}
Therefore $\sign D(P_\pm)=\mp\sign m$, and for $m>0$ ($m<0$) only
$P_+$ ($P_-$) may undergo a Hopf bifurcation.

\emph{The $\epsilon\bar z^2$ case (\S\ref{Sec_Z3}).} Here $m=-3$ and 
$T=2\mu-4ar^2$. $P_-$ is a saddle, but $P_+$ undergoes a Hopf bifurcation
when $T=0$. The condition $T=0$ for $P_+$ gives
\begin{equation}\label{ellipse_cond}
  \sqrt{a\mu+b\nu+1/4-(a\nu-b\mu)^2}=-\frac{1}{2a}
  \Big((a^2-b^2)\mu+2ab\mu+a\Big)>0.
\end{equation}
By squaring and simplifying, we obtain the ellipse
$(b\mu-2a\nu)^2+(a\mu-1)^2=1$ which is tangent to the line $\mu=0$ at
the origin, with its center at $(\mu,\nu)=(2a,b)/(2a^2)$, and whose
elements are:
\begin{align}
 & \text{major semiaxis parallel to }(1-\ell_-)\mu=2ab\nu,
   && \text{length}\quad 1/\sqrt{\ell_-}, \\
 & \text{minor semiaxis parallel to }(1-\ell_-)\nu=-2ab\mu,
   && \text{length}\quad 1/\sqrt{\ell_+},
\end{align}
where $2\ell_\pm=1+4a^2\pm\sqrt{1+8a^2}$. This ellipse has much in common
with the one found in the $\epsilon\bar z$ case, and the eccentricity
$e$ is given by the same expression \eqref{excentricity}.  For
$\alpha_0>\pi/6$, the ellipse is located in the interior of the
parabola of saddle-nodes, for $\alpha_0=\pi/6$ it becomes tangent to
the parabola at a single point, and for $\alpha_0<\pi/6$ it becomes
tangent at the two points
\begin{equation}
\begin{aligned}
 & \mu=\frac{1}{a}\big(1-2a^2-sb\sqrt{1-4a^2}\big),\\
 & \nu=\frac{1}{2a^2}\sqrt{1-4a^2}\big(b\sqrt{1-4a^2}-s(1-2a^2)\big),
   \quad s=\pm1.
\end{aligned}
\end{equation}
These are the points TB$_s$ in figure~\ref{quadratic_barzz}$(b)$. Only
the points on the elliptic arc H$_0$ joining these two points satisfy
\eqref{ellipse_cond}, and along this arc $P_+$ undergoes a Hopf
bifurcation.

\emph{The $\epsilon z\bar z$ case (\S\ref{Sec_zbarz}).} Here $m=-1$
and $T(P_\pm)=-2ar^2<0$, so there are no Hopf bifurcations. Moreover,
$D(P_+)>0$, so it is always stable and $D(P_-)<0$, so it is a
saddle. The only exception is when $r=0$, and this only happens at
$\mu=\nu=0$, the degenerate high-codimension point at the origin.

\emph{The $\epsilon z^2$ case (\S\ref{Sec_z2}).} Here $m=1$, and
$T=-2\mu$ is zero on the line $\mu=0$ inside the parabola. On this
line H$_0$, $P_+$ undergoes a Hopf bifurcation, and the points of
contact with the parabola have $D=T=0$ so they are Takens--Bogdanov
bifurcations (see figure~\ref{quadratic_bifs}b). The Hopf and
Takens--Bogdanov bifurcations are degenerate, as will be discussed in
\S\ref{degenerate_case}.

\subsection{A degenerate Takens-Bogdanov bifurcation}\label{degenerate_case}

In the $\epsilon z^2$ case, numerical simulations of the normal form
\eqref{complexNFzz} show that the Hopf bifurcation H$_0$ and the
Takens--Bogdanov points TB$_\pm$ are degenerate. This can also be
found by direct computation. Let us work out the details for the
TB$_-$ point.

The coordinates of TB$_-$ in parameter space are
$(\mu,\nu)=(0,-0.5/(1+b))=\big(0,1/4\cos^2(\alpha/2)\big)$, where $P_\pm$ are
born in a saddle-node bifurcation, and the fixed points are given by
\begin{equation}
  r^2_\pm=\frac{1}{2(1+b)}=\frac{1}{4\cos^2(\alpha/2)},\quad
  z_\pm=\frac{c+\ci}{2(1+b)}=\frac{\ci\ce^{-\ci\alpha/2}}{4\cos^2(\alpha/2)}.
\end{equation}
In order to obtain the normal form corresponding to the
Takens--Bogdanov point, a translation of the origin plus a convenient
rescaling of $z$ and time is made:
\begin{equation}
  t=4\tau\cos^2(\alpha/2),\ \zeta=2(z-z_\pm)\cos(\alpha/2),\
  \tilde\mu+\ci\tilde\nu=4(\mu+\ci\nu)\cos^2(\alpha/2).
\end{equation}
Substituting in \eqref{complexNFzz} results in 
\begin{equation}
  \dot\zeta=\ci\zeta+\ci\ce^{-2\ci\alpha}\bar\zeta+
  \ce^{\ci\alpha/2}\zeta^2+2\ce^{-3\ci\alpha/2}|\zeta|^2-
  \ci\ce^{-\ci\alpha}\zeta|\zeta|^2.
\end{equation}
In order to obtain the normal form, we introduce the real variables $(x_1,y_1)$
\begin{equation}
 \zeta=(y_1+2\ci x_1)\ce^{-\ci\alpha},
\end{equation}
so that the linear part of the ODE is transformed into Jordan form,
and we obtain
\begin{equation}\label{deg_TB}
\begin{aligned}
  \begin{pmatrix} \dot x_1 \\ \dot y_1 \end{pmatrix}= &
  \begin{pmatrix} y_1 \\ 0 \end{pmatrix} + \cos(\alpha/2)\begin{pmatrix}
    2x_1y_1 \\ 4x_1^2+3y_1^2 \end{pmatrix} + \sin(\alpha/2)\begin{pmatrix}
    -2x_1^2-3y_1^2/2 \\ 4x_1y_1\end{pmatrix} \\
  & -(4x^2+y^2)\begin{pmatrix} x\sin\alpha+(y/2)\cos\alpha \\
    y\sin\alpha-2x\cos\alpha \end{pmatrix} 
\end{aligned}
\end{equation}
Now we can reduce the quadratic and cubic terms to normal form by an
appropriate near-identity quadratic transformation
$(x_1,y_1)\to(x_2,y_2)$.  \citet{kno86} gives explicitely the normal
form coefficients up to and including third order, in terms of the
coefficients of the original ODE (in the form \ref{deg_TB}); a nice
summary is also given in \citet[][\S19.9]{Wig03}.  Using this explicit
transformation, we obtain
\begin{equation}\label{deg_TB_NF}
\left.\begin{array}{l}
  \dot x_2=y_2 \\ \dot y_2=4x_2^2\cos(\alpha/2)+16x_2^3\cos^2(\alpha/2)+O(4)
\end{array}\right\},
\end{equation}
and the $x_2y_2$ term in the normal form of the Takens--Bogdanov
bifurcation is missing, resulting in a degenerate case, the so-called
cusp case, of codimension three. The unfolding of this degenerate case
has been analyzed in detail in \citet{drs87}. Note that the ODE
\eqref{deg_TB_NF} is Hamiltonian at least up to order three, which
helps to explain the continuous family of periodic orbits obtained in
the interior of the homoclinic loop in
figure~\ref{quadratic_bifs}$(b)$.



\begin{thebibliography}{43}
\expandafter\ifx\csname natexlab\endcsname\relax\def\natexlab#1{#1}\fi

\bibitem[Abshagen {\em et~al.\/}(2008)Abshagen, Heise, Hoffmann \&
  Pfister]{AHHP08}
{\sc Abshagen, J., Heise, M., Hoffmann, C. \& Pfister, G.} 2008 Direction
  reversal of a rotating wave in {T}aylor-{C}ouette flow. {\em J.\,Fluid
  Mech.\/} {\bf 607}, 199--208.

\bibitem[Adler(1946)]{Adl46}
{\sc Adler, R.} 1946 A study of locking phenonena in oscillators. {\em Proc.
  Inst. Radio Engineers\/} {\bf 34}, 351--357.

\bibitem[Adler(1973)]{Adl73}
{\sc Adler, R.} 1973 A study of locking phenonena in oscillators. {\em Proc.
  IEEE\/} {\bf 61}, 1380--1385.

\bibitem[Arrowsmith \& Place(1990)]{arpl90}
{\sc Arrowsmith, D.~K. \& Place, C.~M.} 1990 {\em An Introduction to Dynamical
  Systems\/}. Cambridge University Press.

\bibitem[Broer {\em et~al.\/}(2008)Broer, van Dijk \& Vitolo]{BDV08}
{\sc Broer, H., van Dijk, R. \& Vitolo, R.} 2008 Survey of strong
  normal-internal k:l resonances in quasi-periodically driven oscillators for l
  = 1, 2, 3. In {\em SPT 2007: Symmetry and Perturbation Theory, Otranto
  2007\/} (ed. G.~Gaeta, R.~Vitolo \& S.~Walcher), pp. 45--55. World
  Scientific.

\bibitem[Campbell \& Holmes(1992)]{CaHo92}
{\sc Campbell, S.~A. \& Holmes, P.} 1992 Heteroclinic cycles and modulated
  travelling waves in a system with {$D_4$} symmetry. {\em Physica D\/} {\bf
  59}, 52--79.

\bibitem[Chossat \& Iooss(1994)]{ChIo94}
{\sc Chossat, P. \& Iooss, G.} 1994 {\em The {C}ouette--{T}aylor Problem\/}.
  Springer.

\bibitem[Chossat \& Lauterbach(2000)]{ChLa00}
{\sc Chossat, P. \& Lauterbach, R.} 2000 {\em Methods in Equivariant
  Bifurcations and Dynamical Systems\/}. World Scientific.

\bibitem[Chow {\em et~al.\/}(1994)Chow, Li \& Wang]{CLW94}
{\sc Chow, S.~N., Li, C. \& Wang, D.} 1994 {\em Normal Forms and Bifurcations
  of Planar Vector Fields\/}. Cambridge: Cambridge University Press.

\bibitem[Crawford \& Knobloch(1991)]{CrKn91}
{\sc Crawford, J.~D. \& Knobloch, E.} 1991 Symmetry and symmetry-breaking
  bifurcations in fluid dynamics. {\em Ann.\ Rev.\ Fluid Mech.\/} {\bf 23},
  341--387.

\bibitem[Dangelmayr {\em et~al.\/}(1997)Dangelmayr, Hettel \& Knobloch]{DHK97}
{\sc Dangelmayr, G., Hettel, J. \& Knobloch, E.} 1997 Parity-breaking
  bifurcation in inhomogeneous systems. {\em Nonlinearity\/} {\bf 10},
  1093--1114.

\bibitem[Dumortier {\em et~al.\/}(1987)Dumortier, Roussarie \&
  Sotomayor]{drs87}
{\sc Dumortier, F., Roussarie, R. \& Sotomayor, J.} 1987 Generic 3-parameter
  families of vector fields on the plane, unfolding a singularity with
  nilpotent linear part. {T}he cusp case of codimension 3. {\em Ergod.\ Th.\ \&
  Dynam.\ Sys\/} {\bf 7}, 375--413.

\bibitem[Dumortier {\em et~al.\/}(1997)Dumortier, Roussarie \&
  Sotomayor]{drs97}
{\sc Dumortier, F., Roussarie, R. \& Sotomayor, J.} 1997 Bifurcations of
  cuspidal loops. {\em Nonlinearity\/} {\bf 10}, 1369--1408.

\bibitem[Gambaudo(1985)]{Gam85}
{\sc Gambaudo, J.~M.} 1985 Perturbation of a {H}opf bifurcation by an external
  time-periodic forcing. {\em J. Diff. Eqns.\/} {\bf 57}, 172--199.

\bibitem[Golubitsky \& Schaeffer(1985)]{GoSc85}
{\sc Golubitsky, M. \& Schaeffer, D.~G.} 1985 {\em Singularities and Groups in
  Bifurcation Theory, vol. I\/}. Springer.

\bibitem[Golubitsky \& Stewart(2002)]{GoSt02}
{\sc Golubitsky, M. \& Stewart, I.} 2002 {\em The Symmetry Perspective: From
  Equilbrium to Chaos in Phase Space and Physical Space\/}. Birkh{\"{a}}user.

\bibitem[Golubitsky {\em et~al.\/}(1988)Golubitsky, Stewart \&
  Schaeffer]{GSS88}
{\sc Golubitsky, M., Stewart, I. \& Schaeffer, D.~G.} 1988 {\em Singularities
  and Groups in Bifurcation Theory\/}, {\em Appl.\ Math.\ Sci.\/}, vol.~II.
  Springer.

\bibitem[Haragus \& Iooss(2011)]{HaIo11}
{\sc Haragus, M. \& Iooss, G.} 2011 {\em Local Bifurcations, Center Manifolds,
  and Normal Forms in Infinite-Dimensional Dynamical Systems\/}. Springer.

\bibitem[Hirschberg \& Knobloch(1996)]{HiKn96}
{\sc Hirschberg, P. \& Knobloch, E.} 1996 Complex dynamics in the {H}opf
  bifurcation with broken translation symmetry. {\em Physica D\/} {\bf 90},
  56--78.

\bibitem[Keener(1987)]{Kee87}
{\sc Keener, J.~P.} 1987 Propagation and its failure in coupled systems of
  discrete excitable cells. {\em SIAM J. Appl. Math.\/} {\bf 47}, 556--572.

\bibitem[Knobloch(1986)]{kno86}
{\sc Knobloch, E.} 1986 Normal forms for bifurcations at a double-zero
  eigenvalue. {\em Phys.\ Lett.\,A\/} {\bf 115}, 199--201.

\bibitem[Knobloch {\em et~al.\/}(1995)Knobloch, Hettel \& Dangelmayr]{KHD95}
{\sc Knobloch, E., Hettel, J. \& Dangelmayr, G.} 1995 Parity breaking
  bifurcation in inhomogeneous systems. {\em Phys.\ Rev.\ Lett.\/} {\bf 74},
  4839--4842.

\bibitem[Kuznetsov(2004)]{Kuz04}
{\sc Kuznetsov, Y.~A.} 2004 {\em Elements of Applied Bifurcation Theory\/}, 3rd
  edn. Springer.

\bibitem[Lamb \& Wulff(2000)]{LaWu00}
{\sc Lamb, J. S.~W. \& Wulff, C.} 2000 Pinning and locking of discrete waves.
  {\em Physics Letters A\/} {\bf 267}, 167--173.

\bibitem[Lopez \& Marques(2009)]{LoMa09}
{\sc Lopez, J.~M. \& Marques, F.} 2009 Centrifugal effects in rotating
  convection: nonlinear dynamics. {\em J.\,Fluid Mech.\/} {\bf 628}, 269--297.

\bibitem[Marques \& Lopez(2006)]{MaLo06}
{\sc Marques, F. \& Lopez, J.~M.} 2006 Onset of three-dimensional unsteady
  states in small-aspect ratio {T}aylor-{C}ouette flow. {\em J.\,Fluid Mech.\/}
  {\bf 561}, 255--277.

\bibitem[Marques {\em et~al.\/}(2007)Marques, Mercader, Batiste \&
  Lopez]{MMBL07}
{\sc Marques, F., Mercader, I., Batiste, O. \& Lopez, J.~M.} 2007 Centrifugal
  effects in rotating convection: {A}xisymmetric states and three-dimensional
  instabilities,. {\em J.\,Fluid Mech.\/} {\bf 580}, 303--318.

\bibitem[Mercader {\em et~al.\/}(2006)Mercader, Batiste \& Alonso]{MBA06}
{\sc Mercader, I., Batiste, O. \& Alonso, A.} 2006 Continuation of
  travelling-wave solutions of the {Navier-Stokes} equations. {\em Intnl
  J.\,Num.\ Meth.\ Fluids\/} {\bf 52}, 707--721.

\bibitem[Mercader {\em et~al.\/}(2010)Mercader, Batiste \& Alonso]{MBA10}
{\sc Mercader, I., Batiste, O. \& Alonso, A.} 2010 An efficient spectral code
  for incompressible flows in cylindrical geometries. {\em Computers and
  Fluids\/} {\bf 39}, 215--224.

\bibitem[Pacheco {\em et~al.\/}(2011)Pacheco, Lopez \& Marques]{PLM11}
{\sc Pacheco, J.~R., Lopez, J.~M. \& Marques, F.} 2011 Pinning of rotating
  waves to defects in finite {Taylor-Couette} flow. {\em J.\,Fluid Mech.\/}
  {\bf 666}, 254--272.

\bibitem[Pfister {\em et~al.\/}(1992)Pfister, Buzug \& Enge]{PBE92}
{\sc Pfister, G., Buzug, T. \& Enge, N.} 1992 Characterization of experimental
  time series from {T}aylor-{C}ouette flow. {\em Physica D\/} {\bf 58},
  441--454.

\bibitem[Pfister {\em et~al.\/}(1988)Pfister, Schmidt, Cliffe \&
  Mullin]{PSCM88}
{\sc Pfister, G., Schmidt, H., Cliffe, K.~A. \& Mullin, T.} 1988 Bifurcation
  phenomena in {T}aylor-{C}ouette flow in a very short annulus. {\em J.\,Fluid
  Mech.\/} {\bf 191}, 1--18.

\bibitem[Pfister {\em et~al.\/}(1991)Pfister, Schulz \& Lensch]{PSL91}
{\sc Pfister, G., Schulz, A. \& Lensch, B.} 1991 Bifurcations and a route to
  chaos of an one-vortex-state in {T}aylor-{C}ouette flow. {\em Eur. J. Mech.
  B-Fluids\/} {\bf 10}, 247--252.

\bibitem[Saleh \& Wagener(2010)]{SaWa10}
{\sc Saleh, K. \& Wagener, F. O.~O.} 2010 Semi-global analysis of periodic and
  quasi-periodic normal-internal $k:1$ and $k:2$ resonances. {\em
  Nonlinearity\/} {\bf 23}, 2219--2252.

\bibitem[Sanchez {\em et~al.\/}(2002)Sanchez, Marques \& Lopez]{SML02}
{\sc Sanchez, J., Marques, F. \& Lopez, J.~M.} 2002 A continuation and
  bifurcation technique for {N}avier-{S}tokes flows. {\em J.\,Comput.\ Phys.\/}
  {\bf 180}, 78--98.

\bibitem[Shil'nikov {\em et~al.\/}(2001)Shil'nikov, Shil'nikov, Turaev \&
  Chua]{SSTC01}
{\sc Shil'nikov, L.~P., Shil'nikov, A.~L., Turaev, D.~V. \& Chua, L.~O.} 2001
  {\em Methods of Qualitative Theory in Nonlinear Dynamics. Part II.\/}. World
  Scientific.

\bibitem[Strogatz(1994)]{Str94}
{\sc Strogatz, S.} 1994 {\em Nonlinear Dynamics and Chaos\/}. Addison-Wesley.

\bibitem[Thiele \& Knobloch(2006{\natexlab{{\em a\/}}})]{ThKn06_prl}
{\sc Thiele, U. \& Knobloch, E.} 2006{\natexlab{{\em a\/}}} Driven drops on
  heterogeneous substrates: onset of sliding motion. {\em Phys.\ Rev.\ Lett.\/}
  {\bf 97}, 204501.

\bibitem[Thiele \& Knobloch(2006{\natexlab{{\em b\/}}})]{ThKn06_njp}
{\sc Thiele, U. \& Knobloch, E.} 2006{\natexlab{{\em b\/}}} On the depinning of
  a driven drop on a heterogeneous substrate. {\em New Journal of Physics\/}
  {\bf 8}, 313.

\bibitem[Turaev \& Shilnikov(1995)]{TuSh95}
{\sc Turaev, D.~V. \& Shilnikov, L.~P.} 1995 On blue sky catastrophies. {\em
  Doklady Akademii Nauk\/} {\bf 342}, 596--599.

\bibitem[Wagener(2001)]{Wag01}
{\sc Wagener, F.} 2001 Semi-local analysis of the k : 1 and k : 2 resonances in
  quasi-periodically forced systems. In {\em Global Analysis of Dynamical
  Systems\/} (ed. W.~Broer, B.~Krauskopf \& G.~Vegter), pp. 113--129. IOP
  Publishing Ltd.

\bibitem[Westerburg \& Busse(2003)]{WeBu03}
{\sc Westerburg, M. \& Busse, F.~H.} 2003 Centrifugally driven convection in
  the rotating cylindrical annulus with modulated boundaries. {\em Nonlin.
  Proc. Geophys.\/} {\bf 10}, 275--280.

\bibitem[Wiggins(2003)]{Wig03}
{\sc Wiggins, S.} 2003 {\em Introduction to Applied Nonlinear Dynamical Systems
  and Chaos\/}, 2nd edn. Springer-Verlag.

\end{thebibliography}
\end{document}